%% file: pilatusrev.tex
\documentclass[11pt]{article}
\input{macros_angl_pdf.tex}

\renewcommand{\d}{{\mathfrak m}}
\newcommand{\dd}{{\mathfrak n}}
\renewcommand{\D}{\mathbf{d}}
\newcommand{\val}{\nu}

\title{Joint effective equidistribution of partial lattices \\ 
in positive characteristic}
\author{Tal Horesh \and Fr\'ed\'eric Paulin} 
\begin{document}
\bibliographystyle{../alphanum}
\maketitle
\begin{abstract}
Let $\val$ be a place of a global function field $K$ over a finite
field, with associated affine function ring $R_\val$ and completion
$K_\val$, and let $1\leq \d<\D$. The aim of this paper is to prove an
effective triple joint equidistribution result for primitive partial
$R_\val$-lattices $\Lambda$ of rank $\d$ in $K_\val^{\;\D}$ as their
covolume tends to infinity: of their $K_\val$-linear span $V_\Lambda$
in the rank-$\d$ Grassmannian space of $K_\val^{\;\D}$; of their shape
in the modular quotient by $\PGL_\d(R_\val)$ of the Bruhat-Tits
buildings of $\PGL_\d(K_\val)$; and of the shape of $\Lambda^\perp$ in
the similar quotient for $\PGL_{\D-\d}(K_\val)$, where $\Lambda^\perp$
is the orthogonal partial $R_\val$-lattice of rank $\D-\d$ in the dual
space of $K_\val^{\;\D}$. The main tools are a new refined LU
decomposition by blocks of elements of $\SL_\D(K_\val)$, techniques of
Gorodnik and Nevo for counting integral points in well-rounded
families of subsets of algebraic groups, and computations of volumes
of various homogeneous spaces associated with partial
$R_\val$-lattices.
\end{abstract}

\section{Introduction}
\label{sect:intro}

We fix throughout the paper three positive integers $\d$, $\dd$, $\D$
such that $\D=\d+\dd$. A {\it primitive integral vector} in $\RR^\D$
is an element of $\ZZ^\D$ with coprime componants, so that the cyclic
group it generates is a free abelian factor of rank $1$ of
$\ZZ^\D$. More generally, a {\it primitive $\d$-lattice} in $\RR^\D$
is a free abelian factor of rank $\d$ in $\ZZ^\D$. The distribution
problems of primitive integral vectors and of primitive $\d$-lattices
have first been studied by Linnik and Maass (see for instance
\cite{Linnik68} and \cite{Mass59}), and have given rise to a huge
amount of works using various tools, see for instance \cite{Schmidt68,
  Schmidt98, Duke03, Duke07, Marklof10, EllMicVen13, EinMozShaSha16,
  HorKar19}.

Let us define the {\it covolume} $\covol(\Lambda)$ of a primitive
$\d$-lattice $\Lambda$ in $\RR^\D$ as the Lebesgue volume of the
parallelepided generated by any $\ZZ$-basis of $\Lambda$. Let us
define the {\it shape} $\sh(\Lambda)$ of $\Lambda$ as its equivalence
class modulo rotations and homotheties in $\RR^\D$ (or its
``similarity class'' with the terminology of \cite{Schmidt98}), which
belongs to the double coset space $\operatorname{PSh}_\d=
\operatorname{PSO}(\d)\bs\PGL_\d(\RR)/ \PGL_\d(\ZZ)$. Let us denote by
$V_\Lambda$ the $\RR$-linear subspace of $\RR^\D$ generated by
$\Lambda$, which belongs to the Grassmannian space $\Gr_{\d,\D}$ of
$\d$-dimensional $\RR$-linear subspaces of $\RR^\D$. Schmidt in
\cite{Schmidt15} proved that the pairs $(V_\Lambda,\sh(\Lambda))$
equidistribute in $\Gr_{\d,\D}\times \operatorname{PSh}_\d$ (towards
the natural product measure) as the covolumes of the primitive
$\d$-lattices $\Lambda$ tend to infinity in average.  Without
averaging, stronger equidistribution results have been obtained when
fixing the covolume of the primitive $\d$-lattices and letting it go
to $\infty$ (possibly requiring some congruence properties), see for
instance \cite{AkaEinSha16JLMS, AkaEinSha16Inv, EinRuhWir19} when
$\d=1$, \cite{AkaEinWie22} when $\d=\dd=2$, and \cite{BerSha23}.
Schmidt's result in average was strengthen (with an effective version)
in \cite{HorKar23}, by adding in the third factor of
$\Gr_{\d,\D}\times \operatorname{PSh}_\d\times
\operatorname{PSh}_{\dd}$ the equidistribution of the shapes of the
orthogonal $\dd$-lattices $\Lambda^\perp=V_\Lambda^\perp\cap \ZZ^\D$
of the primitive $\d$-lattices $\Lambda$.  Removing the averaging
aspect and under appropriate congruence conditions, this triple
equidistribution result has been extended in \cite{AkaMusWie24}.

In this paper, replacing $\QQ$ by any global function field over a
finite field $\FF_q$ (not only the field $\FF_q(Y)$ of rational
fractions, though already new and interesting), and $\ZZ$ by the
corresponding affine ring for any fixed choice of place at infinity
(not only the polynomial ring $\FF_q[Y]$), we give the first complete
treatment of the triple equidistribution of primitive partial lattices
in positive characteristic. Our results are effective and allow
versions with congruences. Since this affine ring is no longer
principal in general, we will need to adapt our tools. Our result is
not a result in average, we will also fix the covolume of the
primitive $\d$-lattices and let it go to $\infty$ (also requiring an
appropriate congruence properties).
%%
%\todo{\tiny to be seen}
%%
This is a major extension of the case $\D=2$ considered in
\cite{HorPau22}, since only a double equidistribution makes sense in
this dimension, and many moduli spaces constructions and volume
computations were not available.
%We refer to \cite{HorPau24} for the consequences of this work purely in terms
%of equidistribution of rational points in the Grassmannian spaces (and
%in particular to Frobenius numbers).
Class number issues prevent an exact correspondence between primitive
partial lattices and rational points in Grassmannian spaces, that was
satisfied in the real field case, to take place. An analogous (though
different) study of the distribution of the rational points just in
the Grassmannian spaces over some function fields had been conducted
in \cite{Thunder08}.

\medskip
More precisely, referring to \cite{Goss98,Rosen02} and Subsection
\ref{subsect:functionfields} for definitions and complements, we fix a
global function field $K$ of genus $\ggg$ over a finite field $\FF_q$
of order $q$, a (discrete normalized) valuation $\val$ of $K$ and a
uniformizer $\pi_\val$ of $\val$. We denote by $\zeta_K$ the Dedekind
zeta function of $K$, by $K_\val$ the completion of $K$ associated
with $\val$, by $\OOO_\val$ its valuation ring, by $q_\val$ the order
of its residual field, by $|\cdot|= q_\val^{\;-\val(\,\cdot\,)}$ its
(normalized) absolute value, and by $R_\val$ the affine function ring
associated with $\val$ (for instance, $K= \FF_q(Y)$, $\ggg=0$,
$\val(P/Q)=\deg Q-\deg P$ for all $P,Q\in \FF_q[Y]$, $\OOO_\val=
\FF_q[[Y^{-1}]]$, $q_\val=q$ and $R_\val= \FF_q[Y]$).

We endow $K_\val^{\;\D}$ with the supremum norm, and any
$K_\val$-linear subspace of $K_\val^{\;\D}$ with its induced norm and
with its associated normalized Haar measure (giving mass one to its
closed unit ball), see Subsection \ref{subsect:partiallatt}.  A {\it
  partial $R_\val$-lattice $\Lambda$ of rank $\d$} (or $\d$-lattice
for short) in $K_\val^{\;\D}$ is a discrete free $R_\val$-submodule of
rank $\d$ generating an $\d$-dimensional $K_\val$-vector subspace
$V_\Lambda$ of $K_\val^{\;\D}$. We denote by $\covol(\Lambda)$ the
covolume of $\Lambda$ in $V_\Lambda$. We say that $\Lambda$ is {\it
  primitive} if it is a free direct factor of $R_\val^{\;\D}$. Among
all the definitions of primitiveness that were equivalent for partial
$\ZZ$-lattices in $\RR^\D$ and non longer are, this turns out to be the
appropriate one. We denote by $\PL_{\d,\D}$ the set of primitive
$\d$-lattices in $K_\val^{\;\D}$. For every point $y$ of every
measurable space, we denote by $\Delta_y$ the unit Dirac mass at $y$.

\medskip
Our first result is a joint equidistribution result in modular
quotients of Bruhat-Tits buildings. For $k\in\{\d,\dd\}$, let
$\I_{\val,k}$ be the Bruhat-Tits building of the simple algebraic
group $\PGL_k$ over the local field $K_\val$ (see for instance
\cite{BruTit84}).  Its $\PGL_k(K_\val)$-homogeneous set of vertices
$V\!\!\I_{\val,k}$ is the (discrete) set of
$K_\val^{\;\times}$-homothety classes $[L]$ of $\OOO_\val$-lattices
$L$ of $K_\val^{\;k}$. The unimodular group $\GL^1_k(K_\val) = \{g\in
\GL_k(K_\val):|\det g|=1\}$ acts projectively on $V\!\!\I_{\val,k}$
with finitely many orbits. We denote by $V_0\I_{\val,k}$ the orbit by
$\GL^1_k(K_\val)$ of the vertex $[\OOO_\val^{\;k}]\in
V\!\!\I_{\val,k}$, identified with $\GL_k^1(K_\val)/\,
\GL_k(\OOO_\val)$.  The action of the modular group $\wt\Ga_k=
\GL_k(R_\val)$ on $V\!\!\I_{\val,k}$ is proper (with finite
stabilisers $\wt\Ga_{k,x}$ of every vertex $x\in V\!\!\I_{\val,k}$)
with (discrete) infinite quotient $\wt\Ga_k\bs V\!\!\I_{\val,k}$. See
for instance \cite{Serre83} \cite[\S 15.2]{BroParPau19} when $k=2$ and
\cite{Soule79} when $\ggg=0$ for the structure of the quotient complex
of groups (in the sense of \cite{BriHae99}) $\wt\Ga_k\dbs
\I_{\val,k}$.  The measure on $\wt\Ga_k\bs V_0\I_{\val,k}$ induced by
the counting measure on $V_0\I_{\val,k}$ is the finite measure
\[
\mu_{\wt\Ga_k\bs V_0 \I_{\val,k}} =\sum_{[x]\in\wt\Ga_k\bs V_0\I_{\val,k}}
\frac{1}{\card\;\wt\Ga_{k,x}} \Delta_{[x]}
\]
(see \cite[\S 1.5]{BasLub01} when $k=2$).  We identify the quotient
space $\wt\Ga_k\bs V_0\I_{\val,k}$ with the (discrete) double coset space
$\Sh_{k}= \GL_k(\OOO_\val) \,\bs \GL_k^1(K_\val) \,/\,\wt\Ga_k$ by the
map induced by $g\mapsto g^{-1}$ on the double cosets $D\in\Sh_{k}$, that
we still denote by $D\mapsto D^{-1}$.

Let $\Lambda\in \PL_{\d,\D}$ be such that there exists $i\in\ZZ$ with
$\frac{\covol(\Lambda)}{\covol (R_\val^{\;\d})} =
q_\val^{\;\lcm(\d,\dd)\, i}$ and let $g\in\GL_\D (\OOO_\val)$ sending
$V_\Lambda$ to $K_\val^{\;\d}\times\{0\}$ (so that $g\Lambda$ becomes
a full $R_\val$-lattice of $K_\val^{\;\d}$). We define (see Equation
\eqref{eq:defish}) the {\it shape} of $\Lambda$ as the class of the
$\d$-lattice $\Lambda$ modulo scaling and modulo the action of the
maximal compact subgroup $\GL_\D (\OOO_\val)$, that is,
\[
\sh(\Lambda)=\GL_\d(\OOO_\val)\,
\pi_\val^{\;-\frac{\lcm(\d,\dd)} {\d}\,i}\, g\,\Lambda\in \Sh_{\d}\;.
\]
Contrarily to the real field case, the rescaling process is much
harder when the absolute value is discrete, hence the above
restriction on the covolumes.
%%
%\todo{\tiny See Section 4.4 for a way to go around
%  this restriction. ?}
%%
Let ${V_\Lambda}^\perp$ be the subspace of the dual space
$(K_\val^{\;\D})^*$ of $K_\val^{\;\D}$ consisting in the
$K_\val$-linear forms on $K_\val^{\;\D}$ vanishing on $V_\Lambda$ and
let $R_\val^{\;\D,*}$ be the standard $R_\val$-lattice in
$(K_\val^{\;\D})^*$ generated as an $R_\val$-module by the dual basis
of the canonical basis of $K_\val^{\;\D}$. We define the {\it
  orthogonal $R_\val$-lattice} $\Lambda^\perp$ of the primitive
$\d$-lattice $\Lambda$ by $\Lambda^\perp = {V_\Lambda}^\perp\cap
R_\val^{\;\D,*}$, which is a primitive $\dd$-lattice in
$(K_\val^{\;\D})^*$, see Subsection \ref{subsect:ortholatt}. The
following joint equidistribution result of the pairs of shapes of
primitive $\d$-lattices and their orthogonal $\dd$-lattices as their
covolume tends to infinity, in the product of the quotients of the
Bruhat-Tits buildings $\I_{\val,\d}$ and $\I_{\val,\dd}$ by their
modular groups $\wt\Ga_\d$ and $\wt\Ga_\dd$, is a corollary of Theorem
\ref{theo:mainintrotriple}. See the end of Subsection
\ref{subsect:proofmain} for its proof.

Let us define an arithmetic constant $ c'=
(q-1)\,q^{\;(\ggg-1)(\D^2-1-\d\,\dd)}
\prod_{i=1}^{\D-1}\frac{\zeta_K(i+1)}{q_\val^{\;i}-1}$.

\bcoro\label{coro:intro2} For the weak-star convergence of Borel
measures on the locally compact space $\wt\Ga_\d\bs V_0
\I_{\val,\d}\times \wt\Ga_\dd\bs V_0\I_{\val,\dd}$, we have
%%
%\todo{\tiny try to state
%  result purely in terms of Bruthat-Tits building}
%%
\begin{align*}
\lim_{i\ra+\infty}\;\frac{c'}{q_\val^{\;\lcm(\d,\dd) \,\D\,i}}
\sum_{\Lambda\in \PL_{\d,\D}\;:\;\frac{\covol(\Lambda)}{\covol(R_\val^{\;\d})}=
  q_\val^{\;\lcm(\d,\dd) \,i}}& \Delta_{\sh(\Lambda)^{-1}} \otimes
\Delta_{\sh(\Lambda^\perp)^{-1}} \\&= \mu_{\wt\Ga_\d\bs V_0 \I_{\val,\d}}
\otimes \mu_{\wt\Ga_\dd\bs V_0 \I_{\val,\dd}} \;.
\end{align*}
\ecoro

We refer to Equation \eqref{eq:equidisbuilding} for error terms and
for versions with congruences of this corollary.
%%
%\todo{\tiny , and Section 4.4 for a version without
%the restriction on the covolumes. ?}
%%.
The main result of this paper is the following triple joint
equidistribution theorem. We endow the unimodular group
$\GL^1_\d(K_\val)$ with its Haar measure giving mass $1$ to its
maximal compact subgroup $\GL_\d(\OOO_\val)$ and the (discrete,
infinite) double quotient $\Sh_{\d}$ with its induced measure
$\mu_{\Sh_{\d}}$, which is finite (see Subsection
\ref{subsect:shape}). We denote by $\Gr_{\d,\D}$ the compact
Grassmannian space of $\d$-dimensional $K_\val$-linear subspaces of
$K_\val^{\;\D}$, and by $\mu_{\Gr_{\d,\D}}$ its
$\GL_\D(\OOO_\val)$-invariant probability measure.

\btheo\label{theo:mainintrotriple} For the weak-star convergence
of Borel measures on the locally compact space $\Gr_{\d,\D}\times
\Sh_{\d}\times \Sh_{\dd}$, we have
\begin{align}
\lim_{i\ra+\infty}\;\frac{c'}{q_\val^{\;\lcm\{\d,\dd\}\,\D\,i}}
\sum_{\Lambda\in \PL_{\d,\D}\;:\;\frac{\covol(\Lambda)}{\covol(R_\val^{\;\d})}=
  q_\val^{\;\lcm\{\d,\dd\}\, i}}&\Delta_{V_\Lambda}\otimes\Delta_{\sh(\Lambda)}
\otimes\Delta_{\sh(\Lambda^\perp)}\nonumber\\&=
\mu_{\Gr_{\d,\D}}\otimes\mu_{\Sh_{\d}}\otimes\mu_{\Sh_{\dd}}\;.
\label{eq:mainintrotriple}
\end{align}
\etheo

We refer to Corollary \ref{coro:mainshape} for error terms and for
versions with congruences of this theorem.
%%
%\todo{\tiny , and Section 4.4 for a version without
%  the restriction on the covolumes. ?}
%%.
We will actually prove in Theorem \ref{theo:mainsharp} a much stronger
(albeit more technical) equidistribution result. We will define in
Subsection \ref{subsect:latticepair} a (non discrete) moduli space
$\La_{\d,\dd}$ of pairs $(L,L')\in (\GL^1_\d(K_\val)/\GL_\d(R_\val))
\times (\GL^1_\dd(K_\val)/\GL_\dd(R_\val))$ of unimodular $\d$- and
$\dd$-lattices with an appropriate correlation on the determinant of
any of their $R_\val$-basis. Using a refined LU decomposition of
elements of $\SL_\D(K_\val)$ introduced in Subsection
\ref{sect:LUbloc}, to any primitive $\d$-lattice $\Lambda$, we will
associate such a pair $\llbracket \Lambda \rrbracket\in\La_{\d,\dd}$,
under a restriction that the linear subspace $V_\Lambda$ of the
primitive $\d$-lattice $\Lambda$ belongs to the unit ball
$\Gr_{\d,\D}^\flat$ of the lower maximal Bruhat cell of the
Grassmannian space $\Gr_{\d,\D}$ (see Subsection
\ref{subsect:grassmann} for more details). We will then prove in
Theorem \ref{theo:mainsharp} the equidistribution of the pairs
$(V_\Lambda,\llbracket \Lambda\rrbracket)$ in
$\Gr_{\d,\D}^\flat\times\La_{\d,\dd}$ for the primitive $\d$-lattices
$\Lambda$ whose covolume is fixed (satisfying some congruence
property) and tends to $+\infty$. Theorem \ref{theo:mainintrotriple}
will follow by a tricky consideration of compound matrices.

We refer to Subsection \ref{subsect:congruence} for the description of
the appropriate congruence subgroups of $\SL_\D(K_\val)$ that we will
use for the version with congruences of our theorems. A major part of
the paper consists of a fine study of the homogeneous measures on the
various homogeneous spaces $\Gr_{\d,\D}$ (see Subsection
\ref{subsect:grassmann}), $\La_{\d,\dd}$ (see Subsections
\ref{subsect:fulllat} and \ref{subsect:latticepair}), and on the
double coset spaces $\Sh_k$ for $k\in\{\d,\dd\}$ (see Subsection
\ref{subsect:shape}), besides the precise disintegration of the Haar
measure of $\SL_\D(K_\val)$ under the refined LU decomposition in
Subsection \ref{sect:LUbloc}. A key tool of this paper is the counting
result in well-rounded sets of integral points of algebraic groups
over $K$ by Gorodnik and Nevo \cite{GorNev12}. A long study is
necessary in order to introduce the appropriate well-rounded sets, to
prove that they are indeed well-rounded, and to compute their
measures: see Subsection \ref{subsect:correspondance} which gives a
precise relationship between primitive $\d$-lattices in
$K_\val^{\;\D}$ and integral matrices in $\SL_\D(R_\val)$, and
Subsection \ref{subsect:wellround}.

%
%\todo{\tiny Given $\bm{k}=(k_1,\dots, k_\ell)$, generalisation to
%  equidistribution of primitive ``graded'' $\bm{k}$-lattices in other
%  flag spaces that Grassmanians?}

\noindent\mbox{}
%
%\todo{\tiny Analog for $S$-arithmetic groups in characteristic $0$?}
%

\medskip
\noindent{\small {\it Acknowledgements: } The first author thanks the
  Laboratoire de mathématique d'Orsay for visiting financial support, and the
  second author thanks the ETH for visiting financial support.}

\section{Background definitions and notation}
\label{sect:background}

\subsection{On global function fields}
\label{subsect:functionfields}
We refer for instance to \cite{Goss98,Rosen02} and
\cite[Chap.~14]{BroParPau19} for the content of this Section.

Let $\FF_q$ be a finite field of order $q$, where $q$ is a positive
power of a positive prime. Let $K$ be a (global) {\it function field}
over $\FF_q$, that is, the function field of a geometrically connected
smooth projective curve ${\bf C}$ over $\FF_q$, or equivalently an
extension of $\FF_q$ of transcendence degree $1$, in which $\FF_q$ is
algebraically closed. We denote by $\ggg$ the genus of the curve ${\bf
  C}$.

There is a bijection between the set of closed points of ${\bf C}$ and
the set of (normalized discrete) valuations $\val$ of its function
field $K$, where the valuation of a given element $f\in K$ is the order
of the zero or the opposite of the order of the pole of $f$ at the
given closed point. We fix such a valuation $\val$ from now on.

We denote by $K_\val$ the completion of $K$ for the valuation
$\val$, and by 
\[
\OOO_\val= \{x\in K_\val\;:\;\val(x)\geq 0\}
\]
the valuation ring of (the unique extension to $K_\val$) of $\val$.
Let us fix a uniformiser $\pi_\val\in K$ of $\val$, that is, an
element in $K$ with $\val(\pi_\val)=1$.  We denote by $q_\val$ the
order of the residual field $\OOO_\val/\pi_\val\OOO_\val$ of $\val$,
which is a (possibly proper) power of $q$.  We normalize the absolute
value associated with $\val\,$ as usual: for every $x\in K_\val$, we
have the equality
\[
|\,x\,|=q_\val^{\;-\val(x)}\;.
\]

Finally, let $R_\val$ denote the affine algebra of the affine curve
${\bf C} \ssm\{\val\}$, consisting of the elements of $K$ whose only
poles (if any) are at the closed point $\val$ of ${\bf C}$. It is a
Dedeking ring and its field of fractions is equal to $K$. Note that
(see for instance \cite[Eq.~(14.2) and (14.3)]{BroParPau19}
\begin{equation}\label{eq:inversiRv}
  R_\val\cap\OOO_\val=\FF_q\quad\text{and}\quad
  R_\val^{\;\times}=\FF_q^{\,\times}\subset \OOO_\val^{\;\times} \;.
\end{equation}
The {\it Dedekind zeta function} of $K$ (see for instance \cite[\S
  5]{Rosen02}) is defined if $\Re\; s > 1$ by
\[
\zeta_K(s) = \sum_I \frac{1}{N(I)^s}\;,
\]
where the summation is over the nonzero ideals $I$ of $R_\val$, with
norm $N(I)=[R_\val:I]$. By \cite[Theo.~5.9]{Rosen02}), it is a
rational function of $q^{-s}$ and has an analytic continuation on
$\CC\ssm\{0, 1\}$ with simple poles at $s = 0$, $s = 1$. Furthermore,
it has positive values at $s=-i$ for all $i\in\NN\ssm\{0\}$, since by
the functional equation of $\zeta_K$ (see loc.~cit.), we have
\begin{equation}\label{eq:poszetanegent}
\zeta_K(-i)=q^{(\ggg-1)(1+2\,i)}\zeta_K(1+i)>0\;.
\end{equation}

\medskip
The simplest example corresponds to ${\bf C}=\PP^1$ (so that
$\ggg=0$) and $\val=\val_\infty$ the valuation associated with
the point at infinity $[1:0]$. Then

$\bullet$~ $K=\FF_q(Y)$ is the field of rational functions in one
variable $Y$ over $\FF_q$, 

$\bullet$~ $\val_\infty$ is the valuation defined, for all $P,Q\in
\FF_q[Y]$, by
\[
\val_\infty(P/Q)=\deg Q-\deg P\;.
\]

$\bullet$~ $R_{\val_\infty}=\FF_q[Y]$ is the (principal) ring of
polynomials in one variable $Y$ over $\FF_q$,

$\bullet$~ $K_{\val_\infty}= \FF_q((Y^{-1}))$ is the field of formal
Laurent series in one variable $Y^{-1}$ over $\FF_q$,

$\bullet$~ $\OOO_{\val_\infty}= \FF_q[[Y^{-1}]]$ is the ring of formal
power series in one variable $Y^{-1}$ over $\FF_q$, $\pi_{\val_\infty}
=Y^{-1}$ is the usual choice of a uniformizer, and $q_{\val_\infty}=q$.

\subsection{Partial lattices}
\label{subsect:partiallatt}

Let $V$ be a $K_\val$-vector space with finite dimension $D\geq 1$
endowed with an ultrametric norm $\|\;\|$, and let $k\in\llbracket
1,D\rrbracket$. We denote by $B_V(0,1)$ the closed unit ball of
$V$. We endow $V$ with the unique Haar measure $\mu_V$ of the abelian
locally compact topological group $(V,+)$ such that
$\mu_V(B_V(0,1))=1$.  This measure scales as follows under linear
maps: for all $x\in V$ and $g\in\GL(V)$, we have
\begin{equation}\label{eq:homothetyhaar}
d\mu_{V}(g x)=|\det g|\;d\mu_{V}(x)\;.
\end{equation}
When $V$ is $K_\val^{\;\D}$ with its canonical basis $(e_1,\dots,
e_\D)$, we will take the supremum norm $\|\lambda_1e_1+\dots
+\lambda_\D e_\D\| ={\displaystyle \max_{1\leq i\leq \D}}\,
|\lambda_i|$. The Haar measure of $K_\val^{\;\D}$ is then normalized
so that $\mu_{K_\val^{\;\D}} (\OOO_\val^{\;\D}) =1$. On the dual
$K_\val$-vector space $V^*$, we will consider the dual norm $\|\;\|^*$
(which is ultrametric). When $V$ is $K_\val^{\;\D}$, the dual norm on
$V^*$ is the supremum norm with respect to the dual basis
$(e_1^*,\dots, e_\D^*)$ of $(e_1,\dots, e_\D)$.

Recall that for every $g\in \GL(V)$, its (left) action $\widecheck
g:\ell\mapsto \ell\circ g^{-1}$ on the dual space $V^*$ satisfies, for
every $K_\val$-basis $\B$ of $V$ with dual $K_\val$-basis $\B^*$ of
$V^*$, that
\begin{equation}\label{eq:duallinmap}
\operatorname{Mat}_{\B^*}(\widecheck g)=
\;^t\operatorname{Mat}_{\B}(g)^{\,-1}\;.
\end{equation}
For every $K_\val$-vector subspace $W$ of $V$, its {\it orthogonal
space} is the $K_\val$-vector subspace $W^\perp$ of
the dual $K_\val$-vector space $V^*$ defined by
\[
W^\perp =\{\ell\in V^*:\forall \,x\in W,\;\ell(x)=0\}\;.
\]
It is well-known that $\dim(W^\perp)=D-\dim(W)$, that $(W^\perp)^\perp
=W$ and that for every $g\in \GL(V)$, we have $(gW)^\perp= \widecheck
g \,(W^\perp)$.

\medskip
A {\it partial $R_\val$-lattice $\Lambda$ of rank $k$} in $V$, or {\it
  $k$-lattice} for short, is a discrete free $R_\val$-sub\-module of
rank $k$ generating a $k$-dimensional $K_\val$-vector subspace
$V_\Lambda$ of $V$. When $k=D$, we say that $\Lambda$ is {\it full
  $R_\val$-lattice}. We endow $V_\Lambda$ with the restriction of the
norm of $V$, hence with its unique Haar measure $\mu_{V_\Lambda}$ such
that $\mu_{V_\Lambda} (B_V(0,1) \cap V_\Lambda) =1$. We define the
{\it covolume} $\covol(\Lambda)$ of $\Lambda$ as the total mass of the
induced measure (again denoted by $\mu_{V_\Lambda}$) on the quotient
space $V_\Lambda/\Lambda$, that is,
\begin{equation}\label{eq:defcovol}
  \covol(\Lambda)=  \mu_{V_\Lambda}(V_\Lambda/\Lambda)\;.
\end{equation}
The set $\operatorname{Lat}_k(V)$ of $k$-lattices in $V$ is invariant
under the linear action of the linear group $\GL(V)$. This action of
$\GL(V)$ on $\operatorname{Lat}_k(V)$ is transitive, by taking an
$R_\val$-basis in two $k$-lattices, by completing them to two
$K_\val$-basis $\B$ and $\B'$ of $V$, and by taking the $K_\val$-linear
map sending $\B$ to $\B'$. For all $g\in \GL(V)$ and $\Lambda\in
\operatorname{Lat}_k(V)$, we have
\begin{equation}\label{eq:transfocovol}
V_{g\Lambda}= gV_{\Lambda}\quad\text{and}\quad\covol(g\Lambda)=
\frac{dg_*\mu_{V_\Lambda}}{d\mu_{V_{g\Lambda}}}\,\covol(\Lambda)\;.
\end{equation}
In particular, for every $\lambda\in K_\val$, we have
$\covol(\lambda\,\Lambda)= |\lambda|^k\covol(\Lambda)$ and if $k=D$,
then
\begin{equation}\label{eq:detetcovol}
  \covol(g\Lambda)=|\det g|\,\covol(\Lambda)\;.
\end{equation}

An {\it integral structure} (or {\it $R_\val$-structure}) on $V$ is
the choice of a full $R_\val$-lattice in $V$. Alternatively, it is the
choice of an equivalence class of $K_\val$-basis of $V$, where two
$K_\val$-bases are equivalent if their transition matrix belongs to
$\GL_D(R_\val)$. These two definitions agree by identifying the
equivalence class of a $K_\val$-basis $(b_1,\dots, b_D)$ with the
$R_\val$-lattice $R_\val\,b_1+\dots +R_\val\,b_D$ it generates. An
{\it integral $K_\val$-space} is a finite dimensional $K_\val$-vector
space $W$ endowed with an integral structure, denoted by
$W_{R_\val}$. We denote by $\GL(W_{R_\val})$ the subgroup of $\GL(W)$
preserving the integral structure $W_{R_\val}$ of $W$. The dual
$K_\val$-vector space $W^*$ will be endowed with the {\it dual }
integral structure (see the appendix \ref{appen:dualfactorlat} for
developments), denoted by $W^*_{R_\val}$ and defined by
\[
W^*_{R_\val}=\{\ell\in W^*: \forall\;x\in W_{R_\val},\;\;
\ell(x)\in R_\val\}\;.
\]
Equivalently, $W^*_{R_\val}$ is the integral structure on $W^*$ whose
equivalence class of $R_\val$-bases is the set of the dual
$K_\val$-bases of the elements in the equivalence class of the
$R_\val$-bases for $W_{R_\val}$.  This is well defined by Equation
\eqref{eq:duallinmap} since $\GL_\D(R_\val)$ is stable by inversion
and transposition.  Note that $W^{**}_{R_\val}$= $W_{R_\val}$.

For instance, we will endow the product $K_\val$-vector space
$K_\val^{\;\D}$ with its integral structure $R_\val^{\;\D}$ (or
equivalently with the equivalence class of its canonical basis
$(e_1,\dots, e_\D)$). By for instance \cite[Lem.~14.4)]{BroParPau19},
we have
\begin{equation}\label{eq:covolRv}
\covol(R_\val^{\;\D})=(\covol(R_\val))^\D=q^{(\ggg-1)\D}\;.
\end{equation}
We will endow the dual $K_\val$-vector space $(K_\val^{\;\D})^*$ with
the equivalence class of the dual basis $(e_1^*,\dots, e_\D^*)$ of
$(e_1,\dots, e_\D)$ (equivalently with the full $R_\val$-lattice
$R_\val\,e_1^*+\dots +R_\val\,e_\D^*$). For every $k$-lattice
$\Lambda$ in $V$, the pair $(V_\Lambda,\Lambda)$ is an integral
$K_\val$-space with $(V_\Lambda)_{R_\val}=\Lambda$.

\medskip
Since the standard $R_\val$-lattice $R_\val^{\;k}$ in $K_\val^{\;k}$
does not have covolume $1$ (contrarily to the case of the real field),
we define the {\it normalized covolume} of a $k$-lattice
$\Lambda$ in $V$ by
\[
\overline{\covol}(\Lambda)=
\frac{\covol(\Lambda)}{\covol(R_\val^{\;k})}\;.
\]

\medskip
Let $V$ be an integral $K_\val$-space with finite dimension $D$, and
$k\in\llbracket 1,D\rrbracket$. A $k$-lattice in $V$ is

$\bullet$~ {\it unimodular} if its normalized covolume
$\overline{\covol}(\Lambda)$ is equal to $1$;

$\bullet$~ {\it rational} if it is contained in the $K$-vector space
$V_K=V_{R_\val}\otimes K$ generated by the integral structure
$V_{R_\val}$ of $V$ ;

$\bullet$~ {\it integral} if it is contained in $V_{R_\val}$;

$\bullet$~ {\it primitive} if it is integral and satisfies one of the
following equivalent properties:
\begin{enumerate}
\item the $R_\val$-module $\Lambda$ is a free direct factor of
  $V_{R_\val}$ (or equivalently, there exists an $R_\val$-basis
  $(b_1,\dots,b_D)$ of $V_{R_\val}$ such that $(b_1,\dots,b_k)$ is an
  $R_\val$-basis of $\Lambda$),
\item the $R_\val$-module $V_{R_\val}/\Lambda$ is a free $R_\val$-module
  of rank $D-k$.
\end{enumerate}

Note that this definition is appropriate in the setting where $R_\val$
is not necessarily principal, and that definitions that were
equivalent in the case of $(\RR,\QQ,\ZZ)$ instead of
$(K_\val,K,R_\val)$ no longer are. For instance, if $\Lambda$ is a
primitive $k$-lattice, then $V_\Lambda$ determines $\Lambda$, with
\[
\Lambda=V_\Lambda\cap V_{R_\val}\;.
\]
But this equality is no longer sufficient for an integral $k$-lattice
$\Lambda$ to be primitive.

Note that an integral $k$-lattice is a rational $k$-lattice. By taking
an $R_\val$-basis in two rational $k$-lattices, by completing them to
two $K$-bases $\B$ and $\B'$ of $V_K$, and by taking the $K$-linear
map sending $\B$ to $\B'$, we see that the linear group $\GL(V_K)$
acts transitively on the set of rational $k$-lattices in $V$.

\medskip
Let $\La(V)$ be the space of unimodular full $R_\val$-lattices in
$V$. The closed unimodular subgroup
\begin{equation}\label{eq:defiGL1}
\GL^1(V)=\{g\in\GL(V):|\det g\,|=1\}
\end{equation}
acts transitively on the set $\operatorname{Lat}_k(V)$ of $k$-lattices
if $k<D$. It also acts transitively on $\La(V)$ (the determinant of
every element $g\in \GL(V)$ mapping a unimodular full $R_\val$-lattice
to another one has absolute value $1$ by Equation
\eqref{eq:detetcovol}). Note that the discrete group $\GL(V_{R_\val})$
is contained in $\GL^1(V)$, and is exactly the stabilizer in
$\GL^1(V)$ of the full $R_\val$-lattice $V_{R_\val}$. When
$V_{R_\val}$ is unimodular, we hence identify from now on the set
$\GL^1(V)/\GL(V_{R_\val})$ with $\La(V)$ by the map
$g\GL(V_{R_\val})\mapsto gV_{R_\val}$. In particular, we identify
$\La_D=\La(K_\val^{\;D})$ with $\GL_D^1(K_\val)/\GL_D(R_\val)$ (by
taking the matrix of a linear automorphism of $K_\val^{\;D}$ in the
canonical basis of $K_\val^{\;D}$).

Let $\PL_{k}(V)$ be the set of primitive $k$-lattices in $V$, and
$\PL_{\d,\D}=\PL_{\d}(K_\val^{\;\D})$. This set can be described as a
(discrete) homogeneous space as follows. For every commutative ring
$A$ and for all $k,k'\in\NN \smallsetminus\{0\}$, we denote by
$\M_{k,k'}(A)$ (and by $\M_{k}(A)$ when $k=k'$) the $A$-module of
$k\times k'$ matrices with coefficients in $A$.  Let
\begin{equation}\label{eq:defiPplus}
P^+(R_\val)=\Big\{\big(\begin{smallmatrix}
  \alpha &\ga \\ 0 & \delta\end{smallmatrix}\big):
  \alpha\in \GL_\d(R_\val),\;\delta\in \GL_\dd(R_\val),
  \;\ga\in \M_{\d,\dd}(R_\val),\;\det(\alpha)\det(\delta)=1 
  \Big\}\;.
\end{equation}

\blemm\label{lem:descriphomogPL}
The group $\Ga=\SL_\D(R_\val)$ acts transitively on $\PL_{\d,\D}$.
\elemm

The validity of this lemma is one of the main reasons for our
definition of primitive $\d$-lattices, and could be no longer true
with other definitions. Since the stabilizer in $\SL_\D(R_\val)$ of
the first coordinates primitive $\d$-lattice
$R_\val^{\;\d}\times\{0\}$ is equal to $P^+(R_\val)$, we will from now
on, as we may, identify the quotient $\Ga/P^+(R_\val)$ and
$\PL_{\d,\D}$ by the map
\[
gP^+(R_\val)\mapsto \Lambda_g=g(R_\val^{\;\d}\times\{0\})\;.
\]

\medskip
\dem Let $\Lambda\in\PL_{\d,\D}$. By the definition of a primitive
$\d$-lattice in $K_\val^{\;\D}$, there exists an $R_\val$-basis
$(b_1,\dots, b_\D)$ of $R_\val^{\;\D}$ such that $(b_1,\dots, b_\d)$
is an $R_\val$-basis of $\Lambda$. Let $g$ be the transition matrix
from the canonical basis of $R_\val^{\;\D}$ to $(b_1,\dots, b_\D)$.
Note that a priori $g\in\GL_\D(R_\val)$, but since $\det g\in
R_\val^{\;\times}$, up to replacing $b_\D$ by $\lambda \,b_\D$ for
some $\lambda\in R_\val^{\;\times}$, which does not change $\Lambda_g$
since $\d<\D$, we may assume that $g\in\SL_\D(R_\val)$. Then
$\Lambda_g =\Lambda$, and the map $g\mapsto \Lambda_g$ from
$\SL_\D(R_\val)$ to $\PL_{\d,\D}$ is indeed onto.  \cqfd

\subsection{Orthogonal primitive partial lattices}
\label{subsect:ortholatt}

Let $V$ be an integral $K_\val$-space with finite dimension $D$, and
let $k\in\llbracket 1,D-1\rrbracket$. Let $\Lambda$ be a primitive
$k$-lattice in $V$.

The {\it orthogonal $(D-k)$-lattice} of $\Lambda$ is the
$R_\val$-submodule of the dual integral $K_\val$-vector space
$V^*$ defined by
\[
\Lambda^\perp = {V_\Lambda}^\perp\cap V^*_{R_\val}\;.
\]
For instance, if $V$ is $K_\val^{\;\D}$ with its canonical basis
$(e_1, \dots, e_\D)$ (defining its integral structure) and its dual
basis $(e_1^*,\dots, e_\D^*)$, if $\Lambda$ is the primitive
$\d$-lattice $\oplus_{1\leq i\leq \d} R_\val e_i$, that we have
already denoted by $R_\val^{\;\d} \times \{0\}$, then
$\Lambda^\perp=\oplus_{\d+1\leq i\leq \D} R_\val e^*_i$, that we will
also denote by $\{0\}\times R_\val^{\;\dd}$.

\bprop\label{prop:propriortholat} The $R_\val$-submodule
$\Lambda^\perp$ of $V^*$ is a primitive $(D-k)$-lattice in $V^*$. For
every $g\in \GL(V_{R_\val})$, we have
\begin{equation}\label{eq:proportholat}
V_{\Lambda^\perp}= (V_\Lambda)^\perp,\qquad \Lambda=
(\Lambda^\perp)^\perp\qquad\text{and}\qquad
(g\Lambda)^\perp =\widecheck{g}\,\Lambda^\perp\;.
\end{equation}
Furthermore, if we endow $V$ with the supremum norm associated with
any $R_\val$-basis of $V_{R_\val}$ and $V^*$ with its dual norm, then
$\Lambda$ and $\Lambda^\perp$ have the same normalized covolume:
\begin{equation}\label{covollatortholat}
  \overline{\covol}(\Lambda^\perp)=
  \frac{\covol(\Lambda^\perp)}{\covol(R_\val^{\;D-k})}=
  \frac{\covol(\Lambda)}{\covol(R_\val^{\;k})}=\overline{\covol}(\Lambda)\;.
\end{equation}
\eprop

\dem Since the $k$-lattice $\Lambda$ is primitive, there exists an
$R_\val$-basis $(b_1,\dots, b_D)$ of $V_{R_\val}$ such that
$(b_1,\dots, b_k)$ is an $R_\val$-basis of $\Lambda$. The dual
$K_\val$-basis $(b_1^*, \dots, b_D^*)$ of $(b_1,\dots, b_D)$ is also
an $R_\val$-basis of the integral structure $V^*_{R_\val}$ of
$V^*$. We have $V_\Lambda=\oplus_{1\leq i\leq k} K_\val b_i$, hence
${V_\Lambda}^\perp=\oplus_{k+1\leq i\leq D} K_\val b_i^*$.  Therefore
$\Lambda^\perp={V_\Lambda}^\perp\cap V^*_{R_\val}= \oplus_{k+1\leq
  i\leq D} R_\val b_i^*$ is an integral and primitive $(D-k)$-lattice
in $V^*$. Moreover, we have $V_{\Lambda^\perp}=\oplus_{k+1\leq i\leq D}
K_\val b_i^*={V_\Lambda}^\perp$, and $(\Lambda^\perp)^\perp=
\oplus_{1\leq i\leq k} R_\val b_i=\Lambda$.

For every $g\in \GL(V)$, we have
\[
(g\Lambda)^\perp=(V_{g\Lambda})^\perp\cap V^*_{R_\val}=
(gV_{\Lambda})^\perp\cap V^*_{R_\val}=
\widecheck{g}\,(V_{\Lambda})^\perp\cap V^*_{R_\val}\;.
\]
In particular, if $g\in \GL(V_{R_\val})$, then $\widecheck{g} \in
\GL(V^*_{R_\val})$ and we do have $(g\Lambda)^\perp =\widecheck{g}\,
\Lambda^\perp$.

For a proof of Equation \eqref{covollatortholat}, we refer to the end
of the appendix \ref{appen:dualfactorlat}.
\cqfd

\subsection{Congruence properties on primitive partial lattices}
\label{subsect:congruence}

In this subsection, we fix a nonzero ideal $I$ of the Dedekind ring
$R_\val$, and we define a class of primitive partial $R_\val$-lattices
in $K_\val^{\;\D}$ that have specific congruence properties modulo the
ideal $I$.

A primitive partial $R_\val$-lattice $\Lambda$ in $K_\val^{\;\D}$ is
said to be {\it horizontal modulo $I$} if $\Lambda\subset
R_\val^{\;\d}\times I^{\dd}$, as for instance
$R_\val^{\;\d}\times\{0\}$. We will denote by $\PL_{\d,\D}(I)$ the set
of primitive $\d$-lattices in $K_\val^{\;\D}$ that are horizontal
modulo $I$. If $I=R_\val$, then $\PL_{\d,\D}(I)= \PL_{\d,\D}$.

Let $\Ga=\SL_\D(R_\val)$. We consider the following {\it Hecke
  congruence subgroup by blocks}~:
\[
\Ga_I=\{\big(\begin{smallmatrix} \alpha &\ga \\ \beta & \delta
\end{smallmatrix}\big)\in \Ga:\;\; \beta\in\M_{\dd,\d}(I)\}\;.
\]
Note that $\Ga_{R_\val}=\Ga$ and $P^+(R_\val)\subset \Ga_I$ where
$P^+(R_\val)$ is defined in Equation \eqref{eq:defiPplus}.

The first assertion of the following lemma is a congruence version of
Lemma \ref{lem:descriphomogPL}, and implies that the map $gP^+(R_\val)
\mapsto \Lambda_g =g(R_\val^{\;\d}\times\{0\})$ for $g\in\Ga_I$
identifies $\Ga_I/P^+(R_\val)$ with $\PL_{\d,\D}(I)$. The second one
is, when $\D=2$ and $\d=1$, the version of
\cite[Lem.~16.5]{BroParPau19} with $\SL_2(R_\val)$ instead of
$\GL_2(R_\val)$.

\blemm \label{lem:indexHeckesubgrou} (1) The group $\Ga_I$ acts
transitively on $\PL_{\d,\D}(I)$. Furthermore, for every $g\in \Ga$,
we have $\Lambda_g=g(R_\val^{\;\d}\times\{0\})\in \PL_{\d,\D}(I)$ if
and only if $g\in \Ga_I$.

(2) We have
\[
  [\Ga:\Ga_I]=N(I)^{\d\,\dd}\prod_{\ppp\,|I}\prod_{i=1}^\d
  \frac{N(\ppp)^{i}-N(\ppp)^{-\dd}}{N(\ppp)^i-1}\;,
\]
where the first product ranges over the prime factors $\ppp$ of the
ideal $I$.
\elemm

\dem (1) Let $g=\big(\begin{smallmatrix} \alpha &\ga \\ \beta & \delta
\end{smallmatrix}\big)\in \Ga$ with $\beta$ an $\dd\times\d$ matrix,
so that
\[
\Lambda_g=g(R_\val^{\;\d}\times\{0\})=\{(\alpha x,\beta x):x\in
R_\val^{\;\d}\}\;.
\]
Then $\Lambda_g\in \PL_{\d,\D}(I)$ if and only if $\beta x\in I^{\dd}$
for every $x\in R_\val^{\;\d}$, which occurs if and only if $\beta\in
\M_{\dd,\d}(I)$, that is, when $g\in \Ga_I$. The first claim of
Assertion (1) hence follows from Lemma \ref{lem:descriphomogPL}.

\medskip
(2) We denote by $|E|$ the cardinality of a finite set $E$. For every
commutative ring $A$ with unity whose group of invertible elements
$A^\times$ is finite and for every $\ell\in\NN\smallsetminus\{0\}$, it
is well known that $[\GL_\ell(A):\SL_\ell(A)]= |A^\times|$ and that if
$A$ is a finite field, then
\begin{equation}\label{eq:cardGL}
|\GL_\ell(A)|=|A|^{\frac{\ell(\ell-1)}{2}}
\prod_{i=1}^{\ell} (|A|^i-1)\;.
\end{equation}
The group morphism of reduction modulo $I$ from $\SL_\d(R_\val)$ to
$\SL_\d(R_\val/I)$ is onto, by an argument of further reduction to the
various prime power factors of $I$ and of lifting elementary matrices.
The order of the upper triangular subgroup by blocks
\[T_{\d,\D}(I)=\{\big(
\begin{smallmatrix} \alpha &\ga \\ 0 &\delta \end{smallmatrix}
\big)\in \SL_\D(R_\val/I):\;\;\alpha\in\M_{\d}(R_\val/I)\}
\]
of $\SL_\D(R_\val/I)$ is $|(R_\val/I)^\times|^{-1}|\GL_\d(R_\val/I)|\;
|\GL_{\dd}(R_\val/I)|\; |R_\val/I|^{\d\,\dd}$.  Hence
\begin{equation}\label{eq:premcalcindHeckegroup}
[\Ga:\Ga_I]=\frac{|\SL_\D(R_\val/I)|}{|T_{\d,\D}(I)|}
=\frac{|\GL_\D(R_\val/I)|}{|\GL_\d(R_\val/I)|\;|\GL_{\dd}(R_\val/I)|\;
  |R_\val/I|^{\d\,\dd}}\;.
\end{equation}
By the multiplicativity of the norm and by the Chinese remainder
theorem, the result reduces to the case when $I=\ppp^k$ is the $k$-th
power of a fixed prime ideal $\ppp$ of $R_\val$, where $k\in\NN$. Let
$N=N(\ppp)$ so that $N(I)=|R_\val/I|=N^k$, and note that $R_\val/\ppp$
is a field of order $N$. For every $\ell\in\NN\smallsetminus\{0\}$,
the kernel of the morphism of reduction modulo $\ppp$ from
$\GL_\ell(R_\val/\ppp^k)$ to $\GL_\ell(R_\val/\ppp)$ has order
$N^{\ell^2(k-1)}$. Hence by Equation \eqref{eq:cardGL}, we have
\[
|\GL_\ell(R_\val/I)|= N^{\ell^2(k-1)+\frac{\ell(\ell-1)}{2}}
\prod_{i=1}^{\ell} (N^i-1)\;.
\]
Therefore, by Equation \eqref{eq:premcalcindHeckegroup}, we have after
simplifications
\[
  [\Ga:\Ga_I]= N^{\d\,\dd(k-1)}\prod_{i=1}^\d\frac{N^{i+\dd}-1}{N^i-1}
  =N(I)^{\d\,\dd}\prod_{i=1}^\d
  \frac{N(\ppp)^{i}-N(\ppp)^{-\dd}}{N(\ppp)^i-1}\;.
\]
This proves the result.   \cqfd

\subsection{Refined LU decomposition by blocks}
\label{sect:LUbloc}

Let $G=\SL_\D(K_\val)$, which is a unimodular totally disconnected
locally compact topological group. In this subsection, we define some
closed subgroups of $G$ and we study their Haar measures. We will
denote an element $g\in G$ by blocks as $g=\big(\begin{smallmatrix}
  \alpha &\ga \\ \beta & \delta\end{smallmatrix}\big)$ with $\alpha$
  an $\d\times\d$ matrix. For every $k\in\NN \smallsetminus\{0\}$,
  let $I_k$ be the identity $k\times k$ matrix.

We will consider throughout this paper the following subgroups of $G$.
Let
\[
U^-=\big\{\big(\begin{smallmatrix} I_\d & 0 \\ \beta &
  I_\dd\end{smallmatrix}\big):
\beta\in\M_{\dd,\d}(K_\val)\big\}\quad\text{and}\quad
U^+=\big\{\big(\begin{smallmatrix} I_\d & \ga \\ 0 & I_\dd
\end{smallmatrix} \big):\ga\in\M_{\d,\dd}(K_\val)\big\}
\]
be the lower and upper unipotent triangular subgroups by blocks of
$G$.  For every $k\in\NN\smallsetminus\{0\}$ and every subring $A$ of
$K_\val$, we define $\GL_k^1(A)=\{g\in\GL_k(A):|\det g| =1\}$. The
closed subgroup $\GL_k^1(K_\val)$ of $\GL_k(K_\val)$ is a split
extension of its normal closed subgroup $\SL_k(K_\val)$ by the compact
group $\{\big(\begin{smallmatrix} a &
  0\\0&I_{k-1}\end{smallmatrix}\big) :a\in
\OOO_\val^{\;\times}\}$. Let
\[
G''=\big\{\big(\begin{smallmatrix}\alpha &0\\ 0 &\delta
\end{smallmatrix} \big): \alpha\in \GL^1_\d(K_\val),\; \delta\in
\GL^1_\dd(K_\val),\;\;\det\alpha\;\det\delta=1 \big\}
\]
be the intersection with $G$ of the product group $\GL^1_\d(K_\val)
\times \GL^1_\dd(K_\val)$ diagonally embedded by blocks in
$\GL_\D(K_\val)$.  Note that $U^-,G'',U^+$ are closed unimodular
subgroups of $G$, and that $G''$ normalizes $U^-$ and $U^+$.  Let
\[
Z=\big\{\big(\begin{smallmatrix} \pi_\val^{\;r}I_\d & 0 \\ 0 &
  \pi_\val^{\;s}I_{\dd}\end{smallmatrix} \big):r,s\in\ZZ, \;\d r+\dd s=0
\big\}\;,
\]
which is a discrete abelian subgroup of $G$ that centralises $G''$
(actually $G''$ is the centralizer of $Z$ in $G$), and normalizes
$U^-$ and $U^+$.

Let $Z'=\Big\{\left(\begin{smallmatrix}\pi_\val^{\;i} & 0 & 0\\ 0 &
  I_{\D-2} & 0\\ 0 & 0 &\pi_\val^{\;-i}\end{smallmatrix}\right):0\leq
i\leq \lcm\{\d,\dd\}-1\Big\}$ which is a finite subset of order
$\lcm\{\d,\dd\}$ of $G$ (not a subgroup).  Let
\begin{equation}\label{eq:defiUsugG}
\U_G=\big\{\,\big(\begin{smallmatrix}\alpha &\ga\\ \beta & \delta
\end{smallmatrix}\big)\in G:\val(\det(\alpha))\in
\lcm\{\d,\dd\}\;\ZZ\;\big\}
\end{equation}
and $\U^\bullet_G=\big\{\,\big(\begin{smallmatrix}\alpha &\ga\\ \beta & \delta
\end{smallmatrix}\big)\in G:\det(\alpha)=0\;\big\}$, which are disjoint
closed subsets of $G$ (not subgroups), with $\U_G$ open in $G$, such
that
\[
G=\U^\bullet_G\sqcup \bigsqcup_{z'\in Z'} z'\,\U_G= \U^\bullet_G\sqcup
\bigsqcup_{z'\in Z'} \U_G\,z'
\]
is a finite disjoint union of $\U^\bullet_G$ and finitely many left
(or right) translates of $\U_G$. Let $\Scal^\pm_\D$ be the subgroup of
$G$ consisting in the elements of $G$ that act by a permutation and a
possible change of sign on the elements of the canonical basis of
$K_\val^{\;\D}$ (in order for their determinant to be $1$).
Multiplying on the left an element $g\in G$ by an element in
$\Scal^\pm_\D$ amounts to permuting the rows of $g$ by the inverse of
the associated permutation and possibly changing their sign. Since the
rank of the submatrix consisting in the first $\d$ columns of any
invertible $\D\times\D$ matrix is equal to $\d$, we have
$\U^\bullet_G\subset \Scal^\pm_\D(G\smallsetminus\U^\bullet_G)$. Hence
up to left translations by the finitely many elements in $Z'$ or
$\Scal^\pm_\D$, in order to study equidistribution properties of
natural families of points in (products of) homogeneous spaces of $G$,
we may concentrate on the study of these points in the part
corresponding to $\U_G$.

For every closed subgroup $H$ of $G$, we denote by $H(\OOO_\val)$ the
compact-open subgroup $H\cap \GL_{\D}(\OOO_\val)$ of $H$, and by
$\mu_H$ the left Haar measure of $H$ normalized so that
\begin{equation}\label{eq:normalhaarsubgroup}
\mu_H(H(\OOO_\val))=1\;.
\end{equation}
In particular, $G(\OOO_\val)=\SL_{\D}(\OOO_\val)$ and
$\mu_G(G(\OOO_\val)) =1$.  Note that $\GL_\D^1(\OOO_\val)
=\GL_\D(\OOO_\val)$ since $|x|=1$ for every $x\in\OOO_\val^{\;\times}$
(see also Lemma \ref{lem:maxcomptrans}) and similarly
$\GL_\D^1(R_\val)=\GL_\D(R_\val)$.  For all $k,k'\in\NN
\smallsetminus\{0\}$, we endow the locally compact additive group
$\M_{k,k'}(K_\val)$ with its Haar measure $\operatorname{Haar}_{k,k'}$
normalized so that
\[
\operatorname{Haar}_{k,k'}(\M_{k,k'}(\OOO_\val))=1\;.
\]
The group $\GL_k(K_\val)\times\GL_{k'}(K_\val)$ acts linearly on the
$K_\val$-vector space $\M_{k,k'}(K_\val)$ by the action
$\phi(g,h):x\mapsto gxh^{-1}$ for all $x\in\M_{k,k'}(K_\val)$ and
$(g,h)\in\GL_k(K_\val)\times\GL_{k'} (K_\val)$. The following claim is
well-known.

\blemm This action scales the Haar measure
$\operatorname{Haar}_{k,k'}$ as follows:
\begin{equation}\label{eq:haarlin}
  \forall\;(g,h)\in\GL_k(K_\val)\times\GL_{k'}(K_\val),\quad
  \phi(g,h)_*\operatorname{Haar}_{k,k'}=
  |\det g\,|^{k'}\;|\det h\,|^{-k}\operatorname{Haar}_{k,k'}\;.
\end{equation}
\elemm

\dem For every $(g,h)\in\GL_k(K_\val)\times\GL_{k'}(K_\val)$, we have
$\phi(g,h)=\phi(g,\id)\circ\phi(\id,h)$. The element $\phi(g,\id)$
acts on $x\in\M_{k,k'}(K_\val)$ by the diagonal linear action of $g$
on the $k'$ columns of $x$, and $\phi(\id,h)$ acts on $x\in\M_{k,k'}
(K_\val)$ by the transpose of the diagonal linear action by the
transpose-inverse of $h$ on the $k$ columns of $^tx$, since
$^t(x\,h^{-1})=\;^th^{-1}\;^tx$.  The result hence follows from
Equation \eqref{eq:homothetyhaar} and a diagonal by block computation
of determinants.  \cqfd

\medskip
The maps $\uuu^-:\M_{\dd,\d}(K_\val)\ra U^-$ and $\uuu^+:
\M_{\d,\dd}(K_\val) \ra U^+$ defined respectively by $\beta\mapsto
\big(\begin{smallmatrix}I_\d &0\\ \beta & I_{\dd}\end{smallmatrix}
\big)$ and $\ga\mapsto \big(\begin{smallmatrix}I_\d &\ga\\ 0 & I_{\dd}
\end{smallmatrix}\big)$ are topological group isomorphisms, satisfying
\begin{equation}\label{eq:normalisehaar}
  {\uuu^-}_*\operatorname{Haar}_{\dd,\d}=\mu_{U_-}\quad\text{and}\quad
  {\uuu^+}_*\operatorname{Haar}_{\d,\dd}=\mu_{U_+}\,.
\end{equation}

We denote by $\chi_\d$ (respectively $\chi_{\dd}$) the characters from
$Z$ to $K_\val^{\;\times}$ sending $\big(\begin{smallmatrix} \lambda I_\d
  & 0\\ 0 & \mu I_{\dd} \end{smallmatrix} \big)$ to $\lambda$
(respectively $\mu$). Note that ${\chi_\d}^\d\,{\chi_{\dd}}^{\dd}$ is
the trivial character.  For all $z\in Z$, $\beta\in
\M_{\dd,\d}(K_\val)$ and $\ga\in \M_{\d,\dd}(K_\val)$, we have
\begin{align}\label{eq:dilatlarate}
  &z\,\uuu^-(\beta)\,z^{-1}=
  \uuu^-(\,(\chi_{\dd}(z) I_{\dd})\,\beta\,(\chi_\d(z)I_\d)^{-1})
  \;.
\end{align}
The Haar measure $\mu_Z$ on the closed subgroup $Z$ of $G$ is exactly
the counting measure, since $Z(\OOO_\val)=\{I_\D\}$~:
\begin{equation}\label{eq:muZ}
\mu_Z=\sum_{z\in Z} \;\Delta_z\;.
\end{equation}

The next result gives a refined LU decomposition by blocks of $G$ and
the corresponding decomposition of its Haar measure.

\bprop \label{prop:refLUbloc} The product map $(u^-,g'',z,u^+)\mapsto
u^-\,g''\,z\,u^+$ from $U^-\times G''\times Z\times U^+$ to $\U_G$ is
a homeomorphism. With
\[
c_1=\frac{q_\val^{\;\d\,\dd}\;\prod_{i=1}^{\d}
  (q_\val^{\;i}-1)\prod_{i=1}^{\dd}(q_\val^{\;i}-1)}
{\prod_{i=1}^{\D}(q_\val^{\;i}-1)}\leq 1\;,
\]
we have
\[
d\mu_G(u^-\,g''\,z\,u^+)=c_1\;|\chi_\d(z)|^{\D\,\d}\,d\mu_{U^-}(u^-)
\;d\mu_{G''}(g'')\;
d\mu_{Z}(z)\;d\mu_{U^+}(u^+)\;.
\]
\eprop

\dem For every $g=\big(\begin{smallmatrix} \alpha &\ga \\ \beta &
\delta \end{smallmatrix}\big)\in G$ such that $\det\alpha\neq 0$
(which is the case if $g\in \U_G$), we have
\begin{equation}\label{eq:decompgsimp}
g=\big(\begin{smallmatrix} I_\d &0 \\ \beta\alpha^{-1} & I_{\dd}
\end{smallmatrix}\big)\big(\begin{smallmatrix} \alpha &\ga \\ 0 &
  \delta- \beta\alpha^{-1}\ga\end{smallmatrix}\big)=
  \big(\begin{smallmatrix} I_\d &0 \\ \beta\alpha^{-1} & I_{\dd}
\end{smallmatrix}\big)\big(\begin{smallmatrix} \alpha &0 \\ 0 &
  \delta- \beta\alpha^{-1}\ga\end{smallmatrix}\big)
  \big(\begin{smallmatrix}I_\d &\alpha^{-1}\ga\\ 0 & I_{\dd}
\end{smallmatrix}\big)\;.
\end{equation}
In particular, the matrix $\delta- \beta\alpha^{-1}\ga$ is invertible,
and $\det(\alpha)\det(\delta- \beta\alpha^{-1}\ga)=1$, so that
$\val(\det(\delta- \beta\alpha^{-1}\ga))=-\val(\det(\alpha))$. Thus if
$g\in\U_G$, then $\val(\det(\alpha))$ is divisible by $\d$ and
$\val(\det(\delta- \beta\alpha^{-1}\ga))$ is divisible by
$\dd$. Furthermore,
\[
\Big|\det\Big(\pi_\val^{\;-\frac{\val(\det\alpha)}{\d}}\alpha\Big)\Big|=
|\pi_\val^{\;-\val(\det\alpha)}|\,|\det \alpha|=
q_\val^{\;\val(\det\alpha)}\,|\det \alpha|=
|\det\alpha|^{-1}\,|\det \alpha|=1\;.
\]
  
Consider the map $\Xi$ from $\U_G$ to $U^-\times G''\times Z\times U^+$
which to $g=\big(\begin{smallmatrix}\alpha &\ga \\ \beta &
  \delta\end{smallmatrix}\big)\in \U_G$ associates
\begin{align}
&\left(u^-=\begin{pmatrix}I_\d &0
    \\ \beta\alpha^{-1} & I_{\dd}\end{pmatrix},\right. g''=
  \begin{pmatrix}\pi_\val^{\;-\frac{\val(\det\alpha)}{\d}}\alpha &0
    \\ 0 &  \pi_\val^{\;-\frac{\val(\det(\delta-\beta\alpha^{-1}\ga))}{\dd}}
    (\delta-\beta\alpha^{-1}\ga)\end{pmatrix}
  \nonumber\\&z=\begin{pmatrix}\pi_\val^{\;\frac{\val(\det\alpha)}{\d}}I_\d
  & 0\\ 0 &
  \pi_\val^{\;\frac{\val(\det(\delta-\beta\alpha^{-1}\ga))}{\dd}}I_\dd\end{pmatrix},\left.
  u^+=\begin{pmatrix}I_\d &\alpha^{-1}\ga\\ 0 &
  I_{\dd}\end{pmatrix}\right) \;, \label{eq:inverseLUbloc}
\end{align}
which is well defined as we just checked. Let us prove that $\Xi$ is
onto. Let
\[
u^-=\big(\begin{smallmatrix}I_\d &0\\ \beta'&I_{\dd}\end{smallmatrix}
\big)\in U^-,\quad g''=
\big(\begin{smallmatrix}\alpha'&0\\ 0&\delta'\end{smallmatrix}
  \big)\in G'',\quad z=\big(\begin{smallmatrix}\pi_\val^{\;r}I_\d& 0\\
    0 &\pi_\val^{\;s}\end{smallmatrix}\big)\in Z,\quad
  u^+=\big(\begin{smallmatrix}I_\d &\ga'\\ 0 &
  I_{\dd}\end{smallmatrix}\big)\in U^+\;.
\]
Let $\alpha=\pi_\val^{\;r}\alpha'$, $\beta=\beta'\alpha$, $\ga=
\alpha\,\ga'$ and $\delta=\pi_\val^{\;s}\delta'+ \beta\alpha^{-1}\ga$.
The equality $\d\,r+\dd\,s=0$ implies by Gauss Lemma that
$\frac{\dd}{\gcd\{\d,\dd\}}$ divides $r$, hence that $\lcm\{\d,\dd\}=
\frac{\d\,\dd}{\gcd\{\d,\dd\}}$ divides $\d\,r$. Since $\alpha'\in
\GL^1_\d(K_\val)$, we have
\begin{equation}\label{eq:valeurr}
\val(\det\alpha)=\val((\pi_\val^{\;r})^\d\det\alpha')=
\d\,r+\val(\det\alpha')=\d\,r\in\lcm\{\d,\dd\}\,\ZZ\;.
\end{equation}
Thus $g=\big(\begin{smallmatrix}\alpha &\ga\\ \beta&\delta
\end{smallmatrix} \big)$ belongs to $\U_G$ and by construction
$\Xi(g)=(u^-,g'',z,u^+)$. It is immediate to see that $\Xi$ is
continuous on $\U_G$ and is the inverse of the continuous
multiplication map $(u^-,g'',z,u^+)\mapsto u^-\,g''\,z\,u^+$.

Let $P^-=U^-G''Z$, which is a closed subgroup of $G$,
%\footnote{and an index $\d$ subgroup of the parabolic
%subgroup $Z'P^{-}$ of $G$ which is the stabilizer of the
%$K_\val$-linear subspace ${K_\val}^d\times\{0\}$ for
%the linear action of $G$ on $K_\val^{\;\D}={K_\val}^d{K_\val}^{\dd}$}
since $Z$ centralises $G''$ and $G''Z$ normalizes $U^-$, so that
$\U_G = P^-U^+$. By \cite[\S III.1]{Lang75}, since $G$ and $U^+$ are
unimodular, there exists a constant $c_2>0$ such that the restriction
to the open set $\U_G$ of the Haar measure of $G$ satisfies
$d\mu_G(p^-u^+) = c_2\; d\mu_{P^-}(p^-)\;d\mu_{U^+}(u^+)$ for (almost)
all $p^-\in P^-$ and $u^+\in U^+$.

For all $g''= \big(\begin{smallmatrix}\alpha &0\\ 0
  &\delta \end{smallmatrix} \big)\in G''$ and $\beta\in
\M_{\dd,\d}(K_\val)$, we have $g''\uuu^-(\beta){g''}^{-1}= \uuu^-(\delta
\beta \alpha^{-1})$. Therefore, with $\iota_{g''}$ the conjugation map
$x\mapsto g''\,x\,{g''}^{-1}$ by $g''$ on $U^-$, by Equation
\eqref{eq:normalisehaar}, by Equation \eqref{eq:haarlin} with $k=\dd$
and $k'=\d$, and since $|\det \alpha| =|\det \delta|=1$, we have
\begin{align}
  (\iota_{g''})_*\mu_{U^-}&=
  (\iota_{g''}\circ\uuu^-)_*\operatorname{Haar}_{\dd,\d}
  =(\uuu^-\circ\phi(\delta,\alpha))_*\operatorname{Haar}_{\dd,\d}
  \nonumber\\ &=(\uuu^-)_*\operatorname{Haar}_{\dd,\d}=\mu_{U^-}
  \label{eq:conjpreshaar}\;.
\end{align}
Let $z\in Z$. By the relationship between the two characters $\chi_\d$
et $\chi_{\dd}$ of $Z$, we have
\begin{align*}
&|\det(\chi_{\dd}(z) I_{\dd})|^\d \,
|\det(\chi_\d(z)I_\d)|^{-\dd}= |\chi_{\dd}(z)^{\dd}|^\d\,
|\chi_\d(z)^{\d}|^{-\dd}\\ =\;&|\chi_{\d}(z)^{-\d}|^\d\,
|\chi_\d(z)|^{-\d\,\dd}=|\chi_\d(z)|^{-\D\,\d}\;.
\end{align*}
Therefore, with $\iota_{z}$ the conjugation map by $z$ on $U^-$, by
Equations \eqref{eq:normalisehaar} and \eqref{eq:dilatlarate}, by
Equation \eqref{eq:haarlin} with $k=\dd$ and $k'=\d$, we have
\begin{align}
  (\iota_{z})_*\mu_{U^-}&=
  (\iota_{z}\circ\uuu^-)_*\operatorname{Haar}_{\dd,\d}
  =(\uuu^-\circ\phi(\chi_{\dd}(z) I_{\dd},\chi_\d(z)I_\d))_*
  \operatorname{Haar}_{\dd,\d}\nonumber\\ &
  =|\chi_\d(z)|^{-\d\D}\mu_{U^-}\label{eq:conjdilathaar}\;.
\end{align}

The product in the group $P^-=U^-G''Z$ may be written as follows: for
all $(u^-,g'',z)$ and $(\wh u\,^-,\wh g\,'',\wh z)$ in $U^-\times
G''\times Z$, we have
\[
(u^-\,g''\,z)\,(\wh u\,^-\,\wh g\,''\,\wh z\,)=
\big(u^-\,(\,g''\,(z\,\wh u^-\,z^{-1})\,{g''}^{-1})\big)\,
\big(\,g''\,\wh g\,''\big)\big(z \,\,\wh z\,\big)\;.
\]
By Equations \eqref{eq:conjdilathaar} and \eqref{eq:conjpreshaar}, the
image of the measure
\[
|\chi_\d(z)|^{\d\D}\; d\mu_{U^-}(u^-)\; d\mu_{G''}(g'')\;d\mu_{Z}(z)\]
on $U^-\times G''\times Z$ by the product map $(u^-,g'',z)\mapsto
u^-\,g''\,z$ is hence a left Haar measure on $P^-$.  Therefore there
exists a constant $c_3>0$ such that, for (almost) all $u^-\in U^-$,
$g''\in G''$ and $z\in Z$, we have
\[
|\chi_\d(z)|^{\,\d\D}\; d\mu_{U^-}(u^-)\; d\mu_{G''}(g'')\;d\mu_{Z}(z)=
c_3\; d\mu_{P^-}(u^-g''z)\;.
\]
Therefore, with $c_1=c_2c_3^{-1}$, we have
\begin{equation}\label{eq:defc1}
d\mu_G(u^-\,g''\,z\,u^+)=c_1\;|\chi_\d(z)|^{\D\,\d}\,d\mu_{U^-}(u^-)
\;d\mu_{G''}(g'')\;
d\mu_{Z}(z)\;d\mu_{U^+}(u^+)\;.
\end{equation}
In order to compute the constant $c_1$, we evaluate the measures on
both sides of Equation \eqref{eq:defc1} on the compact-open subgroup
\[
H=\Big\{\big(\begin{smallmatrix} \alpha & \ga\\ \beta & \delta
\end{smallmatrix} \big)\in G(\OOO_\val): \begin{array}{l}\alpha\in
I_\d+\pi_\val\M_{\d}(\OOO_\val),\;\delta\in I_{\dd}+\pi_\val\M_{\dd}
(\OOO_\val),\\ \beta\in \pi_\val\M_{\dd,\d}(\OOO_\val),\;
\ga \in \pi_\val\M_{\d,\dd}(\OOO_\val)\end{array}
\Big\}\;.
\]
We are going to need the following well-known result.

\blemm\label{lem:calcindexN} For every $N\in\NN$, let $H_N$ be the
kernel of the morphism from $G(\OOO_{\val})$ to
$\SL_{\D}(\OOO_{\val}/\pi_{\val}^{N}\OOO_{\val})$ of reduction modulo
$\pi_{\val}^{N}\OOO_{\val}$.  Then
\[
[G(\OOO_{\val}):H_N]=q_\val^{\;N(\D^2-1)-\frac{\D(\D+1)}{2}+1}
\prod_{i=2}^{\D}(q_\val^{\;i}-1)\;.
\]
\elemm \dem Since the ring $\OOO_\val$ is principal and by for
instance \cite[page 21]{Shimura71}, the reduction morphism
$\SL_\D(\OOO_{\val}) \to
\SL_\D(\OOO_{\val}/\pi_{\val}^{N}\OOO_{\val})$ is onto. By for
instance \cite[Theo.~2.7]{Han06} (applied with $k=N$ to the finite
local commutative ring $R=\OOO_\val/\pi_\val^{\;k}\OOO_\val$ with
maximal ideal $P=\pi_\val\OOO_\val/\pi_\val^{\;k}\OOO_\val$), and by
Equation \eqref{eq:cardGL}, we have
\begin{align*}
|\GL_\D(\OOO_{\val}/
\pi_{\val}^{N}\OOO_{\val})|&=|\pi_{\val}\OOO_{\val}/
\pi_{\val}^{N}\OOO_{\val}|^{\D^2}\;|\GL_\D(\OOO_{\val}/
\pi_{\val}\OOO_{\val})|=q_\val^{\;(N-1)\D^2+\frac{\D(\D-1)}{2}}
\prod_{i=1}^{\D}(q_\val^{\;i}-1)\;.
\end{align*}
The index of $\SL_\D(\OOO_{\val}/ \pi_{\val}^{\;N}\OOO_{\val})$ in
$\GL_\D(\OOO_{\val}/ \pi_{\val}^{\;N}\OOO_{\val})$ is equal to
$|(\OOO_{\val}/ \pi_{\val}^{N}\OOO_{\val})^\times|=q_\val^{\;N-1}
(q_\val-1)$. The result follows.
\cqfd

\medskip
Since $\mu_G(G(\OOO_\val)) =1$ and by Lemma \ref{lem:calcindexN} with
$N=1$, the group $H$ has Haar measure
\[
\mu_G(H)=\frac{\mu_G(G(\OOO_\val))}{[G(\OOO_\val):H]}=
\frac{1}{q_\val^{\;\frac{\D(\D-1)}{2}}\prod_{i=2}^{\D}(q_\val^{\;i}-1)}\;.
\]

The group $H\cap U^-=\big\{\big(\begin{smallmatrix} I_\d & 0\\ \beta &
I_{\dd}\end{smallmatrix} \big): \beta\in \pi_\val\M_{\dd,\d}(\OOO_\val)
\big\}$ has index $q_\val^{\;\d\,\dd}$ in $U^-(\OOO_\val)$, and so does 
$H\cap U^+=\big\{\big(\begin{smallmatrix} I_\d & \ga\\ 0 &
I_{\dd}\end{smallmatrix} \big): \ga\in \pi_\val\M_{\d,\dd}(\OOO_\val)
\big\}$ in $U^+(\OOO_\val)$. Hence
\[
\mu_{U^-}(H\cap U^-)=\mu_{U^+}(H\cap U^+)=\frac{1}{q_\val^{\;\d\,\dd}}\;.
\]

We have $H\cap Z=\{I_\D\}$, hence $\int_{H\cap Z} |\chi_\d(z)|^{\,\D\,\d}
\;d\mu_{Z}(z)=1$ by Equation \eqref{eq:muZ}.

The index of the subgroup $H\cap G''=\Big\{\big(\begin{smallmatrix}
  \alpha & 0\\ 0 & \delta \end{smallmatrix} \big):
\begin{array}{l}\alpha\in I_\d+\pi_\val\M_{\d}(\OOO_\val),
  \\ \delta\in I_{\dd}+\pi_\val\M_{\dd}(\OOO_\val),
\end{array}\; \det\alpha\;\det\delta=1 \Big\}$ in the group
$G''(\OOO_\val)=\big\{\big(\begin{smallmatrix} \alpha & 0\\ 0 &
  \delta \end{smallmatrix} \big): \alpha\in \GL_\d(\OOO_\val),\;
\delta\in \GL_\dd(\OOO_\val),\; \det\alpha\;\det\delta=1 \big\}$ is
equal to $\frac{1}{|\FF_{q_\val}^\times|}(|\GL_\d(\FF_{q_\val})|\times
|\GL_{\dd}(\FF_{q_\val})|)$. Since $\mu_{G''}(G''(\OOO_\val))=1$, and
by Equation \eqref{eq:cardGL}, we hence have
\[
\mu_{G''}(H\cap G'')=\frac{(q_\val-1)}
{q_\val^{\;\frac{\d(\d-1)}{2}}\prod_{i=1}^{\d}(q_\val^{\;i}-1)\;
q_\val^{\;\frac{\dd(\dd-1)}{2}}\prod_{i=1}^{\dd}(q_\val^{\;i}-1)}\;.
\]
  
Note that $H$ is contained in $\U_G$ since for every $\big(
\begin{smallmatrix} \alpha & \ga\\ \beta & \delta\end{smallmatrix}\big)
\in H$, we have $\val(\det(\alpha))=0$ as $\det\alpha \equiv\det
I_\d\equiv 1\mod \pi_\val$. We then also have $\val(\det(\delta
-\beta\alpha^{-1}\alpha))=0$. It follows from Equation
\eqref{eq:inverseLUbloc} that the product map $(u^-,g'',z, u^+)\mapsto
u^-\,g''\,z\,u^+$ from $U^-\times G''\times Z\times U^+$ to $\U_G$
induces a homeomorphism from $(H\cap U^-)\times (H\cap G'')\times
(H\cap Z)\times (H\cap U^+)$ to $H$. By Equation \eqref{eq:defc1} and
the above computations, we thus have
\begin{align*}
  c_1&=\frac{\mu_G(H)}{\mu_{U^-}(H\cap U^-)\;\mu_{G''}(H\cap G'')
    \;\mu_{U^+}(H\cap U^+)}
  \\&=\frac{q_\val^{\;\d\,\dd}\;\prod_{i=1}^{\d}(q_\val^{\;i}-1)
    \prod_{i=1}^{\dd}(q_\val^{\;\dd}-1)}
    {\prod_{i=1}^{\D}(q_\val^{\;i}-1)}\;.
\end{align*}
Note that $c_1=\prod_{i=1}^{\d}\frac{q_\val^{\;i+\dd}-\,q_\val^{\;\dd}}
{q_\val^{\;i+\dd}-\,1}\leq 1$. This concludes the proof of Proposition
\ref{prop:refLUbloc}.  \cqfd

\subsection{Refined LU decomposition by blocks and partial lattices}
\label{subsect:LUdecpartlatt}

Let $V$ be a $K_\val$-vector space with finite dimension $D\geq 1$,
and let $k\in\llbracket 1,D\rrbracket$. We denote by $\Gr_k(V)$ the
Grassmannian space of $k$-dimensional $K_\val$-linear subspaces of $V$
(endowed with the compact metrisable Chabauty topology, see Subsection
\ref{subsect:grassmann} for its measure theoretic and metric
aspects). We define
\[
\Gr_{\d,\D}=\Gr_\d(K_\val^{\;\D})\;.
\]

An {\it $\OOO_\val$-structure} on $V$ is the choice of a finitely
generated 
%\footnote{Free is automatic since $\OOO_\val$ is a principal ideal
%domain. See \cite[page 139}{GolIwa63}. Equivalently it is a compact
%and open $\OOO_\val$-submodule.}
$\OOO_\val$-submodule $V_{\OOO_\val}$ (which is a free
$\OOO_\val$-module since it is torsion free and $\OOO_\val$ is
principal) generating $V$ as a $K_\val$-vector space, or equivalently
an equivalence class of $K_\val$-basis of $V$, where two
$K_\val$-bases are equivalent if their transition matrix belongs to
$\GL_D(\OOO_\val)$. For instance, we endow $K_\val^{\;D}$ and
$(K_\val^{\;D})^*$ with the $\OOO_\val$-structure defined by their
canonical basis and its dual basis, respectively.  We denote by
$\GL(V_{\OOO_\val})$ the subgroup of $\GL(V)$ preserving
$V_{\OOO_\val}$, and we define $\SL(V_{\OOO_\val})
=\GL(V_{\OOO_\val})\cap\SL(V)$.  The following claim is well-known.

\blemm \label{lem:maxcomptrans} The group $\GL(V_{\OOO_\val})$ is
contained in $\GL^1(V)$ and acts transitively on the Grassmannian
space $\Gr_k(V)$. The group $\SL(V_{\OOO_\val})$ also acts
transitively on $\Gr_k(V)$.
\elemm

\dem The first claim follows from the fact that the determinant of a
matrix in $\GL(V_{\OOO_\val})$, being an element of
$\OOO_\val^{\;\times}$, has absolute value $1$. By for instance
\cite[Theo.~1]{Weil95}, every complete flag $(V_1,\dots, V_D)$ of $V$
admits a $K_\val$-basis $(x_1,\dots, x_D)$ which is both adapted to
this flag (that is, $V_i=K_\val x_1+\dots+K_\val x_i$ for every
$i\in\llbracket 1,D\rrbracket$) and is an $\OOO_\val$-basis of
$V_{\OOO_\val}$. Hence every $k$-dimensional $K_\val$-linear subspace
of $V$ admits a $K_\val$-basis that can be completed to an
$\OOO_\val$-basis of $V_{\OOO_\val}$. Since $\GL(V_{\OOO_\val})$ acts
transitively on the set of $\OOO_\val$-bases of $V_{\OOO_\val}$, the
second claim follows.  The last claim (clear when $k=D$)
follows when $k<D$ by multiplying the last element of the above
$\OOO_\val$-basis by an appropriate element of $\OOO_\val^{\;\times}$.
\cqfd

\medskip
Assume that $V$ is endowed with a $K_\val$-basis $(f_1,\dots, f_D)$
defining both an $R_\val$-structure $V_{R_\val}=R_\val f_1+\dots+
R_\val f_D$, an $\OOO_\val$-structure $V_{\OOO_\val}=\OOO_\val
f_1+\dots +\OOO_\val f_D$ and an ultrametric norm $\|x_1f_1+\dots +x_D
f_D\|= {\displaystyle \max_{1\leq i\leq D}}\, |x_i|$ whose (closed)
unit ball is $V_{\OOO_\val}$. For instance, $V$ could be
$K_\val^{\;D}$ with its canonical $K_\val$-basis, or its dual space
$(K_\val^{\;D})^*$ with its dual $K_\val$-basis, or $K_\val$-linear
subspaces of them generated by $K_\val$-basis elements. We denote the
{\it space of shapes} of unimodular full $R_\val$-lattices of $V$ by
\[
\Sh(V)=\GL(V_{\OOO_\val})\bs\La(V)
=\GL(V_{\OOO_\val})\bs\GL^1(V)/\GL(V_{R_\val})\;,
\]
endowed with the quotient topology (see Subsection \ref{subsect:shape}
for its measure theoretic aspects).  For simplicity, we denote
\[
\Sh_D=\Sh(K_\val^{\;D})\;.
\]

The {\it shape map} of $k$-lattices of $V$ is the map
\[
\sh:\big\{\Lambda\in\operatorname{Lat}_k(V):
\overline{\covol}(\Lambda)\in q_\val^{\;k\ZZ}\big\}\ra \Sh_k
\]
defined as follows. Let $\Lambda\in\operatorname{Lat}_k(V)$ with
$\overline{\covol}(\Lambda)\in q_\val^{\;k\ZZ}$. By Lemma
\ref{lem:maxcomptrans}, choose an element $g\in \GL(V_{\OOO_\val})$
such that we have $g V_\Lambda=K_\val f_1+\dots +K_\val f_k$. Note
that $g$ preserves the covolume of partial $R_\val$-lattices, since it
maps the unit ball $V_\Lambda\cap V_{\OOO_\val}$ of $V_\Lambda$ to the
unit ball $\OOO_\val f_1+\dots +\OOO_\val f_k$ of $K_\val f_1+\dots
+K_\val f_k$. Let $\Theta:K_\val f_1+\dots +K_\val f_k\ra
K_\val^{\;k}$ be the isometric (hence Haar measure preserving)
$K_\val$-linear isomorphism mapping $(f_1,\dots, f_k)$ to the
canonical basis of $K_\val^{\;k}$, that preserves the covolume of full
$R_\val$-lattices. We define
\begin{equation}\label{eq:defish}
\sh(\Lambda)= \GL_k(\OOO_\val)\pi_\val^{\;\frac{1}{k}\log_{q_\val}
  \overline{\covol}(\Lambda)}\Theta(g\Lambda)\;.
\end{equation}
Note that $\pi_\val^{\;\frac{1}{k}\log_{q_\val}\overline{\covol}
  (\Lambda)} \Lambda$ is a unimodular $k$-lattice in $V$, and that
homotheties and linear maps commute. Furthermore, the shape
$\sh(\Lambda)$ of $\Lambda$ does not depend on the choice of $g$ as
above, since given two choices $g_1$ and $g_2$, the linear maps
$\Theta\circ g_1$ and $\Theta\circ g_2$ differ by multiplication on
the left by an element of $\GL_k(\OOO_\val)$. Note that when $\Lambda$
is a unimodular full $R_\val$-lattice in the product space
$V=K_\val^{\;k}$, then Equation \eqref{eq:defish} greatly simplifies
to 
\[
\sh(\Lambda)= \GL_k(\OOO_\val)\Lambda\;,
\]
since we can take $(f_1,\dots, f_k)$ to be the canonical basis of $V$
and $g=\Theta$ to be the identity map of $V$.

\medskip
The next result gives the relationship between the refined LU
decomposition by blocks of elements of $\U_G$ and the partial
$R_\val$-lattices generated by their first $\d$ columns. Let us first
give the notation that will be used. For every $D\in\NN\ssm\{0\}$, we
identify $\GL(K_\val^{\;D})$ and $\GL_D(K_\val)$ (respectively
$\GL((K_\val^{\;D})^*)$ and $\GL_D(K_\val)$) by taking matrices of
linear automorphisms in the canonical basis $(e_1,\dots, e_D)$ of
$K_\val^{\;D}$ (respectively its dual basis $(e_1^*,\dots, e_D^*)$ of
$(K_\val^{\;D})^*$). Recall that the map $h\mapsto
\widecheck{h}=\;^th^{-1}$ is a group isomorphism from
$\GL(K_\val^{\;\D})$ to $\GL((K_\val^{\;\D})^*)$. For practical
reasons, we denote by $R_\val^{\;\d}\times\{0\}$ the $\d$-lattice
$R_\val e_1+\dots+ R_\val e_\d$ of $K_\val^{\;\D}$ and by $\{0\}\times
R_\val^{\;\dd}$ the $\dd$-lattice $R_\val e_{\d+1}^*+\dots+ R_\val
e_\D^*$ of $(K_\val^{\;\D})^*$. For every $g\in \GL(K_\val^{\;\D})$
(respectively $g'\in \GL((K_\val^{\;\D})^*)$~), we define
\[
\Lambda_g=g(R_\val^{\;\d}\times\{0\})\quad\text{(respectively}\quad
\Lambda_{g'}=g'(\{0\}\times R_\val^{\;\dd})\;)\;,
\]
which is the $\d$-lattice generated by the first $\d$ columns of $g$
(respectively the $\dd$-lattice generated by the last $\dd$ columns of
$g'$).  For every element $g=\big(\begin{smallmatrix} \alpha &\ga
  \\ \beta & \delta\end{smallmatrix}\big)\in\U_G$, we denote by
  $(u^-,g''= \big(\begin{smallmatrix} \overline{g} &0 \\ 0
    &\underline{g}
  \end{smallmatrix}\big),z,u^+)$ the decomposition of $g$ given by
  Proposition \ref{prop:refLUbloc}, and by 
\[
t=-\frac{\val(\det\alpha)}{\lcm\{\d,\dd\}} =
\frac{\val(\det(\delta-\beta\alpha^{-1}\ga))}{\lcm\{\d,\dd\}}\in\ZZ
\]
(which is indeed an integer since $g\in\U_G$ by Equation
\eqref{eq:defiUsugG}). By Equation \eqref{eq:inverseLUbloc}, we then
have
\begin{equation}\label{eq:formulezrelatt}
  z=\begin{pmatrix} \pi_\val^{\;-\frac{\lcm\{\d,\dd\}}{\d} \,t}I_\d &0\\ 0 &
    \pi_\val^{\;\frac{\lcm\{\d,\dd\}}{\dd} \,t}I_\dd \end{pmatrix}\;.
\end{equation}
To conclude this list of notation for Proposition
\ref{prop:HKprop4_3}, let us define
\[
G^\sharp=\{g\in\U_G: u^-\in G(\OOO_\val)\}\;.
\]

\bprop\label{prop:HKprop4_3}
For every $g\in \U_G$, we have\smallskip\\
\begin{tabular}{llll}
  \hypertarget{HKprop4_3i}{(i)} &$V_{\Lambda_g}=V_{\Lambda_{u^-}}$&
  \hypertarget{HKprop4_3iperp}{(i)$^\perp$}~
  &$V_{(\Lambda_g)^{\perp}}=V_{\Lambda_{\widecheck{u^-}}}$\\
  \hypertarget{HKprop4_3ii}{(ii)}
  &$\overline{\covol}(\Lambda_g)=q_\val^{\;\lcm\{\d,\dd\} \,t}$
  \hspace*{1.3cm}
  &\hypertarget{HKprop4_3iiperp}{(ii)$^\perp$}~
  &$\overline{\covol}((\Lambda_g)^{\perp})=
  %\overline{\covol}(\Lambda_g)=
  q_\val^{\;\lcm\{\d,\dd\}\,t}$\end{tabular}

\smallskip\noindent and if furthermore $g\in G^\sharp$,  then\smallskip\\
\begin{tabular}{llll}\hypertarget{HKprop4_3iii}{(iii)}
  &$\sh(\Lambda_g)=\sh(\,\overline{g}\,R_\val^{\;\d})$~~~
  &\qquad\quad\;\;\;\hypertarget{HKprop4_3iiiperp}{(iii)$^\perp$}~
  &$\sh((\Lambda_g)^{\perp})=
  \sh(\,\widecheck{\underline{g}\,}\,R_\val^{\;\dd})$.
\end{tabular}
\eprop

\dem Let $g=\big(\begin{smallmatrix} \alpha &\ga \\ \beta &
\delta\end{smallmatrix}\big)\in G\ssm\U_G^\bullet$ (that is, let $g$ be an
element of $G$ whose upper-left $\d\times\d$ submatrix is
invertible), let $u_-=\big(\begin{smallmatrix} I_\d &0\\\beta\alpha^{-1}
& I_\dd\end{smallmatrix}\big)$ and let $t=-\frac{\val(\det\alpha)}
{\lcm\{\d,\dd\}}\in\QQ$, so that these two notations coincide with the above
ones when $g\in \U_G\subset G\ssm\U_G^\bullet$. Let us prove that Assertions
\hyperlink{HKprop4_3i}{(i)}, \hyperlink{HKprop4_3iperp}{(i)$^\perp$},
\hyperlink{HKprop4_3ii}{(ii)}, \hyperlink{HKprop4_3iiperp}{(ii)$^\perp$}
are actually satisfied under this greater generality on $g$. Since we
have $\big(\begin{smallmatrix} \alpha &\ga \\ \beta
  &\delta\end{smallmatrix}\big) \big(\begin{smallmatrix}
    x_{\mbox{}}\\ 0\end{smallmatrix}\big)= \big(\begin{smallmatrix}
      I_\d &0\\\beta\alpha^{-1} & I_\dd
\end{smallmatrix}\big)\big(\begin{smallmatrix}\alpha x \\ 0
\end{smallmatrix}\big)$ for every $x\in\M_{\d,1}(K_\val)$, we have
\begin{equation}\label{eq:explicitLambdag}
\Lambda_g=g(R_\val^{\;\d}\times\{0\})=u^-((\alpha R_\val^{\;\d})\times\{0\})\;.
\end{equation}
Since $\alpha$ is invertible, we hence have $V_{\Lambda_g}= u^-
(K_\val^{\;\d} \times\{0\})$, thus proving \hyperlink{HKprop4_3i}{(i)}.

Let $\delta'=\delta-\beta\alpha^{-1}\ga$ and $u^+=\big(\begin{smallmatrix}
I_\d &\alpha^{-1}\ga \\0 &I_\dd\end{smallmatrix}\big)$. Since $g=u^-
\big(\begin{smallmatrix} \alpha &0 \\ 0 &\delta'\end{smallmatrix}\big)
  u^+$ by Equation \eqref{eq:decompgsimp}, we have $\widecheck{g}=
  \widecheck{u^-}\big(\begin{smallmatrix} \widecheck{\alpha} &0 \\ 0
    &\widecheck{\delta'}\end{smallmatrix} \big) \widecheck{u^+}$. Note
  that $\widecheck{u^+}$, being lower unipotent by blocks, preserves
  $\{0\}\times R_\val^{\;\dd}$. By the last equality in Equation
  \eqref{eq:proportholat}, we hence have
\begin{equation}\label{eq:explicitLambdagperp}
  (\Lambda_g)^{\perp}=\big(g(R_\val^{\;\d}
\times\{0\}) \big)^\perp= \widecheck{g}\,(R_\val^{\;\d}
\times\{0\})^\perp= \widecheck{g}\,(\{0\} \times
R_\val^{\;\dd})=\widecheck{u^-}(\{0\}\times (\widecheck{\delta'}
R_\val^{\;\dd}))\;.
\end{equation}
Since $\widecheck{\delta'}$ is invertible, by the left-hand side of
Equation \eqref{eq:transfocovol}, we thus have
\[
V_{(\Lambda_g)^{\perp}}= \widecheck{u^-}(\{0\}\times
K_\val^{\;\dd})=\widecheck{u^-}\,V_{\{0\}\times R_\val^{\;\dd}}=
V_{\widecheck{u^-}(\{0\}\times R_\val^{\;\dd})}\;,
\]
thereby proving
Assertion \hyperlink{HKprop4_3iperp}{(i)$^\perp$}.

By Equations \eqref{eq:explicitLambdag}, \eqref{eq:transfocovol} (its
right-hand side) and \eqref{eq:detetcovol}, since we have $\det(u^-)=1$
and by the definition of $t=- \frac{\val(\det\alpha)} {\lcm\{\d,\dd\}}
=\frac{\log_{q_\val}|\det\alpha|}{\lcm\{\d,\dd\}}$, we have
\[
\covol(\Lambda_g)=\covol(\alpha R_\val^{\;\d})=
|\det \alpha|\covol(R_\val^{\;\d})=
q_\val^{\;\lcm\{\d,\dd\} \,t}\covol(R_\val^{\;\d})\;,
\]
thus proving \hyperlink{HKprop4_3ii}{(ii)}. Assertion
\hyperlink{HKprop4_3iiperp}{(ii)$^\perp$} follows from Equation
\eqref{covollatortholat}.

Assume from now on in this proof that $g\in G^\sharp$. Note that by
Assertion \hyperlink{HKprop4_3ii}{(ii)}, the $\d$-lattice $\Lambda_g$
has normalized covolume which is an integral power of $q_\val^{\;\d}$,
hence $\sh(\Lambda_g)$ is well defined by Equation \eqref{eq:defish}
with $k=\d$. Recall that the shape of a partial $R_\val$-lattice of
$K_\val^{\;\D}$ is invariant by every homothety and by taking the
image by any element in $\GL_\D(\OOO_\val)$. Again by
Equation \eqref{eq:explicitLambdag}, since $u^-\in G(\OOO_\val)$ when
$g\in G^\sharp$, and since
$\overline{g}\,=\pi_\val^{\;-\frac{\val(\det\alpha)}{\d}}\,\alpha=
\pi_\val^{\;\frac{\lcm(\d,\dd)}{\d}\,t}\,\alpha$ by Equation
\eqref{eq:inverseLUbloc}, we have
\[
\sh(\Lambda_g)=\sh((\alpha R_\val^{\;\d})\times\{0\})=
\sh(\alpha R_\val^{\;\d})=\sh(\,\overline{g}\,R_\val^{\;\d})\;,
\]
thus proving \hyperlink{HKprop4_3iii}{(iii)}. Note that
$\overline{g}\in\GL^1_\d(K_\val)$, so that $\,\overline{g}\,
R_\val^{\;\d}$ is a unimodular full $R_\val$-lattice in $K_\val^{\;\d}$.

By Assertion \hyperlink{HKprop4_3iiperp}{(ii)$^\perp$}, the
$\dd$-lattice $(\Lambda_g)^{\perp}$ has normalized covolume which is
an integral power of $q_\val^{\;\dd}$, hence
$\sh((\Lambda_g)^{\perp})$ is well defined by Equation
\eqref{eq:defish} with $k=\dd$. As previously, since $\widecheck{u^-}$
(which is now upper unipotent by blocks) still belongs to
$G(\OOO_\val)$, and since $\underline{g}\,$ is a scalar multiple of
$\delta'=\delta-\beta\alpha^{-1}\ga$, we have by Equation
\eqref{eq:explicitLambdagperp} that
\[
\sh((\Lambda_g)^{\perp})=\sh(\{0\}\times (\,\widecheck{\delta'}
R_\val^{\;\dd}))=\sh(\,\widecheck{\delta'} R_\val^{\;\dd})=
\sh(\,\widecheck{\underline{g}}\,R_\val^{\;\dd})\;,
\]
thus proving \hyperlink{HKprop4_3iiiperp}{(iii)$^\perp$}. Note that
$\widecheck{\underline{g}}\, R_\val^{\;\dd}$ is a unimodular full
$R_\val$-lattice in $K_\val^{\;\dd}$.
\cqfd

\section{Metric measured moduli spaces of partial lattices}
\label{sect:moduli}

In this section, we define the natural measures and distances on the
moduli spaces $\Gr_{\d,\D}$ (see Subsection \ref{subsect:grassmann}),
$\Sh_\d$ and $\Sh_\dd$ (see Subsection \ref{subsect:shape}), that were
introduced in Subsection \ref{subsect:LUdecpartlatt} and on whose
products the equidistribution results of the Introduction will take
place. We introduce (and similarly analyse) in Subsection
\ref{subsect:latticepair} an avatar $\La_{\d,\dd}$ of the product
$\La_{\d}\times\La_{\dd}$ of the spaces (described in Subsection
\ref{subsect:fulllat}) of $\d$- and $\dd$-lattices, on which our
stronger equidistribution result (Theorem \ref{theo:mainsharp}) will
take place in the subsequent Section \ref{sect:Equidistribution}.

\subsection{The Grassmannian spaces}
\label{subsect:grassmann}

As defined in Subsection \ref{subsect:LUdecpartlatt}, we denote by
$\Gr_{\d,\D}=\Gr_\d(K_\val^{\;\D})$ the Grassmannian space of
$\d$-dimensional $K_\val$-linear subspaces of $K_\val^{\;\D}$.
Recalling that $G=\SL_\D(K_\val)$, we define
\[
Q^+=\big\{\big(\begin{smallmatrix} \alpha &\ga \\ \beta &
  \delta\end{smallmatrix}\big)\in G:\beta=0\big\}\;,
\]
which is a nonunimodular closed subgroup of $G$. It contains the
closed and open subgroup $P^+=G''ZU^+=\big\{\big(\begin{smallmatrix}
\alpha &\ga \\ \beta & \delta\end{smallmatrix}\big)\in Q^+:
\val(\det(\alpha))\in \lcm\{\d,\dd\}\,\ZZ\,\}$ with finite index.  The
compact metrisable group $G(\OOO_\val)$ acts continuously and
transitively by Lemma \ref{lem:maxcomptrans} on the compact metrisable
space $\Gr_{\d,\D}$.  The stabilizer in $G(\OOO_\val)$ of the
$K_\val$-linear subspace $K_\val^{\;\d}\times \{0\}$ of
$K_\val^{\;\D}$ corresponding to the first $\d$ coordinates is exactly
$Q^+(\OOO_\val)$. Therefore the continuous onto orbital map $g\mapsto
g(K_\val^{\;\d}\times\{0\})$ from $G(\OOO_\val)$ to $\Gr_{\d,\D}$
induces a continuous bijection
\[
G(\OOO_\val)/Q^+(\OOO_\val)\ra \Gr_{\d,\D}\;,
\]
which is hence a homeomorphism by compactness arguments. We identify
from now on $G(\OOO_\val)/Q^+(\OOO_\val)$ and $\Gr_{\d,\D}$ by this
map.

By the normalisation convention of the Haar measure of the closed
subgroups of $G$ (see Equation \eqref{eq:normalhaarsubgroup}), the
Haar measures $\mu_{G(\OOO_\val)}$ and $\mu_{Q^+\!(\OOO_\val)}$ are
normalized to be probability measures. We denote by
$\mu_{\Gr_{\d,\D}}$ the unique $G(\OOO_\val)$-invariant probability
%%
%\todo{\tiny I agree there is another rather natural normalization, but
%  let's try to be consistent}
%%
measure on the Grassmannian space $\Gr_{\d,\D}= G(\OOO_\val)/
Q^+\!(\OOO_\val)$. This is in accordance with Weil's normalization
process of measures on homogeneous spaces (see \cite[\S 9]{Weil65}).
Indeed, the probability measure $\mu_{G(\OOO_\val)}$ disintegrates
with respect to the canonical projection $G(\OOO_\val)\ra G(\OOO_\val)
/Q^+\!(\OOO_\val)$ over the measure $\mu_{\Gr_{\d,\D}}$, with
conditionnal measures on the fibers $g\,Q^+\!(\OOO_\val)$ the
probability pushforward measures $g_*\mu_{Q^+\!(\OOO_\val)}$: for
every $f\in C^0(G(\OOO_\val))$, we have
\begin{equation}\label{eq:disintegmugrass}
\int_{g\in G(\OOO_\val)} f(g)\;d\mu_{G(\OOO_\val)}
=\int_{g\,Q^+\!(\OOO_\val)\in \Gr_{\d,\D}} \int_{h\in Q^+\!(\OOO_\val)} f(gh)\;
d\mu_{Q^+\!(\OOO_\val)}(h)\;d\mu_{ \Gr_{\d,\D}}(g\,Q^+\!(\OOO_\val))\;.
\end{equation}
In particular, we have
\begin{equation}\label{eq:normalmesgrass}
\|\mu_{\Gr_{\d,\D}}\|=1\;.
\end{equation}

We denote by $\operatorname{orb}_\d:\M_{\dd,\d}(\OOO_\val)\ra
\Gr_{\d,\D}$ the (continuous injective) map defined by
\[
\operatorname{orb}_\d:\beta\mapsto
\big(\begin{smallmatrix} I_\d & 0 \\ \beta
  &I_\dd \end{smallmatrix}\big)(K_\val^{\;\d}\times\{0\})\;,
\]
and by $\Gr_{\d,\D}^{\flat}=\operatorname{orb}_\d(\M_{\dd,\d}
(\OOO_\val))$ its image in $\Gr_{\d,\D}=G(\OOO_\val)/Q^+(\OOO_\val)$.
Every element $g= \big(\begin{smallmatrix} \alpha &\ga \\ \beta
&\delta\end{smallmatrix}\big)\in G(\OOO_\val)$ such that
$\det(\alpha)\neq 0$ and $\beta\alpha^{-1}\in \M_{\dd,\d} (\OOO_\val)$
belongs to $\big(\begin{smallmatrix} I_\d &0\\ \beta\alpha^{-1}
&I_\dd\end{smallmatrix}\big) Q^+(\OOO_\val)$ by the first equality in
Equation \eqref{eq:decompgsimp}.  Conversely, if an element $g=\big(
\begin{smallmatrix} \alpha &\ga\\ \beta &\delta\end{smallmatrix}
\big)$ belongs to $U^-(\OOO_\val) \,Q^+(\OOO_\val)$, then we have
$\det\alpha\neq 0$ and $\beta\alpha^{-1}\in \M_{\dd,\d}(\OOO_\val)$,
so that
\begin{align*}
\Gr_{\d,\D}^{\flat}&=U^-(\OOO_\val)\, Q^+(\OOO_\val)=
U^-(\OOO_\val)(K_\val^{\;\d}\times\{0\})\\&=
\big\{\big(\begin{smallmatrix} \alpha &\ga
\\ \beta &\delta\end{smallmatrix}\big)\in G(\OOO_\val):
\substack{\det\alpha\neq 0\\\beta\alpha^{-1}\in \M_{\dd,\d}(\OOO_\val)}
\big\}(K_\val^{\;\d}\times\{0\})\;.
\end{align*}
Hence $\Gr_{\d,\D}^{\flat}$ is a compact (as the image by the
continuous map $\operatorname{orb}_\d$ of the compact space
$\M_{\dd,\d}(\OOO_\val)$) and open subset of the open Bruhat cell
$U^-Q^+$ of the Grassmannian space $\Gr_{\d,\D}=G/Q^+$ (corresponding
to the longest element in the Weyl group of $\SL_\D(K_\val)$).
By Equation \eqref{eq:disintegmugrass} with $f$ the
characteristic function of the compact subset $U^-(\OOO_\val)\,
Q^+(\OOO_\val)$ of $G(\OOO_\val)$ for the first equality, by
Proposition \ref{prop:refLUbloc} for the third equality and by the
normalisation in Equation \eqref{eq:normalhaarsubgroup} of the Haar
measures for the last equality, we have
\begin{align*}\mu_{\Gr_{\d,\D}}
  (\Gr_{\d,\D}^{\flat})&=\mu_{G(\OOO_\val)}(U^-(\OOO_\val)\,Q^+(\OOO_\val))
  =\mu_{G}(U^-(\OOO_\val)\,G''(\OOO_\val)\,U^+(\OOO_\val))
  \\&=c_1\;\mu_{U^-}(U^-(\OOO_\val))\,\mu_{G''}(G''(\OOO_\val))\,
  \mu_{U^+}(U^+(\OOO_\val))=c_1\;.
\end{align*}
By the normalization of the Haar measure $\operatorname{Haar}_{\dd,\d}$
of $\M_{\dd,\d}(K_\val)$ so that its restriction $\mu_{\M_{\dd,\d}
  (\OOO_\val)}$ to $\M_{\dd,\d}(\OOO_\val)$ is a probability measure,
we hence have
\begin{equation}\label{eq:relatHaarmnmuGr}
  (\operatorname{orb}_\d)_*(\mu_{\M_{\dd,\d}(\OOO_\val)})= c_1^{\;-1}
  \;(\mu_{\Gr_{\d,\D}})_{\mid\,\Gr_{\d,\D}^{\flat}} \;.
\end{equation}

In order to be able to define locally constant functions on the
Grassmannian space $\Gr_{\d,\D}$ for error term estimates, one way is
to define an appropriate distance on this space. The standard
construction is the following one. We endow the $K_\val$-vector space
$V=K_\val^{\;\D}$ with its usual norm, the maximum of the absolute
values of the coordinates in its canonical $K_\val$-basis $(e_1,\dots,
e_\D)$. For every $k\in\NN\ssm\{0\}$, we endow its $k$-th exterior
power $\text{\Large $\wedge$}\!^k V$ with the corresponding norm, the
maximum of the absolute values of the coordinates in its corresponding
$K_\val$-basis $(e_{i_1} \wedge \dots \wedge e_{i_k})_{1\leq i_1<\dots
  <i_k\leq \D}$. We now endow the projective space $\PP(V)$ with its
usual distance $d$ defined by $d(K_\val x, K_\val y)= \frac{\|x\wedge
  y\|} {\|x\|\;\|y\|}$ for all $x,y\in V\smallsetminus \{0\}$, and
$\Gr_{\d,\D}$ with its induced distance $d$ by the Plücker embedding
$\Gr_{\d,\D}\ra \PP(\text{\Large $\wedge$}\!^\d V)$ defined by
$W\mapsto K_\val (b_1 \wedge \dots \wedge b_\d)$ if $(b_1,\dots ,
b_\d)$ is any $K_\val$-basis of $W$. Since the linear action of
$\GL_\D(\OOO_\val)$ on $V$ preserves the norm, the exterior action of
$\GL_\D(\OOO_\val)$ on $\text{\Large $\wedge$}\!^\d V$ preserves the
norm, hence the projective action of $\GL_\D(\OOO_\val)$ on
$\PP(\text{\Large $\wedge$}\!^\d V)$ preserves the distance. Since the
Plücker embedding is equivariant with respect to the actions of
$\GL_\D(\OOO_\val)$, the action of $\GL_\D(\OOO_\val)$ on
$\Gr_{\d,\D}$ preserves its distance $d$.

For all $D,D',D''\in\NN\smallsetminus \{0\}$, we endow the
$K_\val$-vector space $\M_{D,D'}(K_\val)$ with its supremum norm
$\|\cdot\|$ defined, for every element $X=(X_{i,j})_{1\leq i\leq D,
  1\leq j\leq D'}\in \M_{D,D'}(K_\val)$, by
\[
\|X\|= \max\{|X_{i,j}|:1\leq i\leq D,1\leq j\leq D'\}
\in q_\val^{\;\ZZ} \cup \{0\}\;.
\]
This norm is an ultrametric norm and satisfies the following
properties:

$\bullet$~ The transposition map $A\mapsto \;^t\!A$ from
$\M_{D,D'}(K_\val)$ to $\M_{D',D}(K_\val)$ is a linear isometry for
the norms $\|\;\|$.

$\bullet$~ By the ultrametric property of the absolute value, the norm
$\|\;\|$ is a submultiplicative norm: For all $A\in\M_{D,D'}(K_\val)$
and $B\in\M_{D',D''}(K_\val)$, we have $\|AB\|\leq\|A\|\;\|B\|$.

$\bullet$~ For every $A\in\M_{D,D'}(K_\val)$, we have $\|A\|\leq 1$ if
and only if $A\in\M_{D,D'}(\OOO_\val)$. Hence the unit ball of
$(\M_{D,D'}(K_\val),\|\cdot\|)$ is $\M_{D,D'}(\OOO_\val)$ and the
right and left multiplications by elements of $\M_D(\OOO_\val)$ and
$\M_{D'}(\OOO_\val)$ are $1$-Lipschitz maps on
$\M_{D,D'}(\OOO_\val)$~: For all $A\in\M_D(K_\val)$ and $B\in\M_{D,D'}
(\OOO_\val)$ and $C\in \M_{D'}(\OOO_\val)$, we have $\|ABC\| \leq
\|B\|$. In particular, $\|ABC\|= \|B\|$ if $A\in\GL_D(\OOO_\val)$ and
$C\in\GL_{D'} (\OOO_\val)$.

\blemm \label{lem:isometrygrass}
For all $\beta,\beta'\in \M_{\dd,\d}(\OOO_\val)$, we have
\[
d(\operatorname{orb}_\d(\beta),\operatorname{orb}_\d(\beta'))=
\|\beta-\beta'\|\;.
\]
\elemm

\dem Every matrix $\beta\in\M_{\dd,\d}(\OOO_\val)$ will be seen as a
linear map from $K_\val^{\;\d}\times\{0\}$ to $\{0\}\times
K_\val^{\;\dd}$, so that $\beta e_j=\sum_{1\leq i\leq \dd} \beta_{i,j}
\, e_{i+\d}$ for every $j\in\llbracket1, \d\rrbracket$. A
$K_\val$-basis of the $K_\val$-linear subspace $\operatorname{orb}_\d
(\beta) =\big(\begin{smallmatrix} I_\d & 0\\ \beta
  &I_\dd \end{smallmatrix}\big)(K_\val^{\;\d}\times\{0\})$ of
$K_\val^{\;\D}$ is hence $(e_1+\beta e_1,\dots ,e_\d+\beta e_\d)$.

Let $x_\beta=(e_1+\beta e_1)\wedge\dots \wedge(e_\d+\beta e_\d)$. We
have $\|x_\beta\|=1$ since $\|e_1\wedge\dots \wedge e_\d\|=1$ and
since the entries of $\beta$ have absolute value at most $1$.  For
every $\beta,\beta'\in \M_{\dd,\d}(\OOO_\val)$, we have $x_\beta\wedge
x_{\beta'}=(x_\beta-x_{\beta'})\wedge x_{\beta'}$ and
$x_\beta-x_{\beta'}=v_1 +v_2+\dots +v_\d$ where, separating the terms
according to the number of occurrences of $\beta$'s in them,
\[
v_1=\sum_{1\leq i\leq \d} e_1\wedge\dots \wedge
e_{i-1} \wedge(\beta-\beta')e_i \wedge e_{i+1} \wedge \dots \wedge
e_\d\,,
\]
\begin{align*}
v_2=\sum_{1\leq i<j\leq \d} \;&\big(e_1\wedge\dots \wedge e_{i-1} \wedge
\beta e_i \wedge e_{i+1}\wedge\dots\wedge e_{j-1}\wedge \beta e_j
\wedge e_{j+1} \wedge \dots \wedge e_\d\\ &-e_1\wedge\dots \wedge e_{i-1}
\wedge \beta' e_i \wedge e_{i+1}\wedge\dots\wedge e_{j-1}\wedge \beta' e_j
\wedge e_{j+1} \wedge \dots \wedge
e_\d\big)
\end{align*}
and so on, and $v_\d=\beta e_1\wedge\dots \wedge\beta e_\d-\beta'
e_1\wedge\dots \wedge\beta' e_\d$. By the ultrametric properties of
the norms, we have $\|v_k\|\leq \|\beta-\beta'\|$ for all $k\geq 2$
and $\|v_1\|= \|\beta-\beta'\|$. By considering the elements of the
$K_\val$-basis of $\text{\Large $\wedge$}\!^\d V$ involved in the
formulation of $v_1,v_2,\dots, v_\d$, and again by the ultrametric
properties of the norms, we hence have $\|x_\beta-x_{\beta'}\|=
\|\beta-\beta'\|$.  Furthermore, the coordinate of
$x_\beta-x_{\beta'}$ corresponding to the basis vector $e_1\wedge\dots
\wedge e_\d$ is $0$ while the one of $x_{\beta'}$ is $1$.  Thus
$\|x_\beta\wedge x_{\beta'}\|= \|x_\beta-x_{\beta'}\|$ and the result
follows.
\cqfd

\subsection{The spaces of unimodular full lattices}
\label{subsect:fulllat}

For every unimodular locally compact group $H$ endowed with a Haar
measure $\mu_H$, and for every discrete subgroup $\Ga'$ of $H$, we
again denote by $\mu_H$ the unique left $H$-invariant measure on
$H/\Ga'$ such that the covering map $H\ra H/\Ga'$ locally
preserves the measure.

Let $k\in\NN\ssm\{0\}$. Note that $\SL_k(K_\val)$ is a closed
unimodular subgroup of the unimodular locally compact group
$\GL^1_k(K_\val)$, whose Haar measure $\mu_{\SL_k(K_\val)}$ is
normalized so that $\mu_{\SL_k(K_\val)} (\SL_k(\OOO_\val)) =1$ (as we
did for $k=\D$ in Subsection \ref{sect:LUbloc}).  The restriction to
$\OOO_\val^{\;\times}$ of the Haar measure of the additive group
$(K_\val,+)$ is a Haar measure $\mu_{\OOO_\val^{\;\times}}$ of the
multiplicative group $(\OOO_\val^{\;\times}, \times)$ by Equation
\eqref{eq:homothetyhaar}.  By the normalization of the Haar measure of
$(K_\val,+)$, we have
\begin{equation}\label{eq:haarOvtimes}
\|\mu_{\OOO_\val^{\;\times}}\|= \mu_{K_\val} (\OOO_\val\ssm
\pi_\val\OOO_\val) =1-q_\val^{\;-1}\;.
\end{equation}
We have a split short exact sequence of locally compact groups
\[
1\longrightarrow \SL_k(K_\val)\longrightarrow \GL^1_k(K_\val)
\longrightarrow \OOO_\val^{\;\times}\;,
\]
with section $s_k:\OOO_\val^{\;\times}\ra\GL^1_k(K_\val)$ defined for
instance by $\lambda\mapsto \big(\begin{smallmatrix}\lambda & 0 \\ 0 &
  I_{k-1}\end{smallmatrix}\big)$. We define the Haar measure
$\mu_{\GL^1_k(K_\val)}$ of $\GL^1_k(K_\val)$, for all $g\in
\SL_k(K_\val)$ and $\lambda\in\OOO_\val^{\;\times}$, by
\begin{equation}\label{eq:defihaarGL1}
d\mu_{\GL^1_k(K_\val)}(s_k(\lambda)g)=\frac{q_\val}{q_\val-1}
d\mu_{\OOO_\val^{\;\times}}(\lambda)\;d\mu_{\SL_k(K_\val)}(g)\;.
\end{equation}
In particular, by Equations \eqref{eq:defihaarGL1} and
\eqref{eq:haarOvtimes}, accordingly to the convention in Equation
\eqref{eq:normalhaarsubgroup} when $k=\D$, we have
\begin{equation}\label{eq:normalmuGL1}
    \mu_{\GL^1_k(K_\val)} (\GL_k(\OOO_\val))=1\;.
\end{equation}

After Equation \eqref{eq:defiGL1}, we identified the space
$\La_k=\La(K_\val^{\;k})$ of unimodular full $R_\val$-lattices in
$K_\val^{\;k}$ with the homogeneous space $\GL^1_k(K_\val)/
\GL_k(R_\val)$.  Since $\GL_k(R_\val)$ is a discrete subgroup of the
unimodular group $\GL_k^1(K_\val)$, as said at the beginning of this
subsection \ref{subsect:fulllat}, we endow $\La_k$ with the unique
$\GL_k^1(K_\val)$-invariant measure $\mu_{\La_k}$ such that the
orbital map $\GL_k^1(K_\val)\ra\La_k$ defined by $g\mapsto
g\,R_\val^{\;k}$ locally preserves the measure.

Since the index of $\SL_k(R_\val)$ in $\GL_k(R_\val)$ is equal to
$\card(R_\val^{\;\times})=\card(\FF_q^{\,\times})=q-1$ (see the
formula on the right in Equation \eqref{eq:inversiRv}), and by Equations
\eqref{eq:defihaarGL1} and \eqref{eq:haarOvtimes}, we have
\begin{align}
  \|\mu_{\La_k}\|&=\|\mu_{\GL^1_k(K_\val)/\GL_k(R_\val)}\|=
  \frac{1}{q-1}\;\|\mu_{\GL^1_k(K_\val)/\SL_k(R_\val)}\|\nonumber\\&
  =\frac{1}{q-1}\;\|\mu_{\SL_k(K_\val)/\SL_k(R_\val)}\|\;.
  \label{eq:premreducmulat}
\end{align}

Let us apply \cite[\S 3]{Serre71} in order to compute the total mass
of $\mu_{\La_k}$, using boldface letters in order to denote the
notation of this reference, thus facilitating the reference
process. Let $\bm{L}=\SL_k$, which is a simple simply connected split
algebraic group defined over the global field $\bm{k}=K$, with
relative rank $\bm{\ell}=k-1$, and with exponents of its Weyl group
$\bm{m}_i=i$ for $i\in\llbracket 1,\bm{\ell}\rrbracket$ (see
\cite[page 251]{Bourbaki81}). Let $\bm{S}=\{\bm{v}=\val\}$, which is a
finite nonempty set of places of $\bm{k}$, and note that there are no
archimedean places since $K$ is a function field. The function ring
$\bm{O}_{\bm S}$ defined in \cite[page 123]{Serre71} is then exactly
our function ring $R_\val$, and $\bm{L}_{\bm{0}}=\SL_k$ is a split,
simple, simply connected group scheme over $\bm{O}_{\bm S}=R_\val$ (as
required in \cite[page 157]{Serre71}) such that $\bm{L}=
\bm{L}_{\bm{0}}\otimes_{\bm{O}_{\bm{S}}}\bm{k}$.  The zeta function
$\bm{\zeta}_{\bm{k},\bm{S}}$ of $\bm{k}$ relatively to $\bm{S}$ defined
in \cite[page 156]{Serre71} is exactly our zeta function $\zeta_K$,
and the $\bm{S}$-arithmetic group $\bm{\Ga}_{\bm{S}}$ defined in
\cite[page 157]{Serre71} is exactly our arithmetic group
$\SL_k(R_\val)$.

Let $\bm{G}= \bm{L}(\bm{k}_{\bm{v}})=\SL_k(K_\val)$.  Motivated by the
relationship with the Euler characteristic, Serre defines a canonical
signed measure (with constant sign by homogeneity) $\bm{\mu}_{\bm{G}}$
on $\bm{G}$, whose associated positive measure $|\bm{\mu}_{\bm{G}}|$
is a Haar measure on $\bm{G}$.  We don't need to recall its
definition, only to understand its normalisation. By the second claim
of Theorem 7 (see top of page 151) of \cite{Serre71}, using when $k=1$
the standard convention that an empty product is equal to $1$, and
since the order of the residual field of $\bm{k}_{\bm v}=K_\val$ is
$\bm{q}=q_\val$, we have
\[
|\bm{\mu}_{\bm{G}}|(\SL_k(\OOO_\val))=
\prod_{i=1}^{\bm{\ell}}(\bm{q}^{\bm{m}_i}-1)=
\prod_{i=1}^{k-1}(q_\val^{\;i}-1)\;.
\]
Since our Haar measure of $\SL_k(K_\val)$ is normalized so that
$\mu_{\SL_k(K_\val)}(\SL_k(\OOO_\val))=1$, we hence have
\begin{equation}\label{eq:reducmuSL}
\mu_{\SL_k(K_\val)} =\frac{1}{|\bm{\mu}_{\bm{G}}|(\SL_k(\OOO_\val))}
|\bm{\mu}_{\bm{G}}|=\prod_{i=1}^{k-1}(q_\val^{\;i}-1)^{-1}|\bm{\mu}_{\bm{G}}|
\;.
\end{equation}

By, for instance, \cite{Langlands66b} (see also \cite[Theo
  3.3.1]{Weil61} and \cite[p. 257]{Weil60}),
%%
%\todo{\tiny ou A.~Weil Adeles and algebraic groups, Lecture
%  notes. Princeton, 1961, Theo 3.3.1 + page 257 de Adèles et groupes
%  algébriques, Séminaire N. Bourbaki, 1960, exp. no 186, p. 249-257}
%%
the Tamagawa number $\bm{\tau}$ of $\bm{L}$ is $1$. By the footnote 10
on page 158 of \cite{Serre71}, and since $\zeta_K(-i)>0$ as seen in
Equation \eqref{eq:poszetanegent}, we hence have
\begin{equation}\label{eq:tamagawa}
|\bm{\mu}_{\bm{G}}|(\bm{G}/\bm{\Ga}_{\bm{S}})=
\Big|\bm{\tau}\prod_{i=1}^{\bm{\ell}}\bm{\zeta}_{\bm{k},\bm{S}}(-\bm{m}_i)\Big|=
\prod_{i=1}^{k-1}\zeta_{K}(-i)\;.
\end{equation}
Thus, by Equations \eqref{eq:reducmuSL}, \eqref{eq:tamagawa} and
\eqref{eq:poszetanegent}, we have
\begin{equation}\label{eq:totmassmuSLmodSL}
  \|\mu_{\SL_k(K_\val)/\SL_k(R_\val)}\|=
  \prod_{i=1}^{k-1}\frac{\zeta_K(-i)}{q_\val^{\;i}-1}
  =q^{\;(\ggg-1)(k^2-1)}
  \prod_{i=1}^{k-1}\frac{\zeta_K(1+i)}{q_\val^{\;i}-1}\;.
\end{equation}
Therefore, by Equation \eqref{eq:premreducmulat}, we have
\begin{equation}\label{eq:totmassmulat}
\|\mu_{\La_k}\|
=\frac{1}{q-1}
\prod_{i=1}^{k-1}\frac{\zeta_K(-i)}{q_\val^{\;i}-1}\;.
\end{equation}

\subsection{The space of correlated pairs of lattices}
\label{subsect:latticepair}

We define in this subsection a measured space $\La_{\d,\dd}$ of pairs
of full $R_\val$-lattices in dimensions $\d$ and $\dd$ with correlated
normalized covolume, in which a version of the equidistribution
results stronger than the ones stated in the introduction will take
place (see Theorem \ref{theo:mainsharp}).

For every $k\in\NN\smallsetminus\{0\}$, let $(e_1,\dots,e_k)$ be the
canonical $K_\val$-basis of $K_\val^{\;k}$. Let $\Lambda\in \La_k$.
For every $R_\val$-basis $(b_1,\dots,b_k)$ of $\Lambda$, which is also
a $K_\val$-basis of $K_\val^{\;k}$, let $g\in\GL_k(K_\val)$ be the unique
$K_\val$-linear map sending $(e_1,\dots,e_k)$ to $(b_1,\dots,b_k)$.
Note that $\det g\in\OOO_\val^{\;\times}$ since $|\det g|=
\overline{\covol}(\Lambda) =1$. The image $(\det g)\,\FF_q^{\;\times}$
of $\det g$ by the canonical projection $\OOO_\val^{\;\times}\ra
\OOO_\val^{\;\times}/\FF_q^{\;\times}$ does not depend on the choice
of the $R_\val$-basis of $\Lambda$, since the change of $R_\val$-basis
matrix belongs to $\GL_k(R_\val)$, hence has determinant in
$R_\val^{\;\times}=\FF_q^{\;\times}$. We define
\[
\det \Lambda=
(\det g)\,\FF_q^{\;\times}\in \OOO_\val^{\;\times}/\FF_q^{\;\times}\,.
\]
Note that the ratio $\frac{a}{a'}$ of two elements $a=\lambda
\FF_q^{\;\times}$ and $a'=\lambda'\FF_q^{\;\times}$ of
$\OOO_\val^{\;\times}/\FF_q^{\;\times}$ is a well-defined element
$\frac{a}{a'}=\frac{\lambda}{\lambda'}\, \FF_q^{\;\times}$ of
$\OOO_\val^{\;\times}/\,\FF_q^{\;\times}$, and we will as usual also
denote by $\FF_q^{\;\times}$ the class of $1$ in $\OOO_\val^{\;\times}
/\,\FF_q^{\;\times}$.  Let us define
\[
\La_{\d,\dd}=\Big\{(\Lambda,\Lambda')\in \La_{\d}\times\La_{\dd}:
\frac{\det \Lambda}{\det\Lambda'}= \FF_q^\times\Big\}\;.
\]
For instance, $(R_\val^{\;\d}, R_\val^{\;\dd})$ belongs to
$\La_{\d,\dd}$.

Recall that the map $\delta\mapsto\widecheck{\delta}=\;^t\delta^{-1}$
for every $\delta\in \GL_\dd^1(K_\val)$ is the standard Cartan
involution of $\GL_\dd^1(K_\val)$. The product group $\GL^1_\d(K_\val)
\times \GL^1_\dd(K_\val)$, embedded diagonally by blocs in
$\GL_\D(K_\val)$, acts continuously on the product space
$\La_{\d}\times\La_{\dd}$ by the twisted diagonal action
\begin{equation}\label{eq:actGpplatmn}
\big(\big(\begin{smallmatrix} \alpha&0 \\ 0 &
  \delta\end{smallmatrix}\big),(\Lambda,\Lambda')\big)\mapsto
(\,\alpha\,\Lambda,\widecheck{\delta}\,\Lambda')\;,
\end{equation}
since $\GL_k^1(K_\val)$ preserves the set of unimodular lattices
$\La_k$ for every $k\in\NN\smallsetminus\{0\}$.  Recall that
$G''=\big\{\big(\begin{smallmatrix} \alpha&0 \\ 0 &
  \delta\end{smallmatrix}\big):
  \alpha\in\GL^1_\d(K_\val),\;\delta\in\GL^1_\dd(K_\val),\; \det
  \alpha\;\det\delta=1\big\}$.

\blemm \label{lem:latmnhomogen} The orbit in $\La_{\d}\times\La_{\dd}$
of $(R_\val^{\;\d}, R_\val^{\;\dd})$ by the action defined in Equation
\eqref{eq:actGpplatmn} of the subgroup $G''$ of $\GL^1_\d(K_\val)
\times\GL^1_\dd(K_\val)$ is exactly $\La_{\d,\dd}$.
%, and the restriction to $\La_{\d,\dd}$ of the action of $G''$
%preserves the measure $\mu_{\La_{\d,\dd}}$.
\elemm

\dem For all $k\in\NN\smallsetminus\{0\}$, $g\in\GL_k^1(K_\val)$ and
$\Lambda \in \La_k$, we have $\det(g\Lambda)=\det (g)\det(\Lambda)$.
Since $\det(\,\widecheck{\underline{g}}\,) =(\det\underline{g}
\,)^{-1}$, the action of $G''$ defined in Equation
\eqref{eq:actGpplatmn} preserves $\La_{\d,\dd}$.

Conversely, let $(\Lambda,\Lambda')\in\La_{\d,\dd}$. Since the action
of $\GL_k^1(K_\val)$ on $\La_k$ is transitive for every $k\in\NN
\smallsetminus\{0\}$ and since the map $g\mapsto \widecheck{g}$ is an
automorphism of $\GL_k^1(K_\val)$, there exist $\alpha\in
\GL_\d^1(K_\val)$ and $\delta\in \GL_\dd^1(K_\val)$ such that
$\Lambda=\alpha\,R_\val^{\;\d}$ and $\Lambda'= \widecheck{\delta}
\,R_\val^{\;\dd}$. Since $(\Lambda, \Lambda') \in\La_{\d,\dd}$, the
product $\lambda= \det\alpha \; \det\delta=\frac{\det\alpha}
{\det\widecheck{\delta}}$ belongs to $\FF_q^{\;\times}=
R_\val^{\;\times}$. Hence the $\d\times\d$ diagonal matrix
$s_\d(\lambda)$ with diagonal entries $\lambda,1,\dots, 1$ belongs to
$\GL_\d(R_\val)$, and in particular we have $s_\d(\lambda)^{-1}
R_\val^{\;\d} =R_\val^{\;\d}$. Then the matrix
$\big(\begin{smallmatrix} \alpha\,s_\d(\lambda)^{-1}&0 \\ 0 &
\delta\end{smallmatrix} \big)$ belongs to $G''$ and maps
$(R_\val^{\;\d},R_\val^{\;\dd})$ to $(\Lambda,\Lambda')$ by the action
defined in Equation \eqref{eq:actGpplatmn}. This proves the lemma.
\cqfd

\medskip
We denote by $\pi_{\d,\dd}: G''\ra\La_{\d,\dd}$ the twisted canonical
projection defined by
\begin{equation}\label{eq:defiprok}
\pi_{\d,\dd}:\big(\begin{smallmatrix} \alpha &0 \\ 0 &
  \delta \end{smallmatrix}\big)\mapsto
(\alpha R_\val^{\;\d},\widecheck{\delta} R_\val^{\;\dd})\;.
\end{equation}
Since the stabilizer of $(R_\val^{\;\d},R_\val^{\;\dd})$ by the action
of $G''$ defined in Equation \eqref{eq:actGpplatmn} is exactly the
discrete subgroup $G''(R_\val)$, and since $\pi_{\d,\dd}$ is the
orbital map of the above action at $(R_\val^{\;\d}, R_\val^{\;\dd})$,
we will identify from now on the quotient space $G''/G''(R_\val)$ and
the subspace $\La_{\d,\dd}$ by the $G''$-equivariant homeomorphism
$g\,G''(R_\val)\mapsto \pi_{\d,\dd}(g)= g(R_\val^{\;\d},
R_\val^{\;\dd})$.

Recall that $G''$ is endowed with its Haar measure $\mu_{G''}$
normalized so that we have $\mu_{G''}(G''(\OOO_\val))=1$. As said at the
beginning of Subsection \ref{subsect:fulllat}, we endow $G''/
G''(R_\val)$ with the unique $G''$-invariant measure such that the
covering map $G''\ra G''/G''(R_\val)$ locally preserves the measure.
We denote by $\mu_{\La_{\d,\dd}}$ the corresponding measure on
$\La_{\d,\dd}$, which is $G''$-invariant.

\blemm \label{lem:masstotLamn} We have
$\displaystyle\|\mu_{\La_{\d,\dd}}\|=\frac{1}{q-1}\;
\prod_{i=1}^{\d-1}\frac{\zeta_K(-i)}{q_\val^{\;i}-1}\;
\prod_{i=1}^{\dd-1}\frac{\zeta_K(-i)}{q_\val^{\;i}-1}$.
\elemm

\dem Note that $\GL_\dd(R_\val)$, $\GL_\dd(\OOO_\val)$ and
$\GL_\dd^1(K_\val)$ are stable by the Cartan involution
$\underline{g}\mapsto \widecheck{\underline{g}}=
\;^t\underline{g}^{-1}$, and that this map preserves the Haar measure
$\mu_{\GL_\dd^1(K_\val)}$ (defined in Subsection
\ref{subsect:fulllat}) of the selfadjoint unimodular group
$\GL_\dd^1(K_\val)$. For every $\lambda\in\OOO_\val^{\;\times}$, let
\[
s'_\d(\lambda)
=\big(\big(\begin{smallmatrix} \lambda & 0 \\ 0 &
  I_{\d-1} \end{smallmatrix}\big),I_{\dd}\big) \in\GL_\d(\OOO_\val)
\times \GL_\dd(\OOO_\val)\;.
\]
Let $\iota:G''\ra \GL^1_\d(K_\val) \times \GL^1_\dd(K_\val)$ be the
group monomorphism defined by $\big(\begin{smallmatrix} \overline{g}&
  0 \\ 0 & \underline{g}\end{smallmatrix}\big)\mapsto (\,\overline{g},
\widecheck{\underline{g}}\,)$. Let $\text{\small\textcursive{p}}:
\GL^1_\d(K_\val) \times \GL^1_\dd(K_\val)\ra\OOO_\val^{\;\times}$ be
the group epimorphism defined by $(\alpha,\delta)\mapsto
\frac{\det\alpha}{\det\delta}$. Any element of $\GL^1_\d(K_\val)
\times \GL^1_\dd(K_\val)$ may be written as $s'_\d(\lambda)\,
\iota(g'')$ for unique elements $\lambda\in \OOO_\val^{\;\times}$ and
$g''\in G''$.  Hence we have a split exact sequence
\[
1\longrightarrow G''\stackrel{\iota}{\longrightarrow}\GL^1_\d(K_\val)
\times \GL^1_\dd(K_\val)
\stackrel{\text{\tiny\textcursive{p}}}{\longrightarrow}
\OOO_\val^{\;\times}\longrightarrow1\;,
\]
which induces two split exact sequences
\[
1\longrightarrow G''(\OOO_\val)\stackrel{\iota}{\longrightarrow}
\GL_\d(\OOO_\val)\times \GL_\dd(\OOO_\val)
\stackrel{\text{\tiny\textcursive{p}}}{\longrightarrow}
\OOO_\val^{\;\times}\longrightarrow1\;,
\]
\[
1\longrightarrow G''(R_\val)\stackrel{\iota}{\longrightarrow}
\GL_\d(R_\val)\times \GL_\dd(R_\val)
\stackrel{\text{\tiny\textcursive{p}}}{\longrightarrow}
R_\val^{\;\times}\longrightarrow1\;.
\]
Since the group $s'_\d(\OOO_\val^{\;\times})$ normalizes $\iota(G'')$, the
measure $\mu$ on the product topological group $\GL^1_\d(K_\val) \times
\GL^1_\dd(K_\val)$ defined by $d\mu(s'_\d(\lambda)\iota(g''))=
d\mu_{G''}(g'')\; d\mu_{\OOO_\val^{\;\times}}(\lambda)$ is a Haar
measure. By the normalisation of the measures (see Equations
\eqref{eq:normalmuGL1}, \eqref{eq:normalhaarsubgroup} and
\eqref{eq:haarOvtimes}), the product Haar measure
$\mu_{\GL_\d^1(K_\val)} \otimes \mu_{\GL_\dd^1(K_\val)}$ of
$\GL^1_\d(K_\val) \times \GL^1_\dd(K_\val)$ hence satisfies
\begin{equation}\label{eq:GLGLGL}
d(\mu_{\GL_\d^1(K_\val)}\otimes\mu_{\GL_\dd^1(K_\val)})(s'_\d(\lambda)\iota(g''))
=(1-q_\val^{\;-1})^{-1}\;d\mu_{G''}(g'')\;d\mu_{\OOO_\val^{\;\times}}(\lambda)\;.
\end{equation}
Since we have $\card \;R_\val^{\;\times}=q-1$ by Equation
\eqref{eq:inversiRv} and again by Equation \eqref{eq:haarOvtimes}, we
hence have
\begin{align}\|\mu_{\La_{\d}}\|\;\|\mu_{\La_{\dd}}\|=
\|\mu_{\La_{\d}}\otimes\mu_{\La_{\dd}}\|&=
(1-q_\val^{\;-1})^{-1}\|\mu_{G''/G''(R_\val)}\|\;
\frac{\|\mu_{\OOO_\val^{\;\times}}\|}{\card\;R_\val^{\;\times}}\nonumber
\\&=\frac{1}{q-1}\;\|\mu_{\La_{\d,\dd}}\|\,.\label{eq:muLatnLatmLamn}
\end{align}
The result follows by applying twice Equation \eqref{eq:totmassmulat},
with $k=\d$ and $k=\dd$.
\cqfd

\medskip
As for the Grassmannian spaces $\Gr_{\d,\D}$, in order to be able to
define locally constant functions on $\La_{\d,\dd}$ for our error term
estimates, we now define a natural distance on the space
$\La_{\d,\dd}$.

Let $k\in\NN\!\smallsetminus\!\{0\}$. Since the supremum norm $\|\;\|$
on $\M_k(K_\val)$ is a submultiplicative norm (see above Lemma
\ref{lem:isometrygrass}), the map
\begin{equation}\label{eq:defidistGLk}
d:(g,h)\mapsto \log_{q_\val}(1+\max\{\|\, gh^{-1}-I_k\,\|,
\;\|\, hg^{-1}-I_k\,\|\})
\end{equation}
is well-known to be a distance on the locally compact group
$\GL_k(K_\val)$ (inducing its topology).  By construction, this
distance is invariant by translations on the right by all elements of
$\GL_k(K_\val)$. It is also invariant by translations on the left by
the elements of $\GL_k(\OOO_\val)$, since the supremum norm $\|\;\|$
on $\M_k(K_\val)$ is invariant under conjugation by any element of
$\GL_k(\OOO_\val)$. Since the transposition map $g\mapsto \,^tg$
preserves the supremum norm $\|\;\|$ on $\M_k(K_\val)$ and by the
symmetry of the distance $d$ on $\GL_k(K_\val)$, the map $g\mapsto
\widecheck{g}=\,^tg^{-1}$ is an isometry of $d$. In particular, for
all $g\in \GL_k(K_\val)$ and $\rho>0$, we have
\[
\widecheck{B(g,\rho)}=B(\,\widecheck{g},\rho)\;.
\]
The following lemma will be needed in Subsection
\ref{subsect:wellround}.

\blemm \label{lem:controlnormsupdist}
For all $h,h_0\in\GL_k(K_\val)$, we have
\[
\|h\|\leq \|h_0\|\;q_\val^{\;d(h,\,h_0)}\quad\text{and}\quad
\|h^{-1}\|\leq\|h_0^{\;-1}\|\;q_\val^{\;d(h,\,h_0)}\;.
\]
\elemm

\dem Let $t= d(h,\,h_0)$. We have $\log_{q_\val}(1+\|\, hh_0^{-1}
-I_k\,\|)\leq t$, hence $\|\,hh_0^{-1}-I_k\,\|\leq q_\val^{\;t}-1$, thus
\[
\|h\|-\|h_0\|\leq \|h-h_0\|=\|(hh_0^{-1}-I_k)h_0\|\leq
\|hh_0^{-1}-I_k\|\;\|h_0\| \leq (q_\val^{\;t}-1)\|h_0\|\;.
\]
The first result follows. The second result follows similarly (or
since the map $g\mapsto \widecheck{g}$ is an isometry and the
transposition preserves the supremum norm of matrices).
\cqfd

\medskip
We endow every closed subgroup $H$ of $\GL_k(K_\val)$ with the induced
distance, and, defining $H(R_\val)=H\cap \GL_k(R_\val)$ which is a
discrete subgroup of $H$, we endow the quotient space $H/H(R_\val)$
with the quotient distance
\[
\forall\;g,h\in H,\quad d(g\,H(R_\val),h\,H(R_\val))=
\min_{\ga\in H(R_\val)} d(g,h\ga)\;.
\]
It is easy to check by Equation \eqref{eq:defidistGLk} that for all
$g=\big(\begin{smallmatrix} \alpha &0 \\ 0 &
  \delta\end{smallmatrix}\big),g'=\big(\begin{smallmatrix} \alpha' &0
    \\ 0 & \delta'\end{smallmatrix}\big)\in G''$, we have
\begin{equation}\label{eq:diagodist}
  d_{\GL_\D(K_\val)}(g,g')=\max\{d_{\GL_\d(K_\val)}(\alpha,\alpha'),
  d_{\GL_\dd(K_\val)}(\delta,\delta')\}\;.
\end{equation}

The canonical projection $H\ra H/H(R_\val)$ is $1$-Lipschitz and is a
local isometry since the discrete subgroup $H(R_\val)$ of $H$ acts
isometrically by right-translations on $H$. The action of
$H(\OOO_\val)=H\cap \GL_k(\OOO_\val)$ by translations on the left on
$H/H(R_\val)$ is isometric. This process provides the homogeneous
spaces $\La_{\d}=\GL_\d^1(K_\val)/\GL_\d(R_\val)$, $\La_{\dd}=
\GL_\dd^1(K_\val) /\GL_\dd(R_\val)$ and $\La_{\d,\dd}=G''/G''(R_\val)$
with distances invariant under $\GL_\d(\OOO_\val)$,
$\GL_\dd(\OOO_\val)$ and $G''(\OOO_\val)$ respectively, that from now
on we will consider on these spaces. In particular, the map
$\pi_{\d,\dd}$ from $G''$ to $\La_{\d,\dd}$ defined in Equation
\eqref{eq:defiprok} is a local isometry.

\subsection{The spaces of shapes of unimodular full lattices}
\label{subsect:shape}

Let $k\in\NN\ssm\{0\}$. As defined in Subsection
\ref{subsect:LUdecpartlatt}, the space of shapes of full unimodular
$R_\val$-lattices of $K_\val^{\;k}$ is the locally compact metrisable
separable (actually discrete and countably infinite) quotient space
\[
\Sh_k=\GL_k(\OOO_\val)\bs\La_k
=\GL_k(\OOO_\val)\bs\GL^1_k(K_\val)\,/\GL_k(R_\val)\;.
\]
We endow $\Sh_k$ with the unique finite measure $\mu_{\Sh_k}$ such
that the left invariant finite measure $\mu_{\La_k}$ on the right
homogeneous space $\La_k=\GL^1_k(K_\val)/\GL_k(R_\val)$ disintegrates
with respect to the proper canonical projection $\sh:\La_k\ra
\Sh_k=\GL_k(\OOO_\val)\bs\La_k$ over $\mu_{\Sh_k}$ with conditional
measures on the fibers $\GL_k(\OOO_\val) \Lambda$ the pushforward
measures of the finite Haar measure $\mu_{\GL_k(\OOO_k)}=
{\mu_{\GL^1_k(K_\val)}}_{\mid\GL_k(\OOO_\val))}$ of $\GL_k(\OOO_\val)$
by the orbital maps $g\mapsto g\Lambda$~: for every $f\in C^0(\La_k)$
with compact support, we have
\begin{align*}
&\int_{\Lambda\in\La_k}f(\Lambda)\;d\mu_{\La_k}(\Lambda)\nonumber
\\=\;&\int_{\GL_k(\OOO_\val)\Lambda\in \Sh_k}\int_{g\in\GL_k(\OOO_\val)} f(g\Lambda)
\;d\mu_{\GL_k(\OOO_\val)}(g)\;d\mu_{\Sh_k}(\GL_k(\OOO_\val)\Lambda)
\label{eq:disintegmuLat}\;.
\end{align*}
In particular, using Equation \eqref{eq:normalmuGL1},
we have
\begin{equation}\label{eq:relatmuShmuLat}
  \sh_*\mu_{\La_k}=\|\mu_{\GL_k(\OOO_\val)}\|\;\mu_{\Sh_k}
  =\mu_{\Sh_k}\;,
\end{equation}
and by  Equation \eqref{eq:totmassmulat}, we have
\[
\|\mu_{\Sh_k}\|=\|\mu_{\La_k}\|
  =\frac{1}{q-1}\prod_{i=1}^{k-1}\frac{\zeta_K(-i)}{q_\val^{\;i}-1}\;.
\]

With the notation of Equation \eqref{eq:defish} (that greatly
simplifies as indicated below it), we define a map
\begin{equation}\label{eq:defivarphimn}
\begin{array}{llll}\varphi_{\d,\dd}:&\La_{\d,\dd}&\ra &
  \Sh_\d\times\Sh_\dd=\big(\GL_\d(\OOO_\val)\bs\La_{\d}\big)
  \times\big(\GL_\dd(\OOO_\val)\bs\La_{\dd}\big)\\ &
  (\Lambda,\Lambda')&\mapsto& (\sh(\Lambda),\sh(\Lambda'))=
  (\GL_\d(\OOO_\val)\Lambda,\; \GL_\dd(\OOO_\val)\Lambda')\;.
\end{array}
\end{equation}
We summarize its properties in the following lemma, after giving some
notation.

Since $G''(\OOO_\val)$ is compact and open in $G''$, there exists a
maximal $\rho_0>0$ such that $G''(\OOO_\val)$ contains the (closed)
ball $B_{G''}(I_\D,\rho_0)$ of center $I_\D$ and radius $\rho_0$ for the
distance on $G''$ defined at the end of the previous subsection
\ref{subsect:latticepair}. Thus every map $f:\La_{\d,\dd}\ra\CC$ which
is constant on every left $G''(\OOO_\val)$-orbit in $\La_{\d,\dd}$ is
$\rho_0$-locally constant, that is constant on every ball of radius
$\rho_0$ in $\La_{\d,\dd}$.
%and such that the map $\pi_{\d,\dd}:G''\ra
%\La_{\d,\dd}$ is isometric on $B_{G''}(\id,\rho_0)$.

For every $g\in \U_G$, if $g=u^-\,g''\,z\,u^+$ with
$g''=\big(\begin{smallmatrix} \overline{g} &0 \\ 0 &\underline{g}
\end{smallmatrix}\big)$ is its unique writing given by Proposition
\ref{prop:refLUbloc}, using the map $\pi_{\d,\dd}$ defined in Equation
\eqref{eq:defiprok}, we define the {\it correlated pair of lattices
  $\llbracket \Lambda_g\rrbracket$ associated with} $g$ by
\begin{equation}\label{eq:assoccorrpairlat}
  \llbracket \Lambda_g\rrbracket=\pi_{\d,\dd}(g'')=
  (\;\overline{g} \,R_\val^{\;\d},\,\widecheck{\underline{g}}\,
  R_\val^{\;\dd}\,) \in\La_{\d,\dd}\;.
\end{equation}

\blemm\label{lem:proprivarphi} The map $\varphi_{\d,\dd}$ is proper,
surjective, and satisfies the following properties.
\begin{enumerate}
\item\label{item1:proprivarphi} We have $(\varphi_{\d,\dd})_*\mu_{\La_{\d,\dd}}
  =(q-1)\, \mu_{\Sh_\d}\otimes\mu_{\Sh_\dd}$.
\item\label{item2:proprivarphi} For every $g\in G^{\sharp}$ (as
  defined just above Proposition \ref{prop:HKprop4_3}), we have
  \[\varphi_{\d,\dd}(\llbracket \Lambda_g \rrbracket) =
  (\sh(\Lambda_g),\sh((\Lambda_g)^\perp))\;.
  \]
\item\label{item3:proprivarphi} For all functions $f_1:\Sh_\d\ra\RR$
  and $f_2:\Sh_\dd\ra\RR$ with finite support, denoting by $f_{1}
  \times f_{2}:\Sh_\d\times\Sh_\dd \ra\RR$ their product map $(x,y)
  \mapsto f_1(x)f_2(y)$, the composition function $(f_{1}\times f_{2})
  \circ\varphi_{\d,\dd}: \La_{\d,\dd}\ra\RR$ is compactly supported
  and $\rho_0$-locally constant with $\|(f_{1}\times f_{2}) \circ
  \varphi_{\d,\dd}\|_\infty\leq \,\|f_1\|_\infty
  \,\|f_2\|_\infty$ for the uniform norms $\|\;\|_\infty$ .
\end{enumerate}
\elemm

\dem By the identification of $G''/G''(R_\val)$ with $\La_{\d,\dd}$
given after the proof of Lemma \ref{lem:latmnhomogen}, the map
$\varphi_{\d,\dd}$ is the continuous map
\[
\begin{array}{rcl}
  G''/G''(R_\val)\!\!\!&\ra&\!\!\!(\GL_\d(\OOO_\val) \times\GL_\dd(\OOO_\val))
  \bs(\GL^1_\d(K_\val) \times\GL^1_\dd(K_\val))/
  (\GL_\d(R_\val) \times\GL_\dd(R_\val))\\
  \big(\begin{smallmatrix} \overline{g}&
    0 \\ 0 & \underline{g}\end{smallmatrix}\big)G''(R_\val)
  \!\!\!&\mapsto&\!\!\!
  \big(\GL_\d(\OOO_\val)\,\overline{g}\,\GL_\d(R_\val),
  \GL_\dd(\OOO_\val)\,\check{\underline{g}}\,\GL_\dd(R_\val)\big)\,.
\end{array}
\]
It is onto since, as seen in the proof of Lemma \ref{lem:masstotLamn},
every element of $\GL^1_\d(K_\val) \times \GL^1_\dd(K_\val)$ may be
written as $s'_\d(\lambda)\, \iota(g'')$ with $\lambda\in
\OOO_\val^{\;\times}$ (so that $s'_\d(\lambda) \in(\GL_\d(\OOO_\val)
\times\GL_\dd(\OOO_\val))$) and $g''\in G''$. Since the group
$\GL_\d(\OOO_\val)\times\GL_\dd(\OOO_\val)$ is compact, the map
$\varphi_{\d,\dd}$ is proper.

\medskip\noindent \eqref{item1:proprivarphi} The map
$\varphi_{\d,\dd}$ is the composition of the inclusion map
$\La_{\d,\dd}\hookrightarrow(\La_{\d}\times \La_{\dd})$ with the
canonical projection $\sh\times\sh:(\La_{\d}\times \La_{\dd})\ra
(\Sh_\d\times\Sh_\dd)$, which is measure preserving by Equation
\eqref{eq:relatmuShmuLat}. By Equation \eqref{eq:GLGLGL}, for every
continuous function $f$ with compact support on $\La_{\d}\times
\La_{\dd}$ which is constant on the $(\GL_\d(\OOO_\val) \times
\GL_\dd(\OOO_\val))$-orbits, the integral of the restriction $f_{\mid
  \La_{\d,\dd}} $ with respect to $\mu_{\La_{\d,\dd}}$ is a constant
time the integral of $f$ with respect to $\mu_{\La_{\d}} \otimes
\mu_{\La_{\dd}}$. Hence the measures $(\varphi_{\d,\dd})_*
\mu_{\La_{\d,\dd}}$ and $\mu_{\Sh_\d}\otimes\mu_{\Sh_\dd}$ are
proportional. Since the pushforward maps of measures preserve the
total mass, and by Equations \eqref{eq:relatmuShmuLat} and
\eqref{eq:muLatnLatmLamn}, we have
\[
\|\mu_{\Sh_\d}\otimes\mu_{\Sh_\dd}\|=\|\mu_{\La_{\d}}\otimes\mu_{\La_{\dd}}\|
=\frac{1}{q-1}\;\|\mu_{\La_{\d,\dd}}\|=\frac{1}{q-1}\;
\|(\varphi_{\d,\dd})_*\mu_{\La_{\d,\dd}}\|\,.
\]

\medskip\noindent \eqref{item2:proprivarphi} Let $g\in G^\sharp$. Let
$u^-\in U^-(\OOO_\val)$, $g''=\big(\begin{smallmatrix} \overline{g} &0
  \\ 0 & \underline{g} \end{smallmatrix}\big)\in G''$, $z\in Z$ and
$u^+\in U^+$ be such that $g=u^-\,g''\,z\,u^+$, as in Proposition
\ref{prop:refLUbloc}. Then by Equation \eqref{eq:assoccorrpairlat}, by
the definition of $\varphi_{\d,\dd}$ and by Proposition
\ref{prop:HKprop4_3} \hyperlink{HKprop4_3iii}{(iii)} and
\hyperlink{HKprop4_3iiiperp}{(iii)$^\perp$}, we have
\[
\varphi_{\d,\dd}(\llbracket \Lambda_g\rrbracket) =
\varphi_{\d,\dd}(\,\overline{g} \,R_\val^{\;\d},\,\widecheck{\underline{g}}\,
R_\val^{\;\dd}\,) =(\sh(\,\overline{g} \,R_\val^{\;\d}),\,
\sh(\,\widecheck{\underline{g}}\,R_\val^{\;\dd}\,))=
(\sh(\Lambda_g),\sh((\Lambda_g)^{\perp}))\;.
\]

\medskip\noindent \eqref{item3:proprivarphi} Since the map
$\varphi_{\d,\dd}$ is proper and since the function
$f_{1,2}=f_{1}\times f_{2}$ is compactly supported, the function
$f_{1,2}\circ \varphi_{\d,\dd}$ is compactly supported. Let us prove
that for all points $x,x_0\in \La_{\d,\dd}$ at distance at most
$\rho_0$, we have $f_{1,2} \circ\varphi_{\d,\dd}(x) =
f_{1,2}\circ\varphi_{\d,\dd} (x_0)$. Since
$\|f_{1,2}\circ\varphi_{\d,\dd}\|_\infty= \|f_{1,2}\|_\infty \leq
\|f_1\|_\infty\,\|f_2\|_\infty$, this will prove Assertion
\eqref{item3:proprivarphi}.

Recall that by the end of Subsection \ref{subsect:latticepair}, the
distance on $\La_{\d,\dd}$ is the quotient distance of the distance
$d_{G''}$ on $G''$ by the surjective map $\pi_{\d,\dd}$. Let $\wt{x}_0
\in G''$ be such that $\pi_{\d,\dd}(\wt{x}_0)=x_0$. Since the action
by right translations of $G''$ on itself is isometric and since the
preimages of $\pi_{\d,\dd}$ are the right orbits of $G''(R_\val)$ in
$G''$, there exists $\wt{x}\in G''$ such that $\pi_{\d,\dd}(\wt{x})=x$
and $d_{G''}(\wt{x},\wt{x}_0)\leq \rho_0$. Again since the action by
right translations of $G''$ on itself is isometric and by the
definition of $\rho_0$, we have $g=\wt{x}_0\,\wt{x}^{-1} \in
B_{G''}(\id,\rho_0)\subset G''(\OOO_\val)$. Since $\pi_{\d,\dd}$ is
$G''$-equivariant, we have $g\,x=\pi_{\d,\dd}(g\,\wt{x})=
\pi_{\d,\dd}(\wt{x_0})=x_0$. Since $\varphi_{\d,\dd}$ is constant on
the left orbits of $G''(\OOO_\val)$, we have $\varphi_{\d,\dd}(x)
=\varphi_{\d,\dd}(x_0)$, therefore $f_{1,2} \circ\varphi_{\d,\dd}(x) =
f_{1,2}\circ\varphi_{\d,\dd}(x_0)$.  \cqfd

\section{Joint equidistribution of primitive partial lattices}
\label{sect:Equidistribution}

\subsection{The correspondence between primitive partial lattices
  and integral group elements}
\label{subsect:correspondance}

Let $\Ga=\SL_\D(R_\val)$ be the modular group of integral points of
$G=\SL_\D(K_\val)$. The aim of this subsection is to naturally and
injectively associate elements in $\Ga$ to primitive $\d$-lattices in
$R_\val^{\;\D}$.
%%
%\todo{\tiny check if need introduce left cosets of $\Ga$ for
% primitive lattices with wrong congruence of covolume, as GN still works}
%%
We start by introducing the subsets of the group $G$ and of the moduli
spaces $\PL_{\d,\D}$, $\Gr_{\d,\D}$, $\La_{\d,\dd}$ which will be
technically useful.

We fix from now on a compact-open strict fundamental domain ${\mathcal
  D}$ for the action by translations of $R_\val$ on $K_\val$ (for
instance $\D=\pi_\val\OOO_\val$ when $K=\FF_q(Y)$ and $\val=
\val_\infty$ as defined at the end of Subsection
\ref{subsect:functionfields}), such that for every $x\in {\mathcal
  D}$, the (closed) ball $B(x,q_\val^{\;-1})=x+ \pi_\val\OOO_\val$ in
$K_\val$ is contained in ${\mathcal D}$. This is possible since
$R_\val\cap\pi_\val\OOO_\val =\{0\}$ by Equation \eqref{eq:inversiRv}.
We thus have a compact-open strict fundamental domain
\[
\Box=\big\{(x_{i,j})_{1\leq i\leq \d,1\leq j\leq \dd} \in\M_{\d,\dd} (K_\val):
\forall\;i\in\llbracket 1,\d\rrbracket,\;
\forall\;j\in\llbracket 1,\dd\rrbracket,\;\; x_{i,j}\in{\mathcal D}\big\}
\]
for the action by translations of the $R_\val$-lattice $\M_{\d,\dd}
(R_\val)$ on the $K_\val$-vector space $\M_{\d,\dd}(K_\val)$, for
instance $\Box=\M_{\d,\dd}(\pi_\val\OOO_\val)$ if $K=\FF_q(Y)$ and
$\val= \val_\infty$. We also fix a closed-open strict fundamental
domain $G''_\lozenge$ for the action by translations on the right of
the discrete subgroup $G''(R_\val)$ on $G''$, so that we have $G''=
G''_\lozenge\; G''(R_\val)$ with unique writing.

For every $r\in\ZZ$ and for all measurable subsets $\Psi$ of
$\M_{\dd,\d}(K_\val)$ and $\F$ of $G''$,  we define
\begin{equation}\label{eq:defiUpPhiGseccurlE}
U^-_{\Psi}=\Big\{\big(\begin{smallmatrix}I_{\d} &\; 0\\ \beta &
I_{\dd}\end{smallmatrix}\big)\in U^-:\beta\in \Psi\Big\},\quad
G''_{\F}=G''_\lozenge\cap \F\;,
\end{equation}
\[
Z_r=\Big\{
\big(\begin{smallmatrix}\lambda I_{\d} &\;  0\\ 0  & \lambda' I_{\dd}
\end{smallmatrix}\big)\in Z:\val(\lambda)=\frac{\lcm\{\d,\dd\}}{\d}
\;r\Big\}\quad \text{and}\quad
U^+_{\Box}=\Big\{\big(\begin{smallmatrix}I_{\d} &\;  \ga\\ 0 & I_{\dd}
\end{smallmatrix}\big)\in U^+:\ga\in \Box\Big\}
\;.
\]
Note that $Z=\bigsqcup_{r\in\ZZ} Z_r$ by Equation \eqref{eq:valeurr}.
As defined just before Proposition \ref{prop:HKprop4_3}, let
\[
G^\sharp=\{\big(\begin{smallmatrix} \alpha &\ga \\ \beta &
  \delta \end{smallmatrix}\big)\in G:
\val(\det\alpha)\in\lcm\{\d,\dd\}\;\ZZ,\;\;\beta\alpha^{-1}\in
\M_{\dd,\d}(\OOO_\val)\}\;,
\]
so that the product map $(u^-,g'',z,u^+)\mapsto u^-\,g''\,z\,u^+$ from
$U^-(\OOO_\val)\times G''\times Z\times U^+$ to $G^\sharp$ is a
homeomorphism, by Proposition \ref{prop:refLUbloc}. For every closed
subgroup $H$ of $G$, let
\begin{equation}\label{eq:defiHsharp}
H^\sharp=H\cap G^\sharp\;.
\end{equation}
We also define the corresponding subset of the set $\PL_{\d,\D}$ of
primitive $\d$-lattices in $K_\val^{\;\D}$ by
\[
\PL_{\d,\D}^{\sharp}=\Ga^\sharp(R_\val^{\;\d}\times\{0\})\;.
\]

For every measurable subset $\E$ of $\La_{\d,\dd}$, using the map
$\pi_{\d,\dd}: G''\ra \La_{\d,\dd}$ introduced in Equation
\eqref{eq:defiprok}, we define
\begin{equation}\label{eq:defiwt}
\wt{\E}= \pi_{\d,\dd}^{\;-1}(\E)\subset G''\;,
\end{equation}
which is a measurable subset of $G''$ invariant by the translations on
the right by $G''(R_\val)$.

For every measurable subset $\Phi$ of $\Gr_{\d,\D}^{\flat}$, using the map 
$\operatorname{orb}_\d:\M_{\dd,\d} (\OOO_\val)\ra\Gr_{\d,\D}^{\flat}$
defined in Subsection \ref{subsect:grassmann}, which is
isometric by Lemma \ref{lem:isometrygrass}, we define
\begin{equation}\label{eq:defiwtorb}
  \wt \Phi= \operatorname{orb}_\d^{\,-1}
  (\Phi)\subset \M_{\dd,\d}(\OOO_\val)\;,
\end{equation}
which is a measurable subset of $\M_{\dd,\d}
(\OOO_\val)$.

For every $r\in\ZZ$ and for all measurable subsets $\Phi$ of
$\Gr_{\d,\D}^{\flat}$ and $\E$ of $\La_{\d,\dd}$, using the notation
of Equation \eqref{eq:defiUpPhiGseccurlE} with $\Psi=\wt\Phi$ and
$\F=\wt\E$, we finally define
\begin{equation}\label{eq:defiOmega}
  \Omega=U^-(\OOO_\val)\;G''_{\lozenge}\; Z\;U^+_{\Box}\subset G^\sharp
  \quad\text{and}\quad
  \Omega_{\Phi,\E,r}=U^-_{\wt\Phi}\;G''_{\wt\E}\;
  Z_r\;U^+_{\Box}\subset \Omega\;.
\end{equation}

\blemm\label{lem:calcmesOmegasub} We have
$\mu_G(\Omega_{\Phi,\E,r})=q^{(\ggg-1)\,\d\,\dd}\;
q_\val^{\;\D\lcm(\d,\dd)\,r}\;\mu_{\Gr_{\d,\D}}(\Phi)\;\mu_{\La_{\d,\dd}}(\E)$.
\elemm

\dem By Proposition \ref{prop:refLUbloc}, we have
\begin{equation}\label{eq:volGPhiEFr}
\mu_G(\Omega_{\Phi,\E,r})=c_1\;\mu_{U^-}(U^-_{\wt\Phi})\;
\mu_{G''}(G''_{\wt\E})\;\mu_{U^+}(U^+_{\Box})\;
\int_{z\in Z_r}|\chi_\d(z)|^{\D\,\d}\;d\mu_{Z}(z)\;.
\end{equation}
By Equations \eqref{eq:normalisehaar}, \eqref{eq:defiwtorb} and
\eqref{eq:relatHaarmnmuGr}, and since $\Phi\subset
\Gr_{\d,\D}^\flat$, we have
\begin{equation}\label{eq:volGUmoins}
\mu_{U^-}(U^-_{\wt\Phi})=\operatorname{Haar}_{\dd,\d}(\wt\Phi)=
\operatorname{Haar}_{\dd,\d}(\operatorname{orb}_\d^{\;-1}(\Phi))=
c_1^{\;-1}\;\mu_{\Gr_{\d,\D}}(\Phi)\;.
\end{equation}
By Equations \eqref{eq:normalisehaar} and \eqref{eq:covolRv}, we have
\begin{equation}\label{eq:volGUplus}
\mu_{U^+}(U^+_{\Box})=\operatorname{Haar}_{\d,\dd}(\Box)=
\covol(R_\val^{\;\d\,\dd})=q^{(\ggg-1)\,\d\,\dd}\;.
\end{equation}
Note that $z=\big(\begin{smallmatrix} \lambda I_\d & 0 \\ 0 & \mu
  I_{\dd}\end{smallmatrix} \big)$ belongs to $Z_r$ if and only if
$\chi_\d(z)=\lambda=\pi_\val^{\;-\frac{\lcm\{\d,\dd\}}{\d}\;r}\;$. Since
the Haar measure $\mu_Z$ is normalized so that
$\mu_Z(Z(\OOO_\val))=1$ (see Equation \eqref{eq:muZ}), we have
\begin{equation}\label{eq:volGZ}
\int_{z\in Z_r}|\chi_\d(z)|^{\D\,\d}d\mu_{Z}(z)=
\Big|\pi_\val^{\;-\frac{\lcm\{\d,\dd\}}{\d}\;r}\Big|^{\D\,\d}=
q_\val^{\;\D\lcm(\d,\dd)\,r}\;.
\end{equation}

By Equation \eqref{eq:defiUpPhiGseccurlE} and by the measure
preserving identification of $G''/G''(R_\val)$ and $\La_{\d,\dd}$
done in Subsection \ref{subsect:latticepair} (just before Lemma
\ref{lem:masstotLamn}), we have
\begin{equation}\label{eq:volGEF}
\mu_{G''}(G''_{\wt\E})=\mu_{G''}(G''_{\lozenge}\cap \wt\E)=
\mu_{G''/G''(R_\val)}(\E)=
\mu_{\La_{\d,\dd}}(\E)\;.
\end{equation}
Lemma \ref{lem:calcmesOmegasub} follows from Equation
\eqref{eq:volGPhiEFr} by plugging in it the computations of Equations
\eqref{eq:volGUmoins}, \eqref{eq:volGEF}, \eqref{eq:volGUplus} and
\eqref{eq:volGZ}.  \cqfd

\bigskip
The following result gives a precise 1-to-1 correspondence between
partial $R_\val$-lattices in $\PL_{\d,\D}^{\sharp}=\Ga^\sharp(R_\val^{\;\d}
\times\{0\})$ and appropriate matrices in the discrete group $\Ga=
\SL_\D(R_\val)$. 

\bprop\label{prop:matrepprimlat} The map $g\mapsto \Lambda_g=
g(R_\val^{\;\d}\times\{0\})$ from $\Ga\cap\Omega$ to
$\PL_{\d,\D}^{\sharp}$ is a bijection such that for every nonzero
ideal $I$ of $R_\val$, for every $r\in\ZZ$ and for all measurable
subsets $\Phi$ of $\Gr_{\d,\D}^{\flat}$ and $\E$ of $\La_{\d,\dd}$,
the following two assertions are equivalent
\smallskip
\begin{enumerate}
\item\label{item1:matrepprimlat} the integral matrix $g\in \Ga\cap
  \Omega$ lies in $\Omega_{\Phi,\E,r}\cap\Ga_I$,
\item\label{item2:matrepprimlat} the primitive $\d$-lattice
  $\Lambda_g\in \PL_{\d,\D}^{\sharp}$ satisfies that $\Lambda_g\in
  \PL_{\d,\D}(I)$ (as defined in the beginning of Subsection
  \ref{subsect:congruence}), $V_{\Lambda_g} \in \Phi$ (as defined in
  Subsection \ref{subsect:partiallatt}), $\llbracket
  \Lambda_g\rrbracket \in\E$ (as defined in Equation
  \eqref{eq:assoccorrpairlat}) and $\overline{\covol}(\Lambda_g)=
  \overline{\covol}((\Lambda_g)^{\perp})= q_\val^{\;\lcm\{\d,\dd\}r}$.
\end{enumerate}
\eprop

\dem Since $\Omega\subset G^\sharp$, if $g\in\Ga\cap \Omega$, then
$g\in\Ga^\sharp$, thus $\Lambda_g= g(R_\val^{\;\d}\times\{0\})\in
\PL_{\d,\D}^{\sharp}$. Hence the above map is well defined.

For all $g=\big(\begin{smallmatrix} \alpha &\ga \\ \beta & \delta
\end{smallmatrix}\big)\in G^\sharp$ and $p^+=\big(\begin{smallmatrix}
  \alpha' &\ga' \\ 0 & \delta' \end{smallmatrix}\big)\in P^+=G''ZU^+
\subset G^\sharp$, since the first column of $g\,p^+$ is
$\big(\begin{smallmatrix} \alpha\alpha'\\ \beta\alpha'
\end{smallmatrix}\big)$, since $(\beta\alpha')(\alpha\alpha')^{-1}
=\beta\alpha^{-1}$ and since $\val(\det(\alpha\alpha'))=
\val(\det\alpha)+\val(\det\alpha')$, we have $g\,p^+\in G^\sharp$. In
particular, the action by right translations on $\Ga=\SL_\D(R_\val)$
of its subgroup $P^+(R_\val)$ defined in Equation
\eqref{eq:defiPplus}, which satisfies $P^+(R_\val)=G''(R_\val)\,
U^+(R_\val)$, preserves $\Ga^\sharp=\Ga\cap G^\sharp$.  By the
identification given just after the statement of Lemma
\ref{lem:descriphomogPL}, we thus have
\begin{equation}\label{eq:descripPLmdsharp}
\PL_{\d,\D}^{\sharp}=\Ga^\sharp/P^+(R_\val)\;.
\end{equation}

Let us prove that
\begin{equation}\label{eq:decompGsharp}
  G^\sharp=\coprod_{\ga\in P^+(R_\val)}\Omega\,\ga\;.
\end{equation}
Since $P^+(R_\val)\subset \Ga$, this will imply that $\Ga^\sharp=\Ga\cap
G^\sharp= \coprod_{\ga\in P^+(R_\val)}(\Ga\cap\Omega)\ga$. By Equation
\eqref{eq:descripPLmdsharp}, this will imply that the map $g\mapsto
\Lambda_g= g(R_\val^{\;\d}\times\{0\})$ from $\Ga\cap\Omega$ to
$\PL_{\d,\D}^{\sharp}$ is a bijection.

In order to prove Equation \eqref{eq:decompGsharp}, let us fix $g\in
G^\sharp$. By Proposition \ref{prop:refLUbloc}, there exist unique
elements $u^-\in U^-(\OOO_\val)$, $g''=\big(\begin{smallmatrix}
  \overline{g} &0 \\ 0 & \underline{g} \end{smallmatrix}\big)\in G''$,
$z\in Z$ and $u^+\in U^+$ such that $g=u^-\,g''\,z\,u^+$. By the
definition of the strict fundamental domain $G''_\lozenge$, there
exist a unique $f''=\big(\begin{smallmatrix} \overline{f} &0 \\ 0 &
  \underline{f} \end{smallmatrix}\big)\in G''_{\lozenge}$ and a unique
$\ga''=\big(\begin{smallmatrix} \overline{\ga} &0 \\ 0 &
  \underline{\ga} \end{smallmatrix}\big)\in G''(R_\val)$ such that
\begin{equation}\label{eq:ovundpp}
  \overline{f}\;\overline{\ga}=\overline{g},\quad
  \underline{f}\;\underline{\ga}=\underline{g}
  \quad\text{and}\quad  g''=f''\,\ga''\;.
\end{equation}
Since $Z$ centralizes $G''$, we have $g=u^-\,f''\,z\,\ga''\,u^+$.
Since $G''$ normalizes $U^+$, and by the definition of $\Box$, there
exist unique elements $u^+_0\in U^+_{\Box}$ and $\ga_0\in
U^+(R_\val)$ such that $\ga''\,u^+(\ga'')^{-1}= u^+_0\,\ga_0$.
Defining $\ga=\ga_0\ga''\in P^+(R_\val)$, we have
\begin{equation}\label{eq:decompOmega}
g=u^-\,f''\,z\,(\ga''\,u^+(\ga'')^{-1})\,\ga''=u^-\,f''\,z\,u^+_0\ga
\end{equation}
and $u^-\,f''\,z\,u^+_0\in \Omega$ by the definition of $\Omega$ in
Equation \eqref{eq:defiOmega}.  Since the writing $h=u^+_0\ga$ of an
element $h\in P^+(R_\val)=U^+(R_\val)\, G''(R_\val)$ with $u^+_0\in
U^+(R_\val)$ and $\ga\in G''(R_\val)$ is unique, this proves Equation
\eqref{eq:decompGsharp}.

\medskip
Let us now assume that $g\in\Ga\cap \Omega$. By the uniqueness of the
writing in Equation \eqref{eq:decompOmega}, we may uniquely write
$g=u^-\,f''\,z\,u^+_0$ with $u^-\in U^-(\OOO_\val)$, $f''\in
G''_{\lozenge}$, $z\in Z$ and $u^+_0\in U^+_{\Box}$. By the definition
of $\Omega_{\Phi,\E,r}$ in Equation \eqref{eq:defiOmega}, we have
$g\in \Omega_{\Phi,\E,r}$ if and only if $u^-\in U^-_{\wt \Phi}$,
$f''\in G''_{\wt\E}$ and $z\in Z_r$.

We have $u^-\in U^-_{\wt \Phi}$ if and only if there exists $\beta\in
\wt \Phi$ with $u^-=\big(\begin{smallmatrix}I_{\d} &\; 0\\ \beta &
  I_{\dd}\end{smallmatrix}\big)$, hence if and only if there exists
$\beta\in \wt \Phi={\operatorname{orb}_\d}^{-1}(\Phi)$ with
$\operatorname{orb}_\d(\beta)= u^-(K_\val^{\;\d}\times\{0\})=
V_{\Lambda_{u^-}}$ by the definition of
$\operatorname{orb}_\d$ in Subsection \ref{subsect:grassmann},
therefore if and only if $V_{\Lambda_g}\in \Phi$ by Proposition
\ref{prop:HKprop4_3} \hyperlink{HKprop4_3i}{(i)}.

By Equations \eqref{eq:defiUpPhiGseccurlE} and \eqref{eq:defiwt}, and
by the definition of $\pi_{\d,\dd}$ in Equation \eqref{eq:defiprok},
we have
\[
G''_{\wt\E}=G''_\lozenge\cap \pi_{\d,\dd}^{\;-1}(\E)=
\big\{\big(\begin{smallmatrix} \overline{g} &0
  \\ 0 & \underline{g} \end{smallmatrix}\big)\in
G''_\lozenge:(\;\overline{g}
\,R_\val^{\;\d},\,\widecheck{\underline{g}}\, R_\val^{\;\dd}\,)
\in\E\big\}\;.
\]
Hence by Equation \eqref{eq:assoccorrpairlat} and since $f''\in
G''_{\lozenge}$, we have $f''\in G''_{\wt\E}$ if and only if
$\llbracket \Lambda_g\rrbracket\in \E$.

We have $z\in Z_r$ if and only if $z=\begin{pmatrix}
\pi_\val^{\;-\frac{\lcm\{\d,\dd\}}{\d} \,r}I_\d &0\\ 0 &
\pi_\val^{\;\frac{\lcm\{\d,\dd\}}{\dd} \,r}I_\dd \end{pmatrix}$, hence
if and only if $\overline{\covol}(\Lambda_g)= q_\val^{\;\lcm\{\d,\dd\}
  \,r}$ by Proposition \ref{prop:HKprop4_3}
\hyperlink{HKprop4_3ii}{(ii)} and by Equation
\eqref{eq:formulezrelatt}. Note that we have $\overline{\covol}
((\Lambda_g)^{\perp})= \overline{\covol}(\Lambda_g)$ by Proposition
\ref{prop:HKprop4_3} \hyperlink{HKprop4_3iiperp}{(ii)$^\perp$}.

The fact that $g\in \Ga_I$ if and only if $\Lambda_g\in
\PL_{\d,\D}(I)$ has been shown in Lemma \ref{lem:indexHeckesubgrou} (1).
This concludes the proof of Proposition \ref{prop:matrepprimlat}.
\cqfd

\subsection{Counting in well-rounded families}
\label{subsect:wellround}

A crucial tool of this paper is a counting result of lattice points by
Gorodnik and Nevo \cite{GorNev12}. In this subsection, after the
necessary definitions, we recall from \cite{HorPau22} an adaptation of
the Gorodnik-Nevo result, and we proceed to the construction of the
well-rounded family of subsets to which we will apply it.

Let ${\bf G'}$ be an absolutely connected and simply connected
semi-simple algebraic group over $K_\val$, which is almost
$K_\val$-simple.  Let $G'={\bf G'}(K_\val)$ be the locally compact
group of $K_\val$-points of ${\bf G'}$. Let $\Ga'$ be a
nonuniform\footnote{This implies that ${\bf G'}$ is isotropic over
$K_\val$, as part of the assumptions of \cite{GorNev12}.}  lattice in
$G'$, and let $\mu_{G'}$ be any (left) Haar measure of $G'$. Note that
$G'=G=\SL_\D(K_\val)$ and $\Ga'=\Ga_I$ (defined in Subsection
\ref{subsect:congruence}) satisfy these assumptions for every nonzero
ideal $I$ of $R_\val$.

Let $\rho>0$. Let $(\V'_\epsilon)_{\epsilon>0}$ be a fundamental
system of neighborhoods of the identity in $G'$, which 

$\bullet$~ is symmetric (that is, $x\in \V'_\epsilon$ if and only if
$x^{-1}\in \V'_\epsilon$),

$\bullet$~ is nondecreasing with $\epsilon$ (that is,
$\V'_\epsilon\subset\V'_{\epsilon'}$ if $\epsilon\leq \epsilon'$), and

$\bullet$~ has {\it upper local dimension} $\rho$, that is, there
exist $m_1,\epsilon_1>0$ such that $\mu_{G'}(\V'_\epsilon)\geq
m_1\,\epsilon^\rho$ for every $\epsilon\in\;]0,\epsilon_1[\,$.

Let $C\geq 0$. Let $(\Z_n)_{n\in\NN}$ be a family of measurable
subsets of $G'$.  We define
\[
(\Z_n)^{+\epsilon}= \V'_\epsilon \Z_n\V'_\epsilon=
\bigcup_{g,h\in\V'_\epsilon}g \Z_n h \;\;\;{\rm and}\;\;\;
(\Z_n)^{-\epsilon}= \bigcap_{g,h\in\V'_\epsilon}g \Z_n h\;.
\] 
The family $(\Z_n)_{n\in\NN}$ is said to be {\it $C$-Lipschitz
  well-rounded} with respect to $(\V'_\epsilon)_{\epsilon>0}$ if there
exist $\epsilon_0 >0$ and $n_0\in\NN$ such that for all
$\epsilon\in\; ]0,\epsilon_0[$ and $n\geq n_0$, we have
\[
\mu_{G'}((\Z_n)^{+\epsilon})\leq
(1+C\,\epsilon)\;\mu_{G'}((\Z_n)^{-\epsilon})\;.
\]

We refer to \cite[Theo.~4.1]{HorPau22} for a proof of the following
adaptation of results of Gorodnik-Nevo \cite{GorNev12}.

\btheo\label{theo:GorodnikNevo} For every $\rho>0$, there exists
$\tau(\Ga')\in \;]0, \frac{1}{2(1+\rho)}]$ such that for every
symmetric nondecreasing fundamental system
$(\V'_\epsilon)_{\epsilon>0}$ of neighborhoods of the identity in $G'$
with upper local dimension $\rho$, for every $C\geq 0$, for every
family $(\Z_n)_{n\in\NN}$ of measurable subsets of $G'$ that is
$C$-Lipschitz well-rounded with respect to
$(\V'_\epsilon)_{\epsilon>0}$, and for every $\delta>0$, 
as $n\ra +\infty$, we have 
\[
\Big|\;\card(\Z_{n}\cap\Ga')-
\frac{1}{\|\mu_{G'/\Ga'}\|}\,\mu_{G'}(\Z_{n})\;\Big|=
\bigO\big(\mu_{G'}(\Z_{n})^{1-\tau(\Ga')+\delta}\big)\;,
\]
where the function $\bigO(\cdot)$ depends only on $G',\Ga',\delta,C,
(\V'_\epsilon)_{\epsilon>0},\rho$. \cqfd
\etheo

We will use, as a fundamental system of neighborhoods of the identity
element in $G$, a family of compact-open subgroups of $G(\OOO_\val)$
given by the kernels of the morphisms of reduction modulo $\pi_\val^{\;N}
\OOO_\val$ for appropriate $N\in\NN$. For every $\epsilon>0$, let
$N_\epsilon= \big\lfloor-\log_{q_\val} \epsilon\, \big\rfloor$ so that
$N_\epsilon\geq 1$ if and only if $\epsilon\leq \frac{1}{q_\val}$. Let
$\V_\epsilon=G(\OOO_\val)$ if $\epsilon>\frac{1}{q_\val}$ and
otherwise let
\begin{align}
\V_\epsilon&= \ker(G(\OOO_\val)\ra
\SL_{\D}(\OOO_\val/\pi_\val^{\;N_\epsilon}\OOO_\val)) \nonumber\\ &=
\big\{I_\D+\pi_\val^{\;N_\epsilon}X: X\in \M_\D(\OOO_\val)\big\}\cap G
\;.\label{eq:defiVeps}
\end{align}
The family $(\V_\epsilon)_{\epsilon>0}$ is indeed nondecreasing and we
have $\bigcap_{\epsilon>0} \V_\epsilon=\{\id\}$. Note that for all
$\epsilon_1,\dots, \epsilon_k>0$, we have
\begin{align*}
\min\{N_{\epsilon_1},\cdots, N_{\epsilon_k}\}&\geq
\min\{-\log_{q_\val}\epsilon_1,\cdots, -\log_{q_\val}\epsilon_k\}-1
\\ &\geq -\log_{q_\val}(\epsilon_1+\cdots+\epsilon_k)-1
\geq N_{q_\val(\epsilon_1+\cdots+\epsilon_k)}\,,
\end{align*}
hence
\begin{equation}\label{eq:minNepsilons}
  \V_{\epsilon_1}\V_{\epsilon_2}\cdots\V_{\epsilon_k}\;\;\subset\;\;
  \V_{q_\val(\epsilon_1+\cdots+\epsilon_k)}\,.
\end{equation}
For every subgroup $H$ of $G$, let $\V_\epsilon^H = \V_\epsilon \cap
H$. The index of $\V_\epsilon$ in $G(\OOO_\val)$ is given by Lemma
\ref{lem:calcindexN} with $N=N_\epsilon$.

\medskip
We denote the operator norm of a linear operator $\ell$ of the normed
$K_\val$-algebra $\M_\D(K_\val)$ (for the supremum norm defined before
Lemma \ref{lem:isometrygrass}) by
\[
\|\ell\|=\max\Big\{\frac{\|\ell(X)\|}{\|X\|} :
X\in\M_\D(K_\val)\ssm\{0\}\Big\} \in q_\val^{\;\ZZ} \cup\{0\}\,,
\]
so that $\ell(\M_\D(\OOO_\val))\subset \M_\D(\pi_\val^{\;
  -\log_{q_\val} \|\ell\|} \OOO_\val)$ if $\ell$ is invertible.  For
every $g\in G$, recall that $\operatorname{Ad} g$ is the linear
automorphism $x\mapsto gxg^{-1}$ of $\M_\D(K_\val)$. Also recall that
$P^-=U^-G''Z$ is the lower triangular-by-blocs subgroup of $G$.

\blemm\label{eq:calVdivgroup}
For all $\epsilon\in\;]0,\frac{1}{q_\val}]$ and $g\in G$, we have
\[
g\,\V_\epsilon\, g^{-1}\;\subset\;
\V_{\|\operatorname{Ad} g\,\|\;\epsilon}\;,\;\;\;
\V_\epsilon=\V^{P^-}_\epsilon\;\V^{U^+}_\epsilon\;\;\;{\rm and}\;\;\;
\V^{P^-}_\epsilon=\V^{U^-}_\epsilon\;\V^{G''}_\epsilon\;.
\]
Furthermore, the number $\rho=\D^2-1$ is an upper local dimension
of the family $(\V_\epsilon)_{\epsilon>0}$.
\elemm

\dem The first claim follows from the fact that
\begin{align*}
g\,\V_\epsilon \,g^{-1}&=
\big(I_\D+\pi_\val^{\;N_\epsilon}g\M_\D(\OOO_\val)g^{-1}\big)\cap G\\&
\subset \big(I_\D+\pi_\val^{\;N_\epsilon-\log_{q_\val}\|\operatorname{Ad} g\,\|}
\M_\D(\OOO_\val)\big)\cap G=
\V_{\epsilon\,\|\operatorname{Ad} g\,\|}\;.
\end{align*}

Note that $\V_\epsilon$ is contained in $\U_G$ (defined in Equation
\eqref{eq:defiUsugG}). Indeed, if $g=\big(\begin{smallmatrix} \alpha &\ga
\\ \beta & \delta\end{smallmatrix}\big)\in \V_\epsilon$ then
$\alpha\in I_\d+\pi_{\val}^{N_\epsilon}\M_\d(\OOO_\val)$, hence
$\det(\alpha)\in 1+\pi_{\val}^{N_\epsilon}\OOO_\val$, so that
$\val(\det\alpha) =0\in\lcm\{\d,\dd\}\ZZ$ since $N_\epsilon\geq 1$ and
therefore $g\in\U_G$.

By Proposition \ref{prop:refLUbloc}, we may hence uniquely write any
$g=\big(\begin{smallmatrix} \alpha &\ga \\ \beta & \delta
\end{smallmatrix}\big)\in \V_\epsilon$ as $g=u^-\,g''\,z\,u^+$
with $u^-\in U^-$, $g''\in G''$, $z\in Z$ and $u^+\in U^+$. Since
$\alpha\in\M_\d(\OOO_\val)$ satisfies $\val(\det\alpha) =0$ and by
Equation \eqref{eq:inverseLUbloc}, we have $\alpha\in
\GL_\d(\OOO_\val)$, $u^\pm\in U^\pm(\OOO_\val)$, $g''\in
G''(\OOO_\val)$ and $z=I_\D$. Furthermore, since $g\in \V_\epsilon$,
we have $\alpha=I_\d\mod \pi_\val^{\;N_\epsilon}$, $\ga=0\mod
\pi_\val^{\;N_\epsilon}$, $\beta=0\mod \pi_\val^{\;N_\epsilon}$ and
$\delta=I_\dd\mod \pi_\val^{\;N_\epsilon}$. Therefore again by
Equation \eqref{eq:inverseLUbloc}, we have $u^\pm\in
\V_\epsilon^{U^\pm}$ and $g''\in \V_\epsilon^{G''}$. This proves the
second and third claims.

In order to prove the last claim, let us apply Lemma
\ref{lem:calcindexN} with $N=N_\epsilon=\big\lfloor-\log_{q_\val}
\epsilon\, \big\rfloor$, so that $N_\epsilon-1\leq - \log_{q_\val}
\epsilon$. Since $\mu_G(G(\OOO_\val))=1$, we hence have
\begin{align*}
  \mu_G(\V_\epsilon)&
  =\frac{\mu_G(G(\OOO_\val))}{[G(\OOO_\val):\V_\epsilon]}
  =\frac{q_{\val}-1}{q_{\val}^{\,(N_\epsilon-1)(\D^{2}-1)}}
  \prod_{i=0}^{\D-1}(q_{\val}^{\;\D}-q_{\val}^{\;i})^{-1}\geq
  \epsilon^{\D^{2}-1}(q_{\val}-1)
\prod_{i=0}^{\D-1}(q_{\val}^{\;\D}-q_{\val}^{\;i})^{-1}\;.
\end{align*}
This proves the result. \cqfd

\medskip
We will need the following effective version of the refined LU
decomposition by blocks given in Proposition \ref{prop:refLUbloc}. We
denote by $c:h\mapsto c_h$ the continuous function from $G$ to
$[0,+\infty[$ defined by $c_h= \max\{\|\operatorname{Ad} h\,\|,
    \|\operatorname{Ad} h^{-1}\,\|\}$ for every $h\in G$.

\blemm\label{eq:effectiveLU} For all $\epsilon\in\; ]0,\frac{1}
  {q_\val}]$, $u^-\in U^-$, $g''\in G''$, $z\in Z$ and $u^+\in U^+$,
if $|\chi_\d(z)|\geq 1$ and $g=u^-\,g''\,z\,u^+$, then
\[
\V_\epsilon \;g\;\V_{\epsilon}
\;\subset\; u^-\;\V_{c_{g''}q_\val(c_{u^-g''}+c_{u^+})\,\epsilon}^{U^-}\;g''\;
\V_{q_\val(c_{u^-g''}+c_{u^+})\,\epsilon}^{G''}\;z\;
\V_{q_\val(c_{u^-g''}+2c_{u^+})\,\epsilon}^{U^+} \,u^+\;.
\]
\elemm

\dem In order to simplify the notation, let $p=u^-g''$ and $u=u^+$, so
that $g=pzu$. If $z=\big(\begin{smallmatrix}\lambda I_{\d} & 0\\ 0 &
  \mu I_{\dd}\end{smallmatrix}\big)$, we have
\[
z\big(\begin{smallmatrix}I_{\d} & 0\\
\beta & I_{\dd}\end{smallmatrix}\big) z^{-1}
=\big(\begin{smallmatrix} I_{\d} & 0\\
\mu\,\lambda^{-1}\beta\; & I_{\dd} \end{smallmatrix}\big)
\quad\text{and}\quad z^{-1} \big(\begin{smallmatrix}I_{\d} & \ga\\
0 & I_{\dd}\end{smallmatrix}\big)z  =
\big(\begin{smallmatrix}I_{\d} &\; \lambda^{-1}\mu\,\ga\\ 0 & I_{\dd}
\end{smallmatrix}\big).
\]
Since $|\lambda|=|\chi_\d(z)|\geq 1$ and $\det z=1$ so that
$|\mu|=|\lambda|^{-\frac{\d}{\dd}}\leq 1$, for every $\epsilon'>0$, we
have
\begin{equation}\label{eq:shrinking}
  z\;\V_{\epsilon'}^{U^-}z^{-1}\subset\V_{\epsilon'}^{U^-}\quad\text{and}\quad
  z^{-1}\,\V_{\epsilon'}^{U^+}z\subset\V_{\epsilon'}^{U^+}\;.
\end{equation}

Using for the following sequence of equalities and inclusions respectively

$\bullet$~ the first claim of Lemma \ref{eq:calVdivgroup} for the
first inclusion,

$\bullet$~ the second and third claims of Lemma \ref{eq:calVdivgroup}
for the second equality,

$\bullet$~ the claim on the left in Equation \eqref{eq:shrinking} and
the fact that $Z$ centralises $G''$ for the second inclusion,

$\bullet$~ the fact that $\V_{c_u\epsilon}$ is a normal subgroup of
$G(\OOO_\val)$ that contains $\V_{c_u\epsilon}^{U^-}\,
\V_{c_u\epsilon}^{G''}$ for the third inclusion,

$\bullet$~ the second claim of Lemma \ref{eq:calVdivgroup} for the
fourth equality,

$\bullet$~ twice the claim on the right in Equation \eqref{eq:shrinking}
for the fourth inclusion,

$\bullet$~ twice Equation \eqref{eq:minNepsilons} with $k=2$ and
$k=3$, defining the constants $c''_1=q_\val(c_p+c_u)$ and
$c''_2=q_\val(c_p+2c_u)$ for the fifth inclusion,

$\bullet$~ again the third claim of Lemma \ref{eq:calVdivgroup} for
the sixth equality,

$\bullet$~ the fact that $G''$ normalizes $U^-$ and again the first
claim of Lemma \ref{eq:calVdivgroup} for the last inclusion,

we have
\begin{align*}
\V_\epsilon \,g\,\V_\epsilon&= p\;p^{-1}\V_\epsilon
\,p\;z\;u\,\V_\epsilon \,u^{-1}\;u \;\subset\; p\,\V_{c_p\epsilon}
\,z\,\V_{c_u\epsilon} \,u
%\\ &
=p\,\V_{c_p\epsilon}^{P^-}\,\V_{c_p\epsilon}^{U^+} \,z\;
\V_{c_u\epsilon}^{U^-}\,\V_{c_u\epsilon}^{G''}\,\V_{c_u\epsilon}^{U^+}
\,u \\ &=p\,\V_{c_p\epsilon}^{P^-}\,\V_{c_p\epsilon}^{U^+} \,z\;
\V_{c_u\epsilon}^{U^-}\,z^{-1}\,z\,\V_{c_u\epsilon}^{G''}\,\V_{c_u\epsilon}^{U^+}\,u
%\\ &
\subset p\,\V_{c_p\epsilon}^{P^-}\,\V_{c_p\epsilon}^{U^+} \,
\V_{c_u\epsilon}^{U^-}\,\V_{c_u\epsilon}^{G''}\,z\;\V_{c_u\epsilon}^{U^+}
\,u \\ &\subset
p\,\V_{c_p\epsilon}^{P^-}\,\V_{c_u\epsilon}\,\V_{c_p\epsilon}^{U^+} \,
\,z\;\V_{c_u\epsilon}^{U^+} \,u
=p\,\V_{c_p\epsilon}^{P^-}\,\V_{c_u\epsilon}^{P^-}\V_{c_u\epsilon}^{U^+}\,
\V_{c_p\epsilon}^{U^+}\, \,z\;\V_{c_u\epsilon}^{U^+} \,u
\\ &=p\,\V_{c_p\epsilon}^{P^-}\,\V_{c_u\epsilon}^{P^-}\,z\,z^{-1}\,
\V_{c_u\epsilon}^{U^+}\,z\,z^{-1}\,\V_{c_p\epsilon}^{U^+}\, \,z\;
\V_{c_u\epsilon}^{U^+} \,u
\\ &\subset p\,\V_{c_p\epsilon}^{P^-}\,\V_{c_u\epsilon}^{P^-}\,z\,
\V_{c_u\epsilon}^{U^+}\,\V_{c_p\epsilon}^{U^+}\,\V_{c_u\epsilon}^{U^+} \,u
\subset p\,\V_{c''_1\epsilon}^{P^-}\,z\, \V_{c''_2\epsilon}^{U^+} \,u
\\ &=u^-\,g''\,\V_{c''_1\epsilon}^{U^-}\,\V_{c''_1\epsilon}^{G''}\,z\,
\V_{c''_2\epsilon}^{U^+} \,u
=u^-\,g''\,\V_{c''_1\epsilon}^{U^-}\,{g''}^{-1}\,g''\,
\V_{c''_1\epsilon}^{G''}\,z\, \V_{c''_2\epsilon}^{U^+} \,u
\\ &\subset u^-\,\V_{c_{g''}c''_1\epsilon}^{U^-}\;g''\,
\V_{c''_1\epsilon}^{G''}\,z\, \V_{c''_2\epsilon}^{U^+} \,u
\;,
\end{align*}
as wanted. \cqfd

\medskip
Let $\Phi$ be a closed ball of radius less than $1$ in the metric
space $\Gr_{\d,\D}$, contained in $\Gr_{\d,\D}^\flat$. By Lemma
\ref{lem:isometrygrass} and with the notation of Equation
\eqref{eq:defiwtorb}, the set $\wt\Phi= \operatorname{orb}_\d^{\;-1}
(\Phi)$ is a closed ball of same radius in $\M_{\dd,\d} (\OOO_\val)$.
Let $\E$ be a closed ball in $\La_{\d,\dd}$, small enough so that
there exists a closed ball in the clopen fundamental domain
$G''_\lozenge$ which maps isometrically to $\E$ by the locally
isometric map $\pi_{\d,\dd}: G''\ra \La_{\d,\dd}$ defined in Equation
\eqref{eq:defiprok}. Let $\wt \E=\pi_{\d,\dd}^{-1}(\E)$ and let $r$
vary in $\NN$. Using the notation $\Omega_{\Phi,\E,r}$ defined in
Equation \eqref{eq:defiOmega}, the family of Lipschitz well-rounded
subsets with respect to $(\V_\epsilon)_{\epsilon>0}$ that we will use
in order to apply Theorem \ref{theo:GorodnikNevo} is given by the
following result.

\bprop\label{prop:verifLWR} With $\Phi$ and $\E$ as above, the family
$\big(\Omega_{\Phi,\E,r} \big)_{r\in\NN}$ is $0$-Lipschitz
well-rounded with respect to $(\V_\epsilon)_{\epsilon>0}$.
\eprop

\dem Recall that $\Omega_{\Phi,\E,r}=U^-_{\wt\Phi}\; G''_{\wt\E}\;
Z_r\;U^+_{\Box}$ with the notation at the beginning of Subsection
\ref{subsect:correspondance}. We will actually prove (as allowed by the
ultrametric situation) the stronger statement that given $\Phi$ and
$\E$ as above, if $\epsilon$ is small enough, then for every
$r\in\NN$, we have
\[
\big(U^-_{\wt\Phi}\, G''_{\wt \E}\, Z_r\,
U^+_{\Box}\big)^{-\epsilon}=U^-_{\wt\Phi}\, G''_{\wt\E}\, Z_r\,
U^+_{\Box}=\big(U^-_{\wt\Phi}\, G''_{\wt\E}\, Z_r\,
U^+_{\Box}\big)^{+\epsilon}\;.
\]
Let
\begin{equation}\label{eq:defic}
c=\max\big\{q_\val \max\{c_{g''}(c_{u^-g''}+c_{u^+}), c_{u^-g''}+2\,c_{u^+}\} :
u^-\in U^-_{\wt\Phi}, \;g''\in G''_{\wt\E}, \;u^+\in U^+_{\Box}\big\}\,,
\end{equation}
which is finite since $U^-_{\wt\Phi}$, $G''_{\wt\E}$ and $U^+_{\Box}$ are
compact subsets of $G$. Since $\wt\Phi$ is a ball of radius less than $1$
in $\M_{\dd,\d} (K_\val)$, let $v_{0}\in\M_{\dd,\d}(K_\val)$ and
$k\in\NN\!\smallsetminus\!\{0\}$ be such that
\[
\wt\Phi=v_{0}+\pi_{\val}^{\;k}\M_{\dd,\d}(\OOO_\val)\;.
\]
Let $r_\E$ be the radius of the ball $\E$ (satisfying the assumptions
of Proposition \ref{prop:verifLWR}). By Equation \eqref{eq:diagodist}
and since the map $\underline{g} \mapsto\widecheck{\underline{g}}$ is
an isometry of $\GL_\dd(K_\val)$, there exists $\big(\begin{smallmatrix}
\ov{g}_0 &0 \\ 0 & \underline{g}_0\end{smallmatrix}\big)\in G''$
such that
\begin{equation}\label{eq:explicitGsecwtE}
G''_{\wt\E}=\big\{\big(\begin{smallmatrix}
  \ov{g} &0 \\ 0 & \underline{g}\end{smallmatrix}\big)\in G'':
\max\{d(\ov{g},\ov{g}_0),\;d(\underline{g},\underline{g}_0)\}
\leq r_\E\big\}\;.
\end{equation}
%$x_{0}\in\GL^1_\d(K_{\val})$, $y_{0}\in\GL^1_\dd(K_{\val})$,
%$k,\ell,\ell'\in\ZZ$ with $\ell,\ell'\geq 1$ be such that
%\[
%\Psi=v_{0}+\pi_{\val}^{k}\M_{\dd,\d}(\OOO_\val), \quad
%\X=x_{0}+\pi_{\val}^{\ell}\M_{\d}(\OOO_\val), \quad
%\Y=y_{0}+\pi_{\val}^{\ell'}\M_{\dd}(\OOO_\val)\;.
%\]
%Note that we have $\|x_0\|\geq 1$ since otherwise by the ultrametric
%inequality and the formula of the determinant, we would have $|\det
%x_0|<1$. Similarly we have $\|y_0\|\geq 1$.  Since $\ell,\ell'>0$, we
%have $\X\subset \GL^1_\d(K_{\val})$ and $\Y\subset \GL^1_\dd
%(K_{\val})$. Hence $\X=\GL^1_\d(K_{\val})\cap B(x_0,q_\val^{\;-\ell})$
%and $\Y= \GL^1_\dd (K_{\val})\cap B(y_0,q_\val^{\;-\ell'})$, as wanted
%before the statement of Proposition \ref{prop:verifLWR}.
%
%Since $\E$ is compact, and since the fundamental domain
%$G''_\lozenge$ defined at the beginning of Subsection
%\ref{subsect:correspondance} is closed-open, there exists $k'\in\NN$
%such that for all $g\in G''_{\wt\E}$, $\alpha'\in\M_{\d}(\OOO_\val)$ and
%$\delta'\in\M_{\dd}(\OOO_\val)$, we have $g\big(\begin{smallmatrix}
%  I_\d+\pi_\val^{\;k'}\alpha' &0 \\ 0 & I_\dd+\pi_\val^{\;k'}\delta'
%\end{smallmatrix}\big)\in G''_\lozenge$.

Let us now consider $\epsilon_0= \frac{1}{c}\, q_\val^{-\max\{k,\,
r_\E+\log_{q_\val}(\frac{1}{r_\E}\|\ov{g}_0\|\,\|{\ov{g}_0}^{-1}\|),\,
r_\E+\log_{q_\val}(\frac{1}{r_\E}\|\underline{g}_0\|\,
\|{\underline{g}_0}^{-1}\|)\} -2} >0$, so that for every
$\epsilon\in\;]0, \epsilon_0[\,$, we have
\begin{equation}\label{eq:minoNceps}
N_{c\epsilon}> 1+ \max\big\{k,
r_\E+\log_{q_\val}\big(\frac{1}{r_\E}\|\ov{g}_0\|\,\|{\ov{g}_0}^{-1}\|\big),
r_\E+\log_{q_\val}\big(\frac{1}{r_\E}\|\underline{g}_0\|\,
\|{\underline{g}_0}^{-1}\|\big)\big\}\geq 1\;.
\end{equation}
Let us prove that for every $\epsilon\in\;]0,
%\min\{\epsilon_0,\frac{1}{q_\val}\}
\epsilon_0[\,$, we have
\begin{equation}\label{eq:perturbations}
  U_{\wt\Phi}^{-}\;\V_{c\epsilon}^{U^{-}}=U_{\wt\Phi}^{-},\quad
  G''_{\wt\E}\;\V_{c\epsilon}^{G''}=G''_{\wt\E}\quad\text{and}\quad
  \V_{c\epsilon}^{U^{+}}\,U_{\Box}^{+}=U_{\Box}^{+}\;.
\end{equation}

For all $u\in U_{\wt\Phi}^{-}$ and $u'\in \V_{c\epsilon}^{U^{-}}$, let
$\beta,\beta'\in\M_{\dd,\d}(\OOO_\val)$ be such that $u=\big(
\begin{smallmatrix}  I_{\d} & 0\\ v_0+\pi_\val^{\;k}\beta\; & I_{\dd}
\end{smallmatrix}\big)$ and $u'=\big(\begin{smallmatrix} I_{\d}
&\; 0\\ \pi_\val^{\;N_{c\epsilon}}\beta'\; & I_{\dd}\end{smallmatrix}
\big)$. Then since $N_{c\epsilon}>k$, we have $uu'=\big(
\begin{smallmatrix} I_{\d} &0 \\ v_0+\pi_\val^{\;k}\beta+
\pi_\val^{\;N_{c\epsilon}}\beta'\; & I_{\dd} \end{smallmatrix}\big)\in
U_{\wt\Phi}^{-}$. Therefore we have $U_{\wt\Phi}^{-}\;\V_{c\epsilon}^{U^{-}}
\subset U_{\wt\Phi}^{-}$ and the opposite inclusion is clear.  This
proves the equality on the left-hand side of Formula
\eqref{eq:perturbations}.

The proof of the equality on the right-hand side is similar.  For all
$u\in U_{\Box}^{+}$ and $u'\in \V_{c\epsilon}^{U^{+}}$, let
$\ga\in\Box$ and $\ga'\in\M_{\d,\dd} (\OOO_\val)$ be such that
$u=\big(\begin{smallmatrix} I_{\d} &\; \ga\\ 0 &
  I_{\dd}\end{smallmatrix}\big)$ and $u'=\big(
\begin{smallmatrix} I_{\d} &\; \pi_\val^{\;N_{c\epsilon}}\ga'\\ 0 & I_{\dd}
\end{smallmatrix}\big)$. Then since $N_{c\epsilon}>1$ and since
$\Box+\pi_\val\M_{\d,\dd}(\OOO_\val)=\Box$ by the construction of
the fundamental domain ${\mathcal D}$ at the beginning of Subsection
\ref{subsect:correspondance}, we have $u'u=\big(\begin{smallmatrix}
  I_{\d} &\; \ga+ \pi_\val^{\;N_{c\epsilon}}\ga'\\ 0 &
  I_{\dd} \end{smallmatrix}\big)\in U_{\Box}^{+}$. Therefore we have
$\V_{c\epsilon}^{U^{+}}U_{\Box}^{+} \subset U_{\Box}^{+}$ and the
opposite inclusion is clear.

Let $g=\Big(\begin{smallmatrix} \overline{g}& 0\\ 0 &
  \underline{g} \end{smallmatrix}\Big)\in G''_{\wt\E}$. For every
$g'=\Big(\begin{smallmatrix} \overline{g'}& 0\\ 0 &
  \underline{g'} \end{smallmatrix}\Big)\in \V_{c\epsilon}^{G''}$,
there exist $\alpha\in\M_{\d}(\OOO_\val)$ and $\delta\in
\M_{\dd}(\OOO_\val)$ such that
\[
\quad\overline{g'}=I_{\d} +\pi_\val^{\;N_{c\epsilon}}\alpha\quad \text{and}\quad
\underline{g'}= I_{\dd}+\pi_\val^{\;N_{c\epsilon}}\delta\;.
\]
We have $gg'=\Big(\begin{smallmatrix} \overline{g}\,\overline{g'}&
  0\\ 0 & \underline{g}\,\underline{g'} \end{smallmatrix}\Big)\in
G''$. Let us prove that $d(\overline{g}\,\overline{g'},\,
\overline{g}_0) \leq r_\E$. A similar proof gives that
$d(\underline{g} \,\underline{g'},\,\underline{g}_0)\leq r_\E$, thus
proving that $G''_{\wt\E}\; \V_{c\epsilon}^{G''} \subset G''_{\wt\E}$
by Equation \eqref{eq:explicitGsecwtE}. The opposite inclusion being
clear, this proves the middle equality of Formula
\eqref{eq:perturbations}.

By the submultiplicativity of the supremum norm, since $\alpha\in
\M_\d(\OOO_\val)$ so that $\|\alpha\|\leq 1$, by Lemma
\ref{lem:controlnormsupdist} since $\overline{g}\in B(\overline{g}_0,
r_\E)$ and by Equation \eqref{eq:minoNceps}, we have
\[
\|\pi_\val^{\;N_{c\epsilon}}\,\overline{g}\,\alpha\,{\overline{g}_0}^{-1}\|
\leq q_\val^{\;-N_{c\epsilon}}\|\overline{g}\|\,\|\alpha\|\,
\|{\overline{g}_0}^{-1}\|\leq q_\val^{\;-N_{c\epsilon}+r_\E}\|\overline{g}_0\|\,
\|{\overline{g}_0}^{-1}\|\leq \frac{r_\E}{q_\val}\;.
\]
We also have, by the ultrametric triangle inequality,
\begin{align*}
\|\,\overline{g}\,\overline{g'}\,{\overline{g}_0}^{-1}-I_\d\|&=
\|\,\overline{g}(I_{\d} +\pi_\val^{\;N_{c\epsilon}}\alpha)\,
{\overline{g}_0}^{-1}-I_\d\|=
\|\,\overline{g}\,{\overline{g}_0}^{-1}-I_\d+\pi_\val^{\;N_{c\epsilon}}
\,\overline{g}\,\alpha\,{\overline{g}_0}^{-1}\|\\&\leq
\max\{\|\,\overline{g}\,{\overline{g}_0}^{-1}-I_\d\|,\;
\|\pi_\val^{\;N_{c\epsilon}}\,\overline{g}\,
\alpha\,{\overline{g}_0}^{-1}\|\}\;.
\end{align*}
Thus since $\overline{g}\in B(\overline{g}_0,r_\E)$ and $\ln(1+t)\leq
t$ for every $t\geq 0$, we have
\[
\log_{q_\val}(1+ \|\,\overline{g}\,\overline{g'}\,{\overline{g}_0}^{-1}
-I_\d\|)\leq \max\{r_\E,\log_{q_\val}(1+\frac{r_\E}{q_\val})\}\leq r_\E
\;.
\]
Since $N_{c\epsilon}\geq 1$, the standard formula for the inverse of
$I_\d+X$ when $\|X\|<1$ gives that there exists $\alpha'\in
\M_\d(\OOO_\val)$ such that $\overline{g'}^{\;-1}=I_{\d}
+\pi_\val^{\;N_{c\epsilon}} \alpha'$. A proof similar to the one above
thus gives that $\log_{q_\val}(1+ \|\,{\overline{g}_0}(\,\overline{g}\,
\overline{g'}\, )^{-1} -I_\d\|)\leq r_\E$, which proves as wanted that
$d(\overline{g}\,\overline{g'},\, \overline{g}_0) \leq r_\E$.  This
concludes the proof of Formula \eqref{eq:perturbations}.

Now, for every $r\in\NN$, we have by Lemma \ref{eq:effectiveLU} and by
Equations  \eqref{eq:defic} and \eqref{eq:perturbations} that
\begin{align*}
&\big(U^-_{\wt\Phi}\, G''_{\wt\E}\, Z_r\,
U^+_{\Box}\big)^{+\epsilon}=\V_\epsilon\,
U^-_{\wt\Phi}\, G''_{\wt\E}\, Z_r\,
U^+_{\Box}\,\V_\epsilon\\\subset\;\;& U^-_{\wt\Phi}
\,\V^{U^-}_{c\epsilon}\,G''_{\wt\E}\,\V^{G''}_{c\epsilon}\,
Z_r\,\V^{U^+}_{c\epsilon}\,U^+_{\Box}=
U^-_{\wt\Phi}\, G''_{\wt\E}\, Z_r\, U^+_{\Box}\;.
\end{align*}
Since the converse inclusion is immediate, we have $\big(U^-_{\wt\Phi}\,
G''_{\wt\E}\, Z_r\, U^+_{\Box}\big)^{+\epsilon}= U^-_{\wt\Phi}\,
G''_{\wt\E}\, Z_r \, U^+_{\Box}$.

Since $\V_\epsilon$, being a subgroup, is stable by $g\mapsto g^{-1}$,
this implies that $g\; U^-_{\wt\Phi}\, G''_{\wt\E}\, Z_r \,
U^+_{\Box}\;h$ contains $U^-_{\wt\Phi}\, G''_{\wt\E}\, Z_r \, U^+_{\Box}$
for all $g,h\in \V_\epsilon$, so that $ \big(U^-_{\wt\Phi}\,
G''_{\wt\E}\, Z_r \, U^+_{\Box}\big)^{-\epsilon} \supset U^-_{\wt\Phi}\,
G''_{\wt\E}\, Z_r \, U^+_{\Box}$.  Since the converse inclusion is
immediate, this concludes the proof of Proposition
\ref{prop:verifLWR}.  \cqfd

\subsection{The main statement and its proof}
\label{subsect:proofmain}

Error terms in equidistribution results usually require smoothness
properties on test functions. The appropriate smoothness regularity of
functions defined on ultrametric spaces as $\Gr_{\d,\D}$ and
$\La_{\d,\dd}$ is the locally constant one. The locally constant
regularity on such homogeneous spaces of totally discontinuous groups
could be defined (as for instance in \cite{AthGhoPra12}, \cite[\S
  4.3]{KemPauSha17}) by using the familly of small compact-open
subgroups $(\V_\epsilon)_{\epsilon\,\in\, ]0,1]}$ of $G$ defined in
Subsection \ref{subsect:wellround}, and by defining an
$\epsilon$-locally constant map on $\La_{\d,\dd}$ to be a map which is
constant on every orbit of $\V_\epsilon\cap G''$ on $\La_{\d,\dd}$.
But it turns out to be more convenient in this paper to use a general
purely metric definition. For every ultrametric space $E$ and
$\epsilon \in\;]0,1]$, a bounded map $f : E \ra \RR$ is {\it
$\epsilon$-locally constant} if it is constant on every closed
ball of radius $\epsilon$ in $E$. With $\|f\|_\infty=\sup_{x\in E}
|f(x)|$ the supremum norm of $f$, the {\it $\epsilon$-locally constant
  norm} of $f$ is $\|f\|_\epsilon= \frac{\|f\|_\infty}{\epsilon}$.

\bigskip
The key result of this paper is the following one. Let $\ell=\lcm(\d,
\dd)$.  For every nonzero ideal $I$ of $R_\val$, let
$\PL_{\d,\D}^\sharp(I)=\Ga_I^\sharp(R_\val^{\;\d}\times\{0\})$ (see
Subsection \ref{subsect:congruence} for the definition of $\Ga_I$ and
Equation \eqref{eq:defiHsharp} for the one of $\Ga_I^\sharp$) and
\begin{equation}\label{eq:defcsubI}
  c_I= q^{\;(\ggg-1)(\D^2-1-\d\,\dd)}\,\prod_{i=1}^{\D-1}
  \frac{\zeta_K(i+1)}{q_\val^{\;i}-1}\,N(I)^{\d\,\dd}\prod_{\ppp\,|I}
  \prod_{i=1}^\d\frac{N(\ppp)^{i}-N(\ppp)^{-\dd}}{N(\ppp)^i-1}\;.
\end{equation}

\btheo\label{theo:mainsharp} For every nonzero ideal $I$ of $R_\val$,
for the weak-star convergence of Borel measures on the locally compact
space $\Gr_{\d,\D}^\flat\times \La_{\d,\dd}$, we have
\begin{equation}\label{eq:theomainsharp}
\lim_{i\ra+\infty}\;\frac{c_I}{q_\val^{\;\ell\,\D\,i}}
\sum_{\Lambda\in\PL_{\d,\D}^\sharp(I)\;:\;
  \overline{\covol}\;\Lambda\,=\,q_\val^{\;\ell\,i}}\;
\Delta_{V_\Lambda}\otimes\Delta_{\llbracket\Lambda\rrbracket} =
{\mu_{\Gr_{\d,\D}}}_{\mid\Gr_{\d,\D}^\flat} \otimes\mu_{\La_{\d,\dd}}\;.
\end{equation}
Furthermore, there exists $\tau\in\;]0,\frac{1}{2\,\D^2}]$ such that
for all $\delta\in\;]0,\tau[$ and $\epsilon\in\;]0,1]$, there is an
additive error term of the form $\bigO_{\val,\delta,I}\big(
q_\val^{\;\ell\,\D\,i (-\tau+\delta)} \,\|f\|_\epsilon\,
\|g\|_\epsilon \big)$ in the above equidistribution claim when
evaluated on pairs $(f,g)$ for all compactly supported $\epsilon$-locally
constant maps $f:\Gr_{\d,\D}^\flat \ra\RR$ and
$g:\La_{\d,\dd}\ra\RR\;$: as $i\ra+\infty$, we have 
\begin{align*}
  &\frac{c_I}{q_\val^{\;\ell\,\D\,i}}
\sum_{\Lambda\in\PL_{\d,\D}^\sharp(I)\;:\;\overline{\covol}\;\Lambda\,=\,q_\val^{\;\ell\,i}} 
f(V_\Lambda)\;g(\llbracket\Lambda\rrbracket)\\ =\;& 
\Big(\int_{\Gr_{\d,\D}^\flat} f\,d\mu_{\Gr_{\d,\D}}\Big)
\Big(\int_{\La_{\d,\dd}}g\,d\mu_{\La_{\d,\dd}}\Big)+\bigO_{\val,\delta,I}
\big(q_\val^{\;\ell\,\D\,i(-\tau+\delta)} \,\|f\|_\epsilon\,\|g\|_\epsilon\big)\;.
\end{align*}
\etheo

\dem Let $\Lambda\in\PL_{\d,\D}^\sharp$. Then there exists $g\in
\Ga^\sharp$ (thus $g\in \U_G$) such that we have $\Lambda=\Lambda_g=
g(R_\val^{\;\d}\times\{0\})$, so that if $g=u^-\,g''\,z\,u^+$ is the
decomposition given by Proposition \ref{prop:refLUbloc}, then $u^-\in
U^-(\OOO_\val)$. Hence by Proposition \ref{prop:HKprop4_3}
\hyperlink{HKprop4_3i}{(i)}, we have
\[
V_\Lambda=V_{\Lambda_g}=
V_{\Lambda_{u^-}} \in\Gr_{\d,\D}^\flat=U^-(\OOO_\val)
V_{R_\val^{\;\d}\times\{0\}}\;.
\]
Furthermore, we have $\llbracket\Lambda \rrbracket\in\La_{\d,\dd}$ by
Equation \eqref{eq:assoccorrpairlat}, so that the statement of Theorem
\ref{theo:mainsharp} is well defined.

Let $I$ be a nonzero ideal of $R_\val$. Let $\tau= \tau(\Ga_I)
\in\;]0,\frac{1}{2\,\D^2}]$ be as in Theorem \ref{theo:GorodnikNevo}
applied with $G'=G$, with $\Ga'=\Ga_I$ and with the family
$(\V'_\epsilon)_{\epsilon>0} = (\V_\epsilon)_{\epsilon>0}$ given by
Equation \eqref{eq:defiVeps}, which has an upper local dimension
$\rho=\D^2-1$ according to the final claim of Lemma
\ref{eq:calVdivgroup}.  Let $\delta\in\;]0,\tau[$.

Let $\Phi$ be a closed ball  in $\Gr_{\d,\D}$ of radius $r_\Phi\in\;]0,1]$.
Besides, we assume that $\Phi$ is contained in $\Gr_{\d,\D}^\flat$.
By Lemma \ref{lem:isometrygrass} and with the notation of Equation
\eqref{eq:defiwtorb}, it follows that $\wt\Phi=
\operatorname{orb}_\d^{\;-1}(\Phi)$ is a closed ball of radius
$r_\Phi$ in $\M_{\dd,\d}(\OOO_\val)$. Let $\chi_\Phi$ be the
characteristic function of $\Phi$, which is $r_\Phi$-locally constant
with $\|\chi_\Phi\|_{r_\Phi} = \frac{1}{r_\Phi}\geq 1$. Since
$\mu_{\Gr_{\d,\D}}$ is a probability measure and by Equation
\eqref{eq:volGUmoins}, we have
\[
1\geq \mu_{\Gr_{\d,\D}}(\Phi)=
c_1\,\operatorname{Haar}_{\dd,\d}(\wt\Phi)=c_1\,r_\Phi^{\;\d\,\dd}\;,
\]
so that, since $\tau\leq \frac{1}{2\,\D^2}\leq \frac{1}{\d\,\dd}$, we have
\begin{equation}\label{eq:lcnormchiPhi}
\mu_{\Gr_{\d,\D}}(\Phi)^{-\tau +\delta} \leq
\mu_{\Gr_{\d,\D}}(\Phi)^{-\tau}\leq(c_1\,r_\Phi^{\;\d\,\dd})^{-\frac{1}{\d\,\dd}}
= \bigO(\|\,\chi_\Phi\|_{r_\Phi})\;.
\end{equation}

Let $\E$ be a closed ball in $\La_{\d,\dd}$ of radius $r_\E\in\;]0,1]$
small enough so that $\mu_{\La_{\d,\dd}}(\E)\leq 1$ and there exists a
closed ball $\wt \E_0$ in $G''_\lozenge$ mapping isometrically to $\E$
by $\pi_{\d,\dd}:G''\ra\La_{\d,\dd}$. Let $\wt \E=\pi_{\d,\dd}^{\;-1}
(\E) = \bigsqcup_{\ga\in G''(R_\val)} \wt \E_0\,\ga $. Let $\chi_\E$
be the characteristic function of $\E$, which is $r_\E$-locally
constant with $\|\chi_\E \|_{r_\E} = \frac{1}{r_\E}\geq 1$.  By
Equation \eqref{eq:volGEF} and by the Alhfors regularity of the
homogeneous measure $\mu_{G''}$ of the group $G''$ with dimension
$\dim G''\leq \d^2+\dd^2\leq 2\D^2$ for the distance $d$ defined in
Section \ref{subsect:latticepair} (see in particular Equation
\eqref{eq:defidistGLk}), there exists a constant $c>0$ such that
\[
\mu_{\La_{\d,\dd}}(\E)
= \mu_{G''}(G''_{\wt\E})
= \mu_{G''}(\wt\E_0)
\geq c\, r_\E^{\;\dim G''}\geq c\, r_\E^{\;2\D^2}\;,
\]
so that, since $\tau\leq \frac{1}{2\,\D^2}$ and $\mu_{\La_{\d,\dd}}
(\E) \leq 1$, we have
\begin{equation}\label{eq:lcnormchiE}
\mu_{\La_{\d,\dd}}(\E)^{-\tau +\delta}\leq\mu_{\La_{\d,\dd}}(\E)^{-\tau}
\leq(c\, r_\E^{\;2\D^2})^{-\frac{1}{2\,\D^2}}=
\bigO(\|\,\chi_\E\|_{r_\E})\;.
\end{equation}

For every $r\in\NN$, let us define
\[
\PL_{\d,\D}^\sharp(I,\Phi,\E,r)=
\{\Lambda\in\PL_{\d,\D}^\sharp(I):V_\Lambda\in \Phi,\;\llbracket
\Lambda \rrbracket\in\E, \;\overline{\covol}
\;\Lambda\, =\,q_\val^{\;\ell\,r}\}\;
\]
Using respectively

$\bullet$~ Proposition \ref{prop:matrepprimlat} (where the set
$\Omega_{\Phi,\E,r}$ is defined in Equation \eqref{eq:defiOmega}) for
the first equality,

$\bullet$~ Theorem \ref{theo:GorodnikNevo} applied to the family
$(\Z_r= \Omega_{\Phi,\E,r})_{r\in\NN}$, which is $0$-Lipschitz
well-rounded with respect to $(\V_\epsilon)_{\epsilon>0}$ by
Proposition \ref{prop:verifLWR} for the second equality,

$\bullet$~ Lemma \ref{lem:calcmesOmegasub} for
the third equality,

$\bullet$~ Equations \eqref{eq:lcnormchiPhi} and \eqref{eq:lcnormchiE}
for the last equality (and the fact that $\mu_{\Gr_{\d,\D}}$ and
$\mu_{\La_{\d,\dd}}$ are finite measures),

\noindent we have
\begin{align}
  &\card \PL_{\d,\D}^\sharp(I,\Phi,\E,r)
  =\card (\Omega_{\Phi,\E,r}\cap\Ga_I)\nonumber
  \\=\;&\frac{1}{\|\mu_{G/\Ga_I}\|}\,\mu_{G}(\Omega_{\Phi,\E,r})
  \big(1+ \bigO\big(\mu_{G}(\Omega_{\Phi,\E,r})^{-\tau+\delta}
  \big)\big)\nonumber
  \\=\;&\frac{q^{(\ggg-1)\,\d\,\dd}\;q_\val^{\;\D\,\ell\,r}}{\|\mu_{G/\Ga_I}\|}\,
  \mu_{\Gr_{\d,\D}}(\Phi)\,\mu_{\La_{\d,\dd}}(\E)\big(1+
  \bigO\big(q_\val^{\;\D\,\ell\,r(-\tau+\delta)}
  \mu_{\Gr_{\d,\D}}(\Phi)^{-\tau+\delta}
  \mu_{\La_{\d,\dd}}(\E)^{-\tau+\delta}\big)\big)\nonumber
\\=\;&\frac{q^{(\ggg-1)\,\d\,\dd}\;q_\val^{\;\D\,\ell\,r}}{\|\mu_{G/\Ga_I}\|}\,
  \Big(\mu_{\Gr_{\d,\D}}(\Phi)\,\mu_{\La_{\d,\dd}}(\E)+
  \bigO\big(q_\val^{\;\D\,\ell\,r(-\tau+\delta)}\;
  \|\,\chi_\Phi\|_{r_\Phi}\;\|\,\chi_\E\|_{r_\E}\,\big)\Big)\;.
  \label{eq:maincasboule}
\end{align}
Let $c_I=\frac{\|\mu_{G/\Ga_I}\|}{q^{(\ggg-1)\,\d\,\dd}}=
\frac{[\Ga:\Ga_I]\, \|\mu_{G/\Ga}\|} {q^{(\ggg-1)\,\d\,\dd}}$. With
the value of $[\Ga: \Ga_I]$ given by Lemma \ref{lem:indexHeckesubgrou}
(2) and the value of $\|\mu_{G/\Ga}\|$ given by Equation
\eqref{eq:totmassmuSLmodSL} with $k=\D$, we have
\begin{align*}
  c_I=
  \frac{N(I)^{\d\,\dd}\prod_{\ppp\,|I}\prod_{i=1}^\d
  \frac{N(\ppp)^{i}-N(\ppp)^{-\dd}}{N(\ppp)^i-1}\,q^{\;(\ggg-1)(\D^2-1)}
  \prod_{i=1}^{\D-1}\frac{\zeta_K(1+i)}{q_\val^{\;i}-1}}
       {q^{(\ggg-1)\,\d\,\dd}}\;,
\end{align*}
as wanted in Equation \eqref{eq:defcsubI}.  Note that every compactly
supported $\epsilon$-locally constant map on an ultrametric space is a
finite linear combination of characteristic functions of balls of
radius $\epsilon$.  By a finite bilinearity argument, Theorem
\ref{theo:mainsharp} follows from Equation \eqref{eq:maincasboule}. \cqfd

\bcoro\label{coro:mainshape} For every nonzero ideal $I$ of $R_\val$,
for the weak-star convergence of Borel measures on the locally compact
space $\Gr_{\d,\D}\times \Sh_{\d}\times \Sh_{\dd}$, we have
\begin{align}
\lim_{i\ra+\infty}\;\frac{c_I\,(q-1)}{q_\val^{\;\ell\,\D\,i}}&
\sum_{\Lambda\in\PL_{\d,\D}(I)\;:\;
  \overline{\covol}\;\Lambda\,=\,q_\val^{\;\ell\,i}}\; \Delta_{V_\Lambda}
\otimes\Delta_{\sh(\Lambda)}\otimes\Delta_{\sh(\Lambda^\perp)}\nonumber
\\&= \mu_{\Gr_{\d,\D}} \otimes\mu_{\Sh_\d}\otimes
\mu_{\Sh_\dd}\;. \label{eq:coromainsharp}
\end{align}
Furthermore, there exists $\tau\in\;]0,\frac{1}{2\,\D^2}]$ such that
for all $\delta\in\;]0,\tau[$ and $\epsilon\in\;]0,1]$, there is an
additive error term of the form $\bigO_{\val,\delta,I}\big(
q_\val^{\;\ell\,\D\,i (-\tau+\delta)} \,\|f\|_\epsilon\,
\|f_1\|_\infty \|f_2\|_\infty\big)$ in the above equidistribution
claim when evaluated on $(f,f_1,f_2)$ for every compactly supported
$\epsilon$-locally constant map $f:\Gr_{\d,\D} \ra\RR$ and for
all finitely supported maps $f_1:\Sh_{\d}\ra\RR$ and $f_2:\Sh_{\dd}
\ra\RR\;$: as $i\ra+\infty$, we have
\begin{align*}
  &\frac{c_I\,(q-1)}{q_\val^{\;\ell\,\D\,i}}
\sum_{\Lambda\in\PL_{\d,\D}(I)\;:\;\overline{\covol}\;\Lambda\,=\,q_\val^{\;\ell\,i}} 
f(V_\Lambda)\;f_1(\sh(\Lambda))\;f_2(\sh(\Lambda^\perp))\\ =\;& 
\Big(\int_{\Gr_{\d,\D}} f\,d\mu_{\Gr_{\d,\D}}\Big)
\Big(\int_{\Sh_{\d}}f_1\,d\mu_{\Sh_{\d}}\Big)
\Big(\int_{\Sh_{\dd}}f_2\,d\mu_{\Sh_{\dd}}\Big)\\ & +\bigO_{\val,\delta,I}
\big(q_\val^{\;\ell\,\D\,i(-\tau+\delta)}
\,\|f\|_\epsilon\,\|f_1\|_\infty\|f_2\|_\infty\big)\;.
\end{align*}
\ecoro

Theorem \ref{theo:mainintrotriple} in the Introduction follows from
the first claim of this corollary by taking $I=R_\val$ and
$c'=c_{R_\val}(q-1)$, and by using Equation \eqref{eq:defcsubI}.

\medskip
\dem {\bf Step 1. } We first prove the result with
$\Gr_{\d,\D}^\flat$ instead of $\Gr_{\d,\D}$ and
$\PL_{\d,\D}^\sharp(I)$ instead of $\PL_{\d,\D}(I)$.

Since the map $\varphi_{\d,\dd}$ (defined in Equation
\eqref{eq:defivarphimn}) is proper by Lemma \ref{lem:proprivarphi},
the pushforward map $(\varphi_{\d,\dd})_*$ of Borel measures by
$\varphi_{\d,\dd}$ is linear and weak-star continuous.  Hence applying
the map $(\id\times\varphi_{\d,\dd} )_*$ to Equation
\eqref{eq:theomainsharp}, using Lemma \ref{lem:proprivarphi}
\eqref{item2:proprivarphi} on the left hand side of Equation
\eqref{eq:theomainsharp}, and Lemma \ref{lem:proprivarphi}
\eqref{item1:proprivarphi} on the right hand side of Equation
\eqref{eq:theomainsharp}, we have
\begin{align*}
\lim_{i\ra+\infty}\;&\frac{c_I}{q_\val^{\;\ell\,\D\,i}}
\sum_{\Lambda\in\PL_{\d,\D}^\sharp(I)\;:\;
  \overline{\covol}\;\Lambda\,=\,q_\val^{\;\ell\,i}}\;\Delta_{V_\Lambda}
\otimes\Delta_{\sh(\Lambda)}\otimes\Delta_{\sh(\Lambda^\perp)}
\\& =(q-1)\,{\mu_{\Gr_{\d,\D}}}
_{\mid\Gr_{\d,\D}^\flat} \otimes\mu_{\Sh_\d}\otimes \mu_{\Sh_\dd}\;.
\end{align*}
It follows from Lemma \ref{lem:proprivarphi}
\eqref{item3:proprivarphi} and from the error term in Theorem
\ref{theo:mainsharp} applied with the compactly supported
$\rho_0$-locally constant function $g=(f_1\times f_2)\circ
\varphi_{\d,\dd}$ that we have an additive error term of the form
$\bigO_{\val,\delta,I}\big( q_\val^{\;\ell\,\D\,i (-\tau+\delta)}
\,\|f\|_\epsilon\, \|f_1\|_\infty \|f_2\|_\infty\big)$ in this
equidistribution claim when evaluated on $(f,f_1,f_2)$ for every
compactly supported $\epsilon$-locally constant map
$f:\Gr_{\d,\D}^\flat \ra\RR$ and for all finitely supported maps
$f_1:\Sh_{\d}\ra\RR$ and $f_2:\Sh_{\dd}\ra\RR\;$.

\medskip \noindent {\bf Step 2. } We now explain how to deduce
Corollary \ref{coro:mainshape} from the equidistribution of
$(V_\Lambda, \sh(\Lambda),\sh(\Lambda^\perp))$ in $\Gr_{\d,\D}^\flat
\times \Sh_\d\times \Sh_\dd$ when $\Lambda$ varies in
$\PL_{\d,\D}^\sharp(I)$ with the appropriate covolume. The key
technical lemma is the following one.

Let us denote by $W_\D$ the Weyl subgroup of $\GL_\D(K_\val)$
consisting in the permutation matrices of the canonical basis of
$K_\val^{\;\D}$. Note that $W_\D$ is contained in $\GL_\D(R_\val) \cap
\GL_\D(\OOO_\val)$.

\blemm \label{lem:compound}
For every $g\in \GL_\D(K_\val)$, there exists $\sigma\in W_\D$
such that if $\sigma g= \big(\begin{smallmatrix} \alpha &\ga \\ \beta
& \delta\end{smallmatrix}\big)$ with $\alpha\in\M_\d(K_\val)$, then
$\alpha\in\GL_\d(K_\val)$ and $\beta\,\alpha^{-1}\in\M_{\dd,\d}(\OOO_\val)$.
\elemm

\dem For every $g\in \M_\D(K_\val)$, we denote by $g_{\d\mid}$ the
submatrix of $g$ consisting of its first $\d$ columns. For all
$\alpha\in\M_\d(K_\val)$ and $j,k\in\llbracket 1,\d\rrbracket$, we denote
by $\alpha_{\wh j,\wh k}$ the submatrix of $\alpha$ where the $j$-th
row and $k$-th column have been removed. Recall that the $(j,k)$
coefficient of the comatrix $\operatorname{Comm}(\alpha)$ of $\alpha$
is $\operatorname{Comm} (\alpha)_{j,k}=(-1)^{j+k}\det \alpha_{\wh
  j,\wh k}$.

Note that the statement of Lemma \ref{lem:compound} is invariant by
multiplication on the left of $g$ by an element of $W_\D$. Since
multiplying $g$ on the left by an element of $W_\D$ amounts to
permuting the rows of $g$, up to such a multiplication, we may assume
that the absolute value of the upper-left $\d\times\d$ minor of
$g_{\d\mid}$ (hence of $g$) is maximal over the absolute values of all
$\d\times\d$ minors of $g_{\d\mid}$. Let us then prove that if $g=\big(
\begin{smallmatrix} \alpha &\ga \\ \beta & \delta\end{smallmatrix}
\big)$ with $\alpha\in \M_\d(K_\val)$, then $\alpha\in\GL_\d(K_\val)$
and $\beta\,\alpha^{-1} \in\M_{\dd,\d}(\OOO_\val)$, which proves Lemma
\ref{lem:compound} by taking $\sigma=\id$.

Since $g$ is invertible, the rank of its submatrix $g_{\d\mid}$ is
$\d$. Hence $g_{\d\mid}$ has at least one nonzero $\d\times\d$ minor,
so that $|\det \alpha|\neq 0$ by the above maximality property.

For all $i\in\llbracket 1,\dd\rrbracket$ and $j\in\llbracket
1,\d\rrbracket$, let us prove that the $(i,j)$-coefficient
$(\beta\alpha^{-1})_{i,j}$ of the matrix $\beta\alpha^{-1}\in
\M_{\dd,\d}(K_\val)$ has absolute value at most $1$, which proves
Lemma \ref{lem:compound}. We denote by $A(i,j)\in\M_{\d}(K_\val)$ the
matrix $\alpha$ where its $j$-th row has been replaced by the
$(i+\d)$-th row of $g_{\d\mid}$. By the above maximality property and
since $\det A(i,j)$ is up to a sign an $\d\times\d$ minor of
$g_{\d\mid}$, we have
\[
|\det A(i,j)|\leq |\det\alpha|\;.
\]
Since the $i$-th row of $\beta$ is the $(i+\d)$-th row of
$g_{\d\mid}$, since $\alpha^{-1}= \frac{1} {\det\alpha}
\;^t\!\operatorname{Comm}(\alpha)$, and by the Laplace expansion
formula for the determinant of $A(i,j)$ with respect to its $j$-th
row, we have
\begin{align*}
(\beta\alpha^{-1})_{ij}&=\sum_{k=1}^\d \beta_{i,k}(\alpha^{-1})_{k,j}
=\frac{1}{\det\alpha}\sum_{k=1}^\d g_{i+\d,k}\operatorname{Comm} (\alpha)_{j,k}
\\&=\frac{1}{\det\alpha}\sum_{k=1}^\d (-1)^{j+k}g_{i+\d,k}\det \alpha_{\wh
  j,\wh k}=\frac{\det A(i,j)}{\det\alpha}\;.
\end{align*}
Therefore $|(\beta\alpha^{-1})_{ij}|\leq 1$, as wanted.
\cqfd

\medskip
The linear action of an element $\sigma$ of the Weyl group $W_\D$ on
an element $\Lambda\in\P\L_{\d,\D}$ satisfies the following properties.

$\bullet$~ Since $\sigma\in\GL_\D(R_\val)$, the $\d$-lattice
$\sigma\Lambda$ is primitive, and by Equation \eqref{eq:detetcovol},
we have
\[
\overline{\covol}(\sigma\Lambda)=\overline{\covol}(\Lambda)\;.
\]

$\bullet$~ We have $V_{\sigma\Lambda}=\sigma V_{\Lambda}$ by the left
hand side of Equation \eqref{eq:transfocovol}.

$\bullet$~ Since $\sigma\in\GL_\D(\OOO_\val)$, and by the construction
of the shape map $\sh$ in and above Equation \eqref{eq:defish}, we have
$\sh(\sigma\Lambda)=\sh(\Lambda)$.

$\bullet$~ Let $R_\val^{\;\D,*}$ be the standard full $R_\val$-lattice
of the dual space of $K_\val^{\;\D}$, which is invariant under the
dual action of $\sigma$ since $\widecheck\sigma=\;^t\sigma^{\,-1}\in
\GL_\D(R_\val)$.  As seen in Equation \eqref{eq:proportholat}, we have
$(\sigma\Lambda)^\perp= \widecheck\sigma(\Lambda^\perp)$.  Hence
$\sh((\sigma\Lambda)^\perp)=\sh(\Lambda^\perp)$ since $\widecheck
\sigma \in\GL_\D(\OOO_\val)$.

$\bullet$~ By the $\GL_\D(\OOO_\val)$-invariance of the probability
measure $\mu_{\Gr_{\d,\D}}$ (see Subsection \ref{subsect:grassmann}),
and since $\sigma\in\GL_\D(\OOO_\val)$, we have
$\sigma_*\mu_{\Gr_{\d,\D}}=\mu_{\Gr_{\d,\D}}$.

$\bullet$~ Since $\sigma\in\GL_\D(\OOO_\val)$, the left action of
$\sigma$ on $\Gr_{\d,\D}$ is isometric for the distance $d$ on
$\Gr_{\d,\D}$ constructed in Subsection \ref{subsect:grassmann}.

\medskip
By Lemma \ref{lem:compound}, for every $\Lambda \in\P\L_{\d,\D}
\ssm\P\L_{\d,\D}^\sharp$ with $\overline{\covol} (\Lambda)\in
q_\val^{\;\ell\,\ZZ}$, there exists $\sigma\in W_\D$ such that $\sigma
\Lambda\in\P\L_{\d,\D}^\sharp$ and $V_{\sigma\Lambda} =\sigma
V_{\Lambda} \in \Gr_{\d,\D}^\flat$.  Furthermore, $\sigma$ maps a
small ball centered at $V_{\Lambda}$ contained in $\Gr_{\d,\D}\ssm
\Gr_{\d,\D}^\flat$ to a ball contained in $\Gr_{\d,\D}^\flat$
centered at $V_{\sigma\Lambda}$ of the same radius, by the last point
above.  Hence the equidistribution with error term as $i\ra+\infty$ of
$(V_\Lambda, \sh(\Lambda), \sh(\Lambda^\perp))$ in $\Gr_{\d,\D} \times
\Sh_\d \times \Sh_\dd$ when $\Lambda$ varies in $\PL_{\d,\D}(I)$ with
$\overline{\covol}(\Lambda)=q_\val^{\;\ell\,i}$ follows from the
equidistribution with error term as $i\ra+\infty$ of $(V_\Lambda,
\sh(\Lambda),\sh(\Lambda^\perp))$ in $\Gr_{\d,\D}^\flat \times \Sh_\d
\times \Sh_\dd$ when $\Lambda$ varies in $\PL_{\d,\D}^\sharp(I)$ with
$\overline{\covol}(\Lambda)= q_\val^{\;\ell\,i}$.
\cqfd

\bigskip
We conclude this paper with a proof of Corollary \ref{coro:intro2} in
the introduction.
  
\medskip
\noindent {\bf Proof that Theorem \ref{theo:mainintrotriple} implies
  Corollary \ref{coro:intro2}. }  For $k\in\{\d,\dd\}$, with $\wt\Ga_k=
\GL_k(R_\val)$, let us consider the map $\iota :D\mapsto D^{-1}$
defined in the Introduction from the discrete set $\Sh_{k}=
\GL_k(\OOO_\val) \,\bs \GL_k^1(K_\val) \,/\,\wt\Ga_k$ to the discrete
set $\wt\Ga_k\bs V_0\I_{\val,k}= \wt\Ga_k \,\bs \GL_k^1(K_\val)
\,/\,\GL_k(\OOO_\val)$. By Equation \eqref{eq:normalmuGL1} , this map
satisfies
\[
\iota_*\mu_{\Sh_{k}}= \mu_{\GL_k^1(K_\val)}(\GL_k(\OOO_\val))\,
\mu_{\wt\Ga_k\bs V_0\I_{\val,k}}= \mu_{\wt\Ga_k\bs V_0\I_{\val,k}}\;.
\]
Let $\varphi':\Gr_{\d,\D}\times \Sh_{\d}\times \Sh_{\dd}\ra
\wt\Ga_\d\bs V_0 \I_{\val,\d}\times \wt\Ga_\dd\bs V_0 \I_{\val,\dd}$
be the continuous map defined by $(x,y,z)\mapsto (y^{-1},z^{-1})$.
Since $\Gr_{\d,\D}$ is compact, this map $\varphi'$ is proper, and the
pushforward map $\varphi'_*$ of Borel measures by $\varphi'$ is linear
and weak-star continuous. With $\ell=\lcm(\d,\dd)$, the image by
$\varphi'_*$ of the left hand side of Equation
\eqref{eq:mainintrotriple} is hence
\[
\lim_{i\ra+\infty}\;\frac{c'}{q_\val^{\;\ell\,\D\,i}}
\sum_{\Lambda\in \PL_{\d,\D}\;:\;\overline{\covol}(\Lambda)=
  q_\val^{\;\ell\, i}}\Delta_{\sh(\Lambda)^{-1}}
\otimes\Delta_{\sh(\Lambda^\perp)^{-1}}\;.
\]
Since $\|\mu_{\Gr_{\d,\D}}\|=1$ by Equation \eqref{eq:normalmesgrass},
the image by $\varphi'_*$ of the right hand side of Equation
\eqref{eq:mainintrotriple} is $ \mu_{\wt\Ga_\d\bs V_0 \I_{\val,\d}}
\otimes \mu_{\wt\Ga_\dd\bs V_0 \I_{\val,\dd}}$. Hence Theorem
\ref{theo:mainintrotriple} does imply Corollary \ref{coro:intro2}.
\cqfd

\medskip\rem Proceeding as in the proof of Corollary
\ref{coro:mainshape}, for every nonzero ideal $I$ of $R_\val$, we have
the following version with error term and congruences of Corollary
\ref{coro:intro2}~: There exists $\tau\in\;]0,\frac{1}{2\,\D^2}]$ such
that for all finitely supported maps $f_1:\wt\Ga_\d\bs V_0\I_{\val,\d}
\ra\RR$ and $f_2:\wt\Ga_\dd\bs V_0\I_{\val,\dd}\ra\RR\;$ and for every
$\delta\in\; ]0,\tau[$, we have
\begin{align}
  &\frac{c'}{q_\val^{\;\ell\,\D\,i}}
\sum_{\Lambda\in\PL_{\d,\D}(I)\;:\;\overline{\covol}\;\Lambda\,=\,q_\val^{\;\ell\,i}} 
f_1(\sh(\Lambda)^{-1})\;f_2(\sh(\Lambda^\perp)^{-1})\nonumber\\ =\;& 
\Big(\int f_1\,d\mu_{\wt\Ga_\d\bs V_0 \I_{\val,\d}}\Big)
\Big(\int f_2\,d\mu_{\wt\Ga_\dd\bs V_0 \I_{\val,\dd}}\Big) +\bigO_{\val,\delta,I}
\big(q_\val^{\;\ell\,\D\,i(-\tau+\delta)}
\,\|f_1\|_\infty\|f_2\|_\infty\big)\;.\label{eq:equidisbuilding}
\end{align}

\appendix
\section{Dual and factor partial lattices}
\label{appen:dualfactorlat}

Let $V$ be a $K_\val$-vector space with finite dimension $D\in\NN$
endowed with an ultrametric norm $\|\;\|$, and let $k\in\llbracket
1,D-1\rrbracket$ if $D\geq 2$.

Let $W$ be a $k$-dimensional $K_\val$-vector subspace of $V$, endowed
with the restriction norm. We endow the quotient $(D-k)$-dimensional
$K_\val$-vector space $V/W$ with the quotient norm. We denote by $\pi:
V\ra V/W$ the canonical projection. Then the Haar measures $\mu_V$,
$\mu_W$ and $\mu_{V/W}$ on respectively $V$, $W$, $V/W$, normalized by
these choices of norms as explained in Subsection
\ref{subsect:partiallatt}, satisfy the following Weil's normalization
process (see \cite[\S 9]{Weil65}). For all $\ov{x}\in V/W$ and $x\in
V$ such that $\pi(x)=\ov x$, let $\mu_{\pi^{-1}(\ov{x})}$ be the
measure on the $K_\val$-affine subspace $\pi^{-1}(\ov{x})=x+W$ such
that the translation $\tau_x:y\mapsto x+y$, which is a homeomorphism
from $W$ to $\pi^{-1}(\ov{x})$, satisfies $(\tau_x)_*\mu_W=
\mu_{\pi^{-1}(\ov{x})}$ (this does not depend on the choice of $x$ in
$\pi^{-1}(\ov{x})$). Then {\it Weil's normalisation} is asking that we
have the following disintegration property of the measure $\mu_V$ over
the measure $\mu_{V/W}$ by $\pi$:
\begin{equation}\label{eq:weildisinteg}
d\mu_V=\int_{\ov{x}\in V/W}d\mu_{\pi^{-1}(\ov{x})}\;d\mu_{V/W}(\ov{x})\;.
\end{equation}
This formula implies that the normalizations of $\mu_V$ and $\mu_W$
uniquely determine the normalization of $\mu_{V/W}$. In order to check
that this normalization coincides with the one coming from the
quotient norm on $V/W$, we apply the above formula on $B_{V}(0,1)$,
noting that $\pi(B_V(0,1))=B_{V/W}(0,1)$ and that, by the ultrametric
property, for every $x\in B_{V}(0,1)$, we have $-x+B_{V}(0,1)\cap
(x+W) = B_{W}(0,1)$.

\medskip
The ultrametric norm $\|\,\|$ on $V$ is {\it integral} if its set of
values is exactly $\{0\}\cup q_\val^\ZZ$ unless $D=0$. For instance,
the usual norm on $K_\val^{\;D}$ is integral. If $\|\,\|$ is integral,
then the restriction of $\|\,\|$ to any $K_\val$-linear subspace of
$V$ is integral, the quotient norm of $\|\,\|$ on any quotient $V/W'$
of $V$ by a $K_\val$-linear subspace $W'$ of $V$ is integral, and the
dual norm $\|\;\|^*:\ell\mapsto\max_{x\in V\smallsetminus\{0\}}
\frac{|\ell(x)|} {\|x\|}$ on $V^*$ of $\|\;\|$ is integral. An {\it
  orthogonal} $K_\val$-basis of $V$ is a $K_\val$-basis $(e_1,\dots,
e_D)$ of $V$ such that we have $\|\sum_{i=1}^D\lambda_i b_i\| =
\max_{1\leq i\leq D} |\lambda_i|\,\|b_i\|$ for all $\lambda_1,\dots,
\lambda_k\in K_\val$. It is {\it orthonormal} if furthermore $\|b_i\|
=1$ for every $i\in\llbracket 1,D\rrbracket$.  An orthogonal
$K_\val$-basis of $V$ exists by for instance
\cite[Prop.~1.1]{GolIwa63}.  If the norm $\|\,\|$ is integral, then
there exists an orthonormal $K_\val$-basis of $V$, by renormalizing an
orthogonal basis.

\medskip
Let $\Lambda$ be a $k$-lattice in $V$. Its {\it dual lattice} is the
$R_\val$-submodule of the dual $K_\val$-vector space $V_\Lambda^{\;*}$
defined by
\[
\Lambda^*=\{\ell\in V_\Lambda^{\;*}:\forall\;x\in \Lambda,\;\;\ell(x)\in
R_\val\}\;.
\]

\blemm\label{lem:propriduallat} The dual lattice $\Lambda^*$ is a full
$R_\val$-lattice in $V_\Lambda^{\;*}$, and we have $(\Lambda^*)^*=
\Lambda$.  If the norm of $V$ is integral, then
\[
\covol(\Lambda^*) \covol(\Lambda)=q^{2k(\ggg-1)}\;.
\]
%\begin{equation}\label{eq:covoldual}
%  \covol(\Lambda)\covol(\Lambda^*)=q^{2k(\ggg-1)}\;.
%\end{equation}
\elemm

\dem Let $(b_1,\dots, b_k)$ be an $R_\val$-basis of $\Lambda$. Then
$(b_1,\dots, b_k)$ is a $K_\val$-basis in $V_\Lambda$, and we denote
by $(b_1^*,\dots,b_k^*)$ its dual $K_\val$-basis in $V_\Lambda^{\;*}$.
A linear form $\ell=\sum_{i=1}^k \lambda_i b_1^*\in V_\Lambda^{\;*}$
takes integral values on all elements of $\Lambda$ if and only if it
takes integral values on $b_1,\dots, b_k$, that is, if and only if its
coordinates $\lambda_1,\dots,\lambda_k$ are integral. Thus $\Lambda^*=
\oplus_{i=1}^k R_\val b_i^*$, which is a full $R_\val$-lattice in the
$k$-dimensional vector space $V_\Lambda^{\;*}$. Since the dual basis
in $(V_\Lambda^{\;*})^*=V_\Lambda$ of the $K_\val$-basis $(b_1^*,
\dots, b_k^*)$ of $V_\Lambda^{\;*}$ is the $K_\val$-basis $(b_1,\dots,
b_k)$ of $V_\Lambda$, we have $(\Lambda^*)^*=\Lambda$.

Assume now that the ultrametric norm $\|\,\|$ on $V$ is integral. Let
$(e_1,\dots, e_k)$ be an orthonormal $K_\val$-basis of $V_\Lambda$.
Let $P\in\GL_k(K_\val)$ be the transition matrix from the
$K_\val$-basis $(e_1,\dots, e_k)$ to the $K_\val$-basis $(b_1,\dots,
b_k)$ of $V_\Lambda$. Then the dual $K_\val$-basis $(e_1^*,\dots,
e_k^*)$ of $V_\Lambda^{\;*}$ is an orthonormal $K_\val$-basis of
$V_\Lambda^{\;*}$ for the dual norm $\|\;\|^*:\ell\mapsto \max_{x\in
  V_\Lambda \smallsetminus \{0\}} \frac{|\ell(x)|}{\|x\|}$ on
$V_\Lambda^{\;*}$. The transition matrix from the $K_\val$-basis
$(e^*_1,\dots, e^*_k)$ to the $K_\val$-basis $(b^*_1, \dots, b^*_k)$
of $V_\Lambda^{\;*}$ is $^tP^{-1}$. Let $\Lambda_e=\oplus_{i=1}^k
R_\val e_i$ and $\Lambda^*_e =\oplus_{i=1}^kR_\val e_i^*$. We have
$\overline{\covol} (\Lambda_e)= \overline{\covol} (\Lambda^*_e)=1$
since the coordinate maps from $V_\Lambda$ with $K_\val$-basis
$(e_1,\dots, e_k)$ to $K_\val^{\;k}$ and from $V_\Lambda^{\;*}$ with
$K_\val$-basis $(e_1^*,\dots, e_k^*)$ to $K_\val^{\;k}$ are isometries
sending $\Lambda_e$ and $\Lambda^*_e$ to $R_\val^{\;k}$, and by
Equation \eqref{eq:covolRv}. Hence using twice Equation
\eqref{eq:detetcovol}, we have

\[
\overline{\covol}(\Lambda^*) \;\overline{\covol}(\Lambda)
=(\,|\det \,^tP^{-1}|\;\overline{\covol}(\Lambda^*_e)\,)
(\,|\det P\,|\;\overline{\covol}(\Lambda_e)\,)=1\,,
\]
as wanted by Equation \eqref{eq:covolRv}.
\cqfd

\medskip
Assume that $V$ is also endowed with an integral structure
$V_{R_\val}$.  Let $\Lambda$ be a primitive $k$-lattice in the
integral $K_\val$-space $V$. The {\it factor lattice} of $\Lambda$ is
the $R_\val$-submodule $\Lambda^\pi$ of the quotient $K_\val$-vector
space $V/V_\Lambda$ which is the image of $V_{R_\val}$ by the
canonical projection $\pi:V\ra V/V_\Lambda$.

\blemm\label{lem:proprifaclat} The factor lattice $\Lambda^\pi$ is a
full $R_\val$-lattice in $V/V_\Lambda$. The canonical $K_\val$-linear
isomorphism $V/V_\Lambda\ra (V_\Lambda^\perp)^*$ maps $\Lambda^\pi$ to
$(\Lambda^\perp)^*$. We have
\[
\covol(\Lambda^\pi)=\covol((\Lambda^\perp)^*) \quad\text{and}\quad 
\covol(\Lambda^\pi) \covol(\Lambda)=\covol(V_{R_\val})
\]
%\begin{equation}\label{eq:covolfactor}
%  \covol(\Lambda)\covol(\Lambda^\pi)=\covol(V_{R_\val})\;.
%\end{equation}
%
%\begin{equation}\label{eq:covolperfacdual}
%  \covol(\Lambda^\pi)=\covol((\Lambda^\perp)^*)=?????????\;.
%\end{equation}
\elemm

\dem Since the $k$-lattice $\Lambda$ is primitive, there exists an
$R_\val$-basis $(b_1,\dots, b_D)$ of $V_{R_\val}$ such that
$(b_1,\dots, b_k)$ is an $R_\val$-basis of $\Lambda$, hence a
$K_\val$-basis of $V_\Lambda$. Then $(\pi(b_{k+1}),\dots, \pi(b_D))$
is a $K_\val$-basis of $V/V_\Lambda$, and an $R_\val$-basis of
$\Lambda^\pi$ by definition. Hence $\Lambda^\pi$ is a full
$R_\val$-lattice in the $(D-k)$-dimensional vector space
$V/V_\Lambda$.

Identifying a $K_\val$-vector space $W$ with its bidual $(W^*)^*$ by
the map $x\mapsto (\ell\mapsto \ell(x))$ as usual, the map $\wt
{\Theta''}:V\ra (V_\Lambda^\perp)^*$ defined by $x\mapsto x_{\mid
  V_\Lambda^\perp}$ induces a $K_\val$-linear isomorphism
$\Theta'':V/V_\Lambda\ra (V_\Lambda^\perp)^*$.  With $(b_1,\dots,
b_D)$ as above, we have seen in the proof of Proposition
\ref{prop:propriortholat} that $(b_{k+1}^*, \dots, b_D^*)$ is an
$R_\val$-basis of $\Lambda^\perp$, hence a $K_\val$-basis of
$V_{\Lambda^\perp} =V_{\Lambda}^\perp$. As seen in the proof of Lemma
\ref{lem:propriduallat}, the dual $K_\val$-basis $({b_{k+1}^*}^*,
\dots, {b_D^*}^*)$ of $(b_{k+1}^*, \dots, b_D^*)$ is an $R_\val$-basis
of $(\Lambda^\perp)^*$.  But $({b_{k+1}^*}^*, \dots, {b_D^*}^*)$ is
exactly $\big(\Theta''(\pi(b_{k+1})),\dots, \Theta''(\pi(b_{k+1}))
\big)$.  Hence $\Theta''(\Lambda^\pi)= (\Lambda^\perp)^*$.

When $V/V_\Lambda$ is endowed with the quotient norm, and
$(V_\Lambda^\perp)^*$ with the dual norm of the restriction to
$V_\Lambda^\perp$ of the dual norm on $V^*$, the above map $\Theta''$ is
an isometry. Hence $\covol(\Lambda^\pi)=\covol((\Lambda^\perp)^*)$.

Let $F$ be a clopen strict fundamental domain for the action of
$R_\val$ on $K_\val$. The formula $\covol(V_{R_\val})=\covol(\Lambda)
\covol(\Lambda^\pi)$ follows by integrating Equation
\eqref{eq:weildisinteg} with $W=V_\Lambda$ on the strict fundamental
domain $F\,b_1+\dots+ F\,b_D$ of $V_{R_\val}$.
\cqfd

\medskip\noindent{\bf Proof of Equation \eqref{covollatortholat}. }
Under its assumptions, we may assume that $V={K_\val}^D$ and
$V_{R_\val} = {R_\val}^D$ and that the norm of $V$ is the standard
supremum norm of ${K_\val}^D$, which is integral. We then have $\covol
(V_{R_\val})= q^{D(\ggg-1)}$ by Equation \eqref{eq:covolRv}. Hence
respectively by Lemma \ref{lem:propriduallat} (recalling that
$\Lambda^\perp$ is a $(D-k)$-lattice in $V^*$), by the first equality
in Lemma \ref{lem:proprifaclat}, and by the second equality in Lemma
\ref{lem:proprifaclat}, we have
\begin{align*}
\covol(\Lambda^\perp)&=q^{2(D-k)(\ggg-1)}\covol((\Lambda^\perp)^*)^{-1}
=q^{2(D-k)(\ggg-1)}\covol(\Lambda^\pi)^{-1}\\&=
q^{2(D-k)(\ggg-1)}\covol(\Lambda)\;q^{-D(\ggg-1)}
=\covol(\Lambda)\;q^{(D-2k)(\ggg-1)}\;,
\end{align*}
as wanted by Equation \eqref{eq:covolRv}. \cqfd

{\small \bibliography{../biblio} }

\bigskip
{\small
\noindent \begin{tabular}{l} 
Department of Mathematics, Rämistrasse 101\\
ETH Zurich, 8092 Zurich, Switzerland\\
{\it e-mail: tal.horesh@math.ethz.ch}
\end{tabular}
\medskip

\noindent \begin{tabular}{l}
Laboratoire de math\'ematique d'Orsay, UMR 8628 CNRS\\
Universit\'e Paris-Saclay,
91405 ORSAY Cedex, FRANCE\\
{\it e-mail: frederic.paulin@universite-paris-saclay.fr}
\end{tabular}
}

\end{document}

%% file: macros_angl_pdf.tex
\usepackage[utf8]{inputenc}
\usepackage{amsmath,amssymb,epsfig,bbm}
\usepackage{stmaryrd,mathabx,frcursive}
\usepackage{comment}
\usepackage{color}
\usepackage[T1]{fontenc}

\usepackage[textsize=small]{todonotes}
\usepackage{enumitem}
\usepackage{varwidth}
\setlist{nolistsep}
\usepackage{hyperref}
\usepackage{bm,mathtools}

\newtheorem{defi}{Definition}[section]
\newtheorem{prop}[defi]{Proposition}
\newtheorem{theo}[defi]{Theorem}
\newtheorem{conj}[defi]{Conjecture}
\newtheorem{lemm}[defi]{Lemma}
\newtheorem{coro}[defi]{Corollary}
\newtheorem{rema}[defi]{Remark}
\newtheorem{exem}[defi]{Example}
\newtheorem{exems}[defi]{Examples}

\newcommand{\bdefi}{\begin{defi}}
\newcommand{\edefi}{\end{defi}}
\newcommand{\bprop}{\begin{prop}}
\newcommand{\eprop}{\end{prop}}
\newcommand{\btheo}{\begin{theo}}
\newcommand{\etheo}{\end{theo}}
\newcommand{\blemm}{\begin{lemm}}
\newcommand{\brema}{\begin{rema}}
\newcommand{\erema}{\end{rema}}
\newcommand{\bexer}{\begin{exem}}
\newcommand{\eexer}{\end{exem}}
\newcommand{\bexems}{\begin{exems}}
\newcommand{\eexems}{\end{exems}}
\newcommand{\bconj}{\begin{conj}}
\newcommand{\econj}{\end{conj}}
\newcommand{\elemm}{\end{lemm}}
\newcommand{\bcoro}{\begin{coro}}
\newcommand{\ecoro}{\end{coro}}
\newcommand{\dem}{\noindent{\bf Proof. }}
\newcommand{\rem}{\noindent{\bf Remark. }}

\usepackage{mathrsfs}
\renewcommand\mathcal{\mathscr}

\newcommand{\B}{{\cal B}}

\newcommand{\D}{{\cal D}}
\newcommand{\E}{{\cal E}}
\newcommand{\F}{{\cal F}}

\newcommand{\I}{{\cal I}}

\renewcommand{\L}{{\cal L}}
\newcommand{\M}{{\cal M}}

\newcommand{\OOO}{{\cal O}}
\renewcommand{\P}{{\cal P}}

\newcommand{\Scal}{{\cal S}}

\newcommand{\U}{{\cal U}}
\newcommand{\V}{{\cal V}}

\newcommand{\Z}{{\cal Z}}

\newcommand{\maths}[1]{{\mathbb #1}}  

\newcommand{\CC}{\maths{C}}

\newcommand{\FF}{\maths{F}}

\newcommand{\NN}{\maths{N}}

\newcommand{\PP}{\maths{P}}
\newcommand{\QQ}{\maths{Q}}
\newcommand{\RR}{\maths{R}}

\newcommand{\ZZ}{\maths{Z}}

\renewcommand{\ggg}{{\mathfrak g}}

\newcommand{\ppp}{{\mathfrak p}}

\newcommand{\uuu}{{\mathfrak u}}

\newcommand{\ra}{\rightarrow}
\newcommand{\bs}{\backslash}

\newcommand{\ov}[1]{{\overline #1}} 
\newcommand{\wt}[1]{{\widetilde{#1}}}
\newcommand{\wh}[1]{{\widehat{#1}}}

\newcommand{\ga}{\gamma}
\newcommand{\Ga}{\Gamma}

\newcommand{\cqfd}{\hfill$\Box$}

\newcommand{\bigO}{\operatorname{O}}
\newcommand{\card}{{\operatorname{Card}}}

\newcommand{\covol}{\operatorname{Covol}}

\newcommand{\dbs}{\backslash\!\!\backslash}

\newcommand{\Gr}{\operatorname{Gr}}

\newcommand{\id}{\operatorname{id}}

\newcommand{\La}{\operatorname{Lat}^1}

\newcommand{\lcm}{\operatorname{lcm}}
\renewcommand{\ln}{\operatorname{ln}}

\newcommand\PL{\operatorname{\P\!\L}}

\renewcommand{\Re}{{\operatorname{Re}}}

\newcommand{\ssm}{{\smallsetminus}}

\newcommand{\sh}{\operatorname{sh}}
\newcommand{\Sh}{\operatorname{Sh}^1}

\newcommand{\SL}{\operatorname{SL}}
\newcommand{\GL}{\operatorname{GL}}

\newcommand{\PGL}{\operatorname{PGL}}

\usepackage{ulem}